\documentclass[10pt]{article}
\usepackage[utf8]{inputenc}
\usepackage{amsmath,amsthm,amsfonts}
\usepackage{graphicx}
\usepackage{color}
\usepackage{multirow}
\usepackage{todonotes}

\newtheorem{theorem}{Theorem}

\usepackage{subcaption}
\graphicspath{ {./figs/}}

\newcommand{\rev}[1]{{\color{black}#1}}

\begin{document}

\title{Multiscale method for image denoising using nonlinear diffusion process: local denoising and spectral multiscale basis functions}

\author{
Maria Vasilyeva
\thanks{Department of Mathematics \& Statistics, Texas A\&M University - Corpus Christi, Corpus Christi, TX, USA.  
Email: {\tt maria.vasilyeva@tamucc.edu}.} 
\and
Aleksei Krasnikov
\thanks{Department of Mathematics \& Statistics, Texas A\&M University - Corpus Christi, Corpus Christi, TX, USA.  
Email: {\tt akrasnikov@islander.tamucc.edu}.} 
\and
Kelum Gajamannage
\thanks{Department of Mathematics and Applied Mathematical Sciences, University of Rhode Island, Kingston, RI, USA.  
Email: {\tt kelum.gajamannage@uri.edu}.} 
\and
Mehrube Mehrubeoglu 
\thanks{Department of Engineering, Texas A\&M University - Corpus Christi, Corpus Christi, TX, USA.  
Email: {\tt ruby.mehrubeoglu@tamucc.edu}.} 
}

\maketitle

\begin{abstract}
We consider image denoising using a nonlinear diffusion process, where we solve unsteady partial differential equations with nonlinear coefficients. The noised image is given as an initial condition, and nonlinear coefficients are used to preserve the main image features. In this paper, we present a multiscale method for the resulting nonlinear parabolic equation in order to construct an efficient solver. To both filter out noise and preserve essential image features during the denoising process, we utilize a time-dependent nonlinear diffusion model. Here, the noised image is fed as an initial condition and the denoised image is stimulated with given parameters. We numerically implement this model by constructing a discrete system for a given image resolution using a finite volume method and employing an implicit time approximation scheme to avoid time-step restriction. However, the resulting discrete system size is proportional to the number of pixels which leads to computationally expensive numerical algorithms for high-resolution images. In order to reduce the size of the system and construct efficient computational algorithms, we construct a coarse-resolution representation of the system. We incorporate local noise reduction in the coarsening process to construct an efficient algorithm with fewer denoising iterations. We propose a computational approach with two main ingredients: (1) performing local image denoising in each local domain of basis support; and (2) constructing multiscale basis functions to construct a coarse resolution representation by a Galerkin coupling. We present numerical results for several \rev{classic and high-resolution image datasets} to demonstrate the effectiveness of the proposed multiscale approach with local denoising and local \rev{multiscale} representation.
\end{abstract}

\section{Introduction}

Image denoising is a fundamental problem in image processing and computer vision in which the noise is filtered out from a noisy image to recover the noise-free version of it \cite{gajamannage2024image}. Various image-denoising techniques have been developed \cite{fan2019brief, rudin1992nonlinear, gajamannage2022efficient, gajamannage2023geodesic}. Among them, the Partial Differential Equation (PDE)-based approaches have gained significant attention due to their theoretical foundation and effectiveness in preserving important image features including edges while reducing noise. PDE-based methods formulate image denoising as a problem of solving a linear or nonlinear time-dependent equation within a given image resolution that describes the evolution of pixel intensities \cite{alvarez1993axioms, alvarez1992image, aubert2006mathematical}. These methods are grounded by a nonlinear diffusion process, which smooths the image while preserving essential features like edges. The Perona-Malik equation describes the classic nonlinear diffusion process, which introduces a nonlinear coefficient for edge detection based on the image gradient and addresses the edge-blurring issue occurring in the linear diffusion process \cite{perona1994anisotropic, perona1990scale}. Various modifications generalize the Perona-Malik model by allowing the diffusion tensor to vary spatially and anisotropically to achieve high-quality denoising results \cite{wei1999generalized, yuan2016perona, wang2018hybrid, ahmed2008segmentation}. 

The use of PDEs for image denoising provides a solid framework for reducing noise in images while preserving important features. When dealing with high-resolution images, specialized techniques are required to address the increased computational complexity. The underlying nonlinear time-dependent problems are defined on the grid that is related to the initial noisy image resolution, then the number of unknowns is directly proportional to the number of pixels \cite{sturmer2008fast}. The larger number of pixels increases computational complexity and necessitates the development of special numerical techniques  \cite{acton1998multigrid, brito2010multigrid, sturmer2008fast, zhang2013nonlinear}. To reduce computational complexity of solving diffusion-based equation, various techniques were developed that include preconditioned iterative techniques to solve resulting large system of linear equations instead of usage of direct solvers. The other class of techniques includes the construction of the accurate coarse-scale representation of the considered problem and is related to the homogenization and multiscale techniques. 
Multiscale methods and homogenization techniques are used to construct efficient and accurate macroscale model representations that incorporate underlying fine-scale heterogeneity into a coarse-scale model.
In the image denoising problem, the heterogeneity is induced by the gradient of the image and should be accurately addressed in the coarse-scale model construction process.
Homogenization methods in material science involve deriving effective macroscopic material properties by averaging microscopic variations \cite{allaire1992homogenization, bakhvalov2012homogenisation}. The system is then solved on the macroscale level while preserving the overall behavior of the material. Multiscale methods extend the concept of homogenization by explicitly considering multiple scales in the mathematical model. Most multiscale methods are based on constructing accurate multiscale space for a coarse-scale approximation via generating multiscale basis functions \cite{efendiev2009multiscale}. For example, in the Multiscale Finite Element Method (MsFEM) \cite{hou1997multiscale, efendiev2009multiscale}, the local problems are solved in each local domain of the basis functions support. Standard linear basis functions are used as boundary conditions for a local problem \cite{allaire2005multiscale}. In the Generalized Multiscale Finite Element Method (GMsFEM)  \cite{efendiev2011multiscale, efendiev2013generalized}, the MsFEM concept is extended to high-contrast problems by deriving multiple multiscale basis functions in each local domain by solving local spectral problem. Both homogenization and multiscale methods provide powerful frameworks for analyzing complex systems with multiple scales. These techniques allow for the derivation of simplified models that capture the system's essential features while accounting for the effects of fine-scale variations.

In this work, we construct a coarse scale approximation for the nonlinear time-dependent equation using GMsFEM. In our previous works, the GMsFEM was applied to a wide range of applications  (Li-ion batteries, geothermal reservoirs, unconventional reservoirs, composite materials, seismic wave propagation), including different types of heterogeneity (perforated, high-contrast, fractured, networks)  \cite{chung2018multiscale, vasilyeva2018multiscale, akkutlu2018multiscale, vasilyeva2019upscaling, vasilyeva2021multiscale, vasilyeva2024decoupled, vasilyeva2024generalized}. In GMsFEM, we construct a multiscale basis function to capture the behavior of the solution at a fine scale.
The GMsFEM approach involves two stages: offline and online. In the offline stage, we define local domains (subdomains), construct multiscale basis functions by solving the local eigenvalue problems, and generate a projection matrix. In the online stage, we project fine-resolution representation onto a fine grid using a projection matrix and solve a reduced-order problem on the coarse grid. The solution obtained using GMsFEM accurately represents the multiscale behavior on the coarse-scale grid by introducing a spectral multiscale basis function. 
The accuracy of the approximation properties of the constructed multiscale space highly depends on the nonlinear coefficient that depends on an initial given noised image. To address this issue, we propose an additional local denoising process for a local image that can significantly improve the basis representation and capture the 'right' behavior related to the global denoising iterations. Accurate multiscale basis functions are then based on the denoised local image and coupled Galerkin approach on a coarse resolution. This provides good denoising results with faster solutions on the coarse grid and better denoising results with fewer iterations. We present the construction of the spectral basis functions to illustrate the influence of noised images on local spectral behavior and show the local denoising effect on the resulting multiscale basis functions. The construction is given for a greyscale image representation and extended to the color images. Numerical results are given for several \rev{image dataset} with different noise levels. 

The paper is organized as follows:  In Section 2, we present a problem formulation which includes introducing the nonlinear parabolic equation used in the image denoising process and giving an approximation by space and time on a given image resolution.  In Section 3, we present a multiscale method, describe the main steps used to construct an accurate low-resolution representation using a local denoising process, and introduce local spectral multiscale basis functions. A numerical investigation is presented in Section 4 for several test images with different levels of noise.  Finally, the conclusion is presented in Section 5.

\section{Overview}

\rev{
Let $I_0(x)$ represent observed noised image with pixel values $x \in \Omega  \subset \mathbb{R}^d$. Did we define $\Omega$ ($d=2$ be grey-scale image and $d=3$ for color image),  $I(x)$ represent desired (unknown) denoised image and
\[
I_0(x) = I(x) + N(x),
\]
where $N(x)$ represent spatially distributed additive Gaussian noise with standard deviation $\sigma$. In general, noise can be represented within other types, e.g., salt-and-pepper, speckle noise, etc. However, we concentrate on Gaussian noise in this study. Next, we start with a quick overview of classical and state-of-the-art image denoising methods.
}

\subsection{\rev{Denoising methods}}

\rev{
Image denoising processes have been extensively represented in the literature over past years and generally can be roughly classified as follows \cite{diwakar2018review, goyal2020image, buades2005review, fan2019brief, tian2020deep}:
\begin{itemize}
\item Spatial domain methods.

In general, classical methods can be divided into two categories: (i) spatial domain filtering and (ii) variational denoising methods. 

(i) Spatial domain filtering. Linear filters, such as mean filtering and Weiner filtering, fail to preserve image textures and can over-smooth image \cite{gonzalez2009digital, jain1989fundamentals}. Nonlinear filtering, such as median filtering, weighted median filtering, and bilateral filtering, can better preserve edges \cite{pitas2013nonlinear, yin1996weighted}. 

(ii) Variational denoising is designed to minimize the energy function needed to calculate the denoised image. In the case of Gaussian noise, the objective function contains total variation and regularization terms. 
In total variational (TV) methods, $L_2$ or $ L_1$-based regularization can be used, where $L_1$-norm regularization is used for sparse representation.  
It has significantly succeeded due to effectively calculating the optimal solution and retaining sharp edges. 
A TV-based approach leads to a partial differential equation (PDE) solution and can be formulated as a nonlinear diffusion-reaction equation with an anisotropic diffusion coefficient. Furthermore, the iterative solution of the anisotropic diffusion model can be considered as a time-dependent problem, where time-marching schemes converge to the steady state solution or denoised image, and the noised image is given as an initial condition. 
Moreover, for a more accurate approximation of the small-scale features, a class of non-local approaches is introduced by non-local spatial connections; for example, each pixel value depends on the regions centered at the estimated pixel with some given radius and can be used as a regularization method \cite{buades2005non, kheradmand2014general, sutour2014adaptive}. 

\item Transform domain methods.

Transform domain methods are classified based on the chosen basis transform functions and can be divided into data-adaptive and non-data-adaptive classes. 
Data adaptive transform methods transform the given noisy image to another domain, for example, using the wavelet-based approach and then applying a denoising procedure based on the different characteristics of image and noise in the transform domain \cite{muresan2003adaptive}. 
Non-data transform methods can be further divided into spatial-frequency domain and wavelet domain. The spatial-frequency domain uses low-pass filtering based on the information that the main image is represented in a low-frequency domain, while noise belongs to the high-frequency domain. The wavelet transform leverages the sparsity of natural images in the wavelet domain but heavily depends on selecting wavelet bases \cite{mallat1989theory, malfait1997wavelet}. It contains three main steps: wavelet transform, thresholding, and inverse wavelet transform. 
BM3D (Block-Matching and 3D Filtering) is a state-of-the-art algorithm for image denoising \cite{dabov2007image}. It uses the similarity of image patches to suppress noise effectively while preserving image details and edges. BM3D contains three main steps: block matching, collaborative filtering, and two-step denoising, which include hard thresholding and Weiner filtering. 

\item Deep learning-based methods. 

This class of methods belongs to the extension of the class of spatial methods and is typically represented by the application of convolutional neural networks (CNN). In contrast with TV-based methods, which can related to the model-based optimization methods, CNN-based denoising is designed to learn a mapping function by optimizing a loss function on a training set that constrains pairs of noised and denoised images \cite{kim2016accurate, nah2017deep, jain2008natural, vincent2008extracting, xie2012image, fan2019brief, tian2020deep}. Moreover, various deep learning-based methods are developing rapidly nowadays. The main disadvantage of machine learning is the lack of image pairs for training. 

\end{itemize}

A deep learning approach requires an extensive set of data to train to learn an image denoise process in the general case, which is not only complicated by the construction of the training dataset but also requires extensive computational resources for training a high-resolution data-based on the CNN and usually done using multiple GPUs.  We will not consider a DNN model in this work and leave it for future investigation. 
The transform domain methods are classic algorithms that will be used in the numerical results section to compare with the presented approach. 
}

\subsection{\rev{Variational denoising approach}}

\rev{
In this work, we consider a class of spatial model-based methods that is related to total variational optimization and introduced by L. Rudin, S. Osher, and E. Fatemi (ROF model) \cite{rudin1992nonlinear, rudin1994total}. 
The ROF model is formulated as a constrained minimization problem
\[
\min_I J(I), \quad J(I) = \int_{\Omega} ||\nabla I(x)|| \ dx,
\]
subject to constraints
\[
\int_{\Omega} I(x) \ dx = \int_{\Omega} I_0(x) dx, \quad 
\int_{\Omega} \frac{1}{2} (I(x) - I_0(x))^2 dx = \sigma^2,
\]
where $J(I) = \int_{\Omega} ||\nabla I(x)|| \ dx$ is called the total variation of $I(x)$ over the domain $\Omega$  that promotes smoothness while preserving edges and for smooth data is equivalent to the integral of gradient magnitude. Here, constraints signify that noise $N(x)$ has a zero mean with standard deviation $\sigma$. 
This equivalent approach is to use a Lagrange multiplier $\mu$ to solve the unconstrained minimization problem
\[
\min_I  \int_{\Omega} \left( ||\nabla I(x)|| + \frac{\mu}{2}  (I(x) - I_0(x))^2 \right)  dx.
\]
This gives us the nonlinear elliptic equation (Euler-Lagrange equation)
\[
- \nabla \cdot \left(  \frac{1}{||\nabla I(x)||} \nabla I(x) \right) + \mu(I(x) - I_0(x)) = 0, \quad x \in \Omega, 
\]
with  free boundary conditions $\nabla I(x) \cdot n = 0$ on $\partial \Omega$.

The solution procedure can be considered as a solution of the parabolic equation for $I = I(x, t)$
\begin{equation}
\label{tv1}
I_t - \nabla \cdot \left( \frac{1}{||\nabla I(x)||}\nabla I \right) + \mu(I(x) - I_0(x)) = 0, \quad x \in \Omega,  \quad 0 < t \leq T_{max},
\end{equation}
with given initial condition 
\[
I(x, 0) = I_0(x), \quad x \in \Omega, \quad t = 0.
\]
Here, as time increases, we approach the denoised image. 

In equation \eqref{tv1}, the term $\frac{1}{||\nabla I(x)||}$ should be slightly perturbed and replaced by $ \frac{1}{\sqrt{||\nabla I(x)||^2 + \beta^2}}$ 
\begin{equation}
\label{tv2}
I_t - \nabla \cdot \left( \frac{1}{\sqrt{||\nabla I(x)||^2 + \beta^2}}\nabla I \right) + \mu(I(x) - I_0(x)) = 0, \quad x \in \Omega,  \quad 0 < t \leq T_{max},
\end{equation}
and leading to the following unconstrained minimization problem \cite{chan1999nonlinear}
\[
\min_I  \int_{\Omega} \left( \sqrt{||\nabla I(x)||^2 + \beta^2} + \frac{\mu}{2}  (I(x) - I_0(x))^2 \right)  dx.
\]

In this work, we consider the Perona-Malik model for image denoising \cite{perona1990scale, perona1994anisotropic}. The Perona-Malik model differs slightly from the ROF model in formulating the nonlinear diffusion coefficient, but both aim to preserve edges while reducing noise. 
A nonlinear parabolic equation with anisotropic diffusion describes the Perona-Malik (PM) model
\[
I_t - \nabla \cdot (c(||\nabla I||) \nabla I) = 0, \quad t > 0, \quad x \in \Omega,
\]
where  $c(||\nabla I||)$ is an edge-stopping function, often chosen as
\[
c(||\nabla I||) = \frac{1}{1 + \left(\frac{||\nabla I||}{\lambda}\right)^2} 
\quad 
\text{or} 
\quad 
c(||\nabla I||) = \exp\left(-\left(\frac{||\nabla I||}{\lambda}\right)^2\right)
\]
and $\lambda$ controls edge sensitivity. This model performs selective diffusion based on gradient magnitude that diffuses more in flat regions (low gradients) and less near edges (high gradients), preserving edge details while reducing noise. 
Both PM and ROF models aim to preserve edges while denoising and giving an anisotropic diffusion based on the magnitude of the gradient. The ROF model this through the \( ||\nabla I|| \) term, while the PM model uses the edge-stopping function \( c(||\nabla I||) \) to control diffusion near edges. 

Furthermore, for the PM model, we may also introduce the regularization term \cite{voronin2023monolithic, nordstrom1990biased} and consider it in the frame of the following minimization problem \cite{charbonnier1997deterministic, atlas2014perona}
\[
\min_I  \int_{\Omega} \left( \ln \left(1 + \frac{||\nabla I(x)||^2}{\lambda^2}\right) + \frac{\mu}{2}  (I(x) - I_0(x))^2 \right) dx.
\]

Therefore, by combining the PM and ROF models and considering them as a solution to the time-dependent problem, we introduce the following unsteady nonlinear reaction-diffusion equation
\[
I_t(x, t) 
- \nabla \cdot (c(||\nabla I(x, t)||) \nabla I(x, t)) 
+ \mu(I(x, t) - I_0(x))
= 0, \quad t > 0, \quad x \in \Omega,
\]
with two edge-preserving nonlinear anisotropic diffusion coefficients: 
\[
c(||\nabla I||) = \frac{1}{1 + \left(\frac{||\nabla I||}{\lambda}\right)^2} \quad \text{(PM)}, 
\quad \quad 
c(||\nabla I||) = \frac{1}{\sqrt{||\nabla I||^2 + \beta^2}} \quad  \text{(ROF)}.
\]
and given initial conditions $I(x, 0) = I_0(x)$ in $\Omega$.
}

\subsection{\rev{Solution Techniques}}

\rev{
The classical solution technique for nonlinear anisotropic diffusion is based on explicit schemes. Explicit schemes are stable for very small time step sizes, which leads to poor efficiency \cite{weickert1998efficient}. 
In the ROF approach \cite{rudin1992nonlinear}, the authors propose explicit (forward Euler) time approximation to obtain a gradient descent scheme for a semi-discrete problem arising after spatial discretization with a finite-difference scheme. Similarly, in the PM approach \cite{perona1994anisotropic}, a conditionally stable explicit scheme is constructed for time approximation with finite difference space approximation. In both methods, due to the conditional stability of the explicit scheme, we have time step restriction, $\tau \leq 1/4$. 
Furthermore, the highly nonlinear and possible non-differentiable diffusion coefficient leads to the additional restriction for parameters choice $\lambda$ or $\beta$, where various approaches were introduced that make the anisotropic diffusion enhance edges and preserving stability \cite{perona1994anisotropic, van2018linear, tsiotsios2013choice}. Such approaches usually include specific algorithms for explicit schemes to calculate the dynamical value of the parameter that adapts to the current noise level \cite{canny1986computational, black1998robust, voci2004estimating}.
Moreover, solution techniques include various iterative methods in the context of iterative solutions of convex minimization problems. In \cite{knyazev2016accelerated}, the preconditioned conjugate gradient (PCG) acceleration of nonlinear iterative smoothing and Nesterov's acceleration are presented. In optimization, Nesterov's accelerated gradient and FISTA (fast iterative shrinkage/thresholding algorithm) are proposed to speed up convergence \cite{nesterov2018lectures, beck2009fast}. This technique can still be related to explicit schemes with acceleration in choice of time step size.
In this work, to avoid this additional complexity and avoid stability restriction to time step size, we consider stable semi-implicit approximation with constant parameters $\lambda$ or $\beta$.

Implicit time approximation leads to stable schemes and allows the use of a larger time marching step size and construct efficient algorithms with faster convergence to the final solution. 
An explicit scheme as a method of steepest descent with undesirable convergence properties  \cite{vogel1996iterative, luenberger1984linear}. Better convergence properties are obtained for Newton's method. In \cite{vogel1996iterative}, a rapid linear convergence was shown for lagged fixed point iterations related to a semi-implicit scheme with a broad range of time step size and $\beta$. For implicit schemes, a nonlinear primal-dual method is introduced for total variation denoising, and results show that the technique improves global convergence behavior more than the primal Newton's method \cite{vogel1996iterative, chan1999nonlinear}. 
%
In \cite{weickert1998efficient}, the authors present an additive operator splitting (AOS) for semi-implicit approximation of the nonlinear diffusion model. The split approach divides the multidimensional problem into a sequence of one-dimensional problems \cite{vabishchevich2013additive}. 
A multigrid method is widely used for solving resulting linear problems in each iteration \cite{vogel1996iterative, voronin2023monolithic}.

In this work, we consider solving the time-dependent nonlinear equation used for image denoising. We construct an explicit scheme and implicit scheme for time approximation with a finite-volume approximation by space. Stability estimates for both explicit and implicit time approximations are presented. We show that the implicit scheme is unconditionally stable and allows calculations with a large time step size but requires the solution of the large system of linear equations in each time iteration. To reduce the size of the system and perform a fast solution using an implicit scheme, we introduce a multiscale method.
}

\section{Problem formulation}

%

\rev{In this work, we separate a nonlinear flux in order to construct an approximation using a finite volume scheme. We use a finite volume scheme with a two-point flux approximation leading to the approximation of the nonlinear flux on the interface between two cells (pixels). Moreover, the grid is related to the image resolution given as a number of pixels in the vertical and horizontal directions $N_1 \times N_2$. Each pixel has a corresponding value $I(x) \in [0, 255]$ in a grey-scale representation and color value $I(x) = (I_1(x), I_2(x), I_3(x))$ with $I_i \in [0, 255]$ that correspond to the color model, such as RGB, YCbCr and other.

Let $I(x)$ be an image with $x = (x_1, x_2) \in  \Omega \subset \mathbb{R}^2$. Then, we consider the following nonlinear parabolic equation   in the computational domain $\Omega = [0, N_1] \times [0, N_2]$}
\begin{equation} 
\label{eq:1}
I_t + \nabla \cdot q(I, x) = 0, \quad x \in \Omega,  \quad 0 < t \leq T_{max}, 
\end{equation}
where $q$ is a nonlinear flux
 \[
q(I, x) = - c(||\nabla I||) \nabla I,
\]
with a given nonlinear coefficient $c$  
\rev{
\[
c(||\nabla I||) = \frac{1}{1 + \left(\frac{||\nabla I||}{\lambda}\right)^2} \quad \text{(PM)}, 
\quad \quad 
c(||\nabla I||) = \frac{1}{\sqrt{||\nabla I||^2 + \beta^2}} \quad  \text{(ROF)}.
\]
}
and $T_{max}$ is the final time,  $\nabla$ and $\nabla \cdot$ are the gradient and divergence operators. 

We consider \eqref{eq:1} with  zero flux boundary conditions 
\rev{
\[
q \cdot n = 0, \quad x \in \partial \Omega, \quad 0 < t \leq T_{max}, 
\]
}
and the initial condition
\[
I(x, 0) = I_0(x), \quad x \in \Omega, \quad t = 0, 
\]
and $I_0(x)$ is the initial noised image.

\begin{figure}[h!]
\centering
\includegraphics[width=0.4\linewidth]{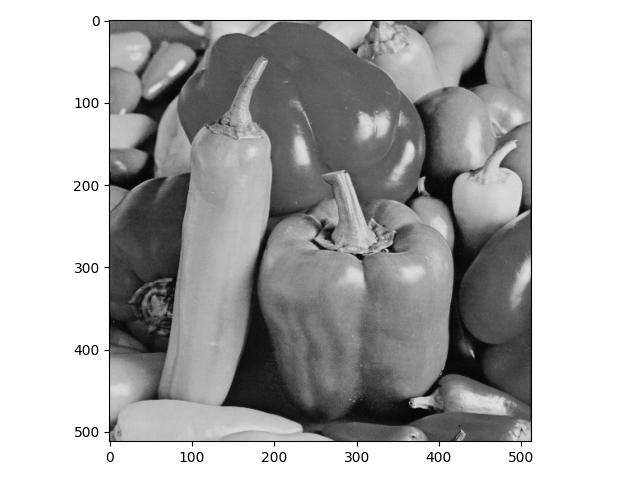}
\includegraphics[width=0.4\linewidth]{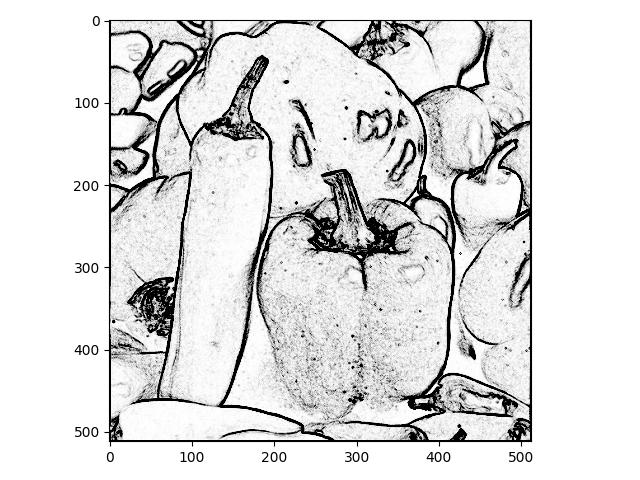}
\caption{Illustration of the reference image $I(x)$  (left) and $||\nabla I(x)||^2$ (right)}
\label{fig:img}
\end{figure}

\rev{
The nonlinear coefficient depends on the gradient of the image, i.e., $||\nabla I||^2 = (\nabla I, \nabla I)$, which is crucial in image denoising to preserve the edges of the image. Here, the term $||\nabla I||^2$ is discontinuous even if image $I(x,t)$ is smooth (see Figure \ref{fig:img}). Figure \ref{fig:img} presents an illustration for greyscale image $I(x)$ with resolution $N_1 = N_2 = 512$ and edge detection properties of image gradient $||\nabla I(x)||^2$. 
The problem formulation with a flux definition gives a better technique to linearize a nonlinear diffusion, leading to a more stable solution \cite{ewing2001modified,  chan1999nonlinear}. It was shown that this technique leads to the mixed formulation of the problem \cite{boffi2013mixed}. Furthermore, this formulation has some similarities to the primal-dual optimization method and gives a better convergence compared with a primal formulation, especially when term $\beta$ is very small \cite{chan1999nonlinear}. Here, the primal variable is $I(x,t)$, and the dual variable is $q(x, I)$. 
}

\subsection{\rev{Approximation by space}}

To solve the given nonlinear equation \eqref{eq:1}, we construct a discrete system at each time step using an appropriate approximation by space and time variables. An image can be represented on a structured grid with square cells where each cell represents a pixel so the value in the cell is the intensity of that pixel determining the color. We call a fully resolved approximation by space a fine grid. Later, we will introduce a coarse grid associated with a lower-resolution representation of the discrete problem. 

\rev{
Let $K_l$ be the fine grid cell (pixel of the image) and $\mathcal{T}_h = \cup_{l=1}^{N_h} K_l$, where $N_h$ is the number of cells such that $N_h = N_1 \times N_2$ and 
$x_{i,j}$ is the center of a cell  $K_l = K_{i,j}$ 
\[
x_{i,j} = ((i - 0.5) \cdot h, \, (j - 0.5) \cdot h),
\quad i = 1,..,N_1,
\quad j = 1,..,N_2,
\]
where $l=l(i,j)$ is the global cell index and $h$ is the distance between grid nodes (grid size). For example, we can set $l(i,j) = i \cdot N_2+j$. 
To approximate flows we introduce a grid $\mathcal{T}’_h$  with nodes
\[
x’_{i-0.5,j} = ( (i-0.5) \cdot h, \, j \cdot h)
\quad \text{ and } \quad
x’_{i, j-0.5} = ( i \cdot h, \, (j-0.5) \cdot h).
\]
We note that the computational domain $\Omega = [0, N_1] \times [0, N_2]$, so we have unit square cell $K_i$ (pixel) with mesh size $h = 1$.

\begin{figure}[h!]
\centering
\includegraphics[width=1\linewidth]{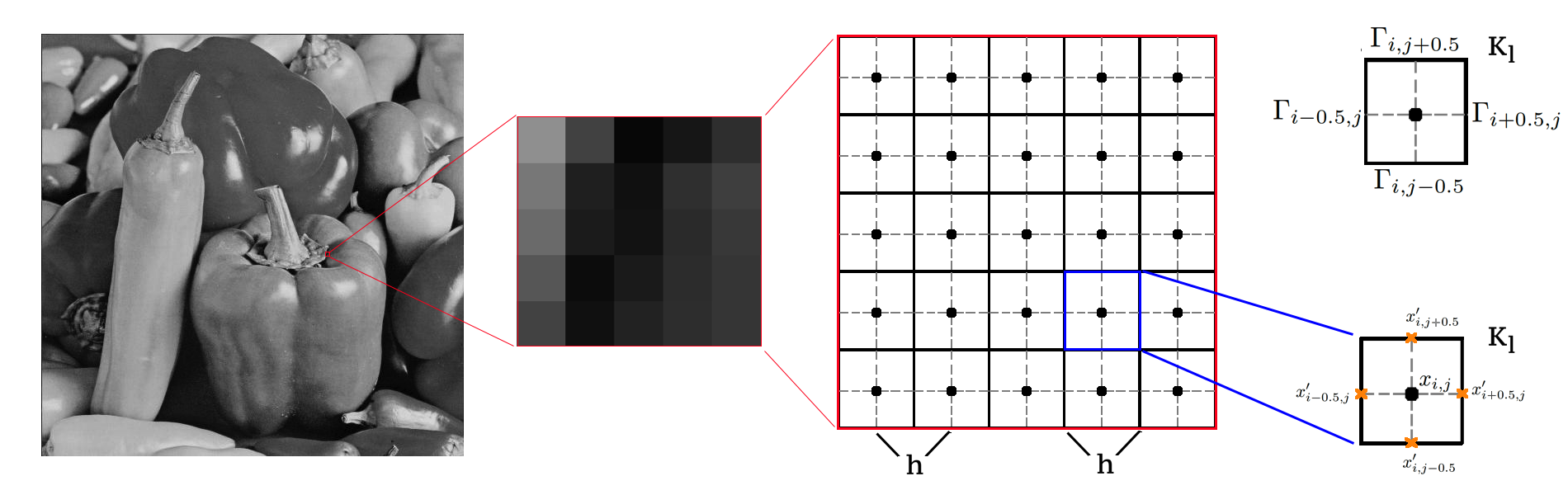}
\caption{\rev{Computational grid $\mathcal{T}_h$ and cell $K_l$, where $l=l(i,j)$ is the cell index}}
\label{fvm}
\end{figure}
A finite volume method can be considered as a cell-centered finite-difference scheme with a specific approximation of the discontinuous coefficient defined in each cell $K_l$. Next, we discuss the details of space approximation with relation to the finite-difference scheme and consider conservative approximation \cite{vasilyeva2018conserv}. 
To construct a conservative finite-difference scheme, we use the integral-interpolation method (balance method) \cite{samarskii2001theory}.
Let $q_{i-0.5,j}$, $q_{i+0.5,j}$, $q_{i, j-0.5}$ and $q_{i, j+0.5}$ be a fluxes on the cell interfaces, where $n$ is the outward normal to boundaries of the cell, $\Gamma_{i-0.5,j}$, $\Gamma_{i+0.5,j}$, $\Gamma_{i, j-0.5}$ and $\Gamma_{i, j+0.5}$ are the left, right, lower and upper bounds of the cell $K_l = K_{ij}  = [(i-1) \cdot h, \, i\cdot h] \times [(j-1) \cdot h, \, j \cdot h]$ with the center at the point $x_{i,j}$  (see Figure \ref{fvm}).

We write the conservation (mass balance) equation for $I_{i,j} = I(x_{i,j})$ on each cell $K_{ij}$
\begin{equation}
\label{m3}
(I_{i,j})_t |K_{ij}|  + 
q_{i+0.5,j} - q_{i-0.5,j} + q_{i, j+0.5} - q_{i, j-0.5}
= 0, 
\end{equation}
with 
\[
q_{i - 0.5,j} =
\int_{\Gamma_{i - 0.5,j}} q \cdot n \, ds, \quad
q_{i, j - 0.5} = 
\int_{\Gamma_{i, j - 0.5}} q \cdot n \, ds,
\]
and $|K_{ij}|$ is the volume of cell $K_{ij}$, $|K_{ij}| = h^2$.

Then we integrate flux $q(I, x) = - c(||\nabla I||) \nabla I$ over the cell interface
\[
\int_{\Gamma_{i - 0.5,j}} c^{-1} q \cdot n \, ds =
- \int_{\Gamma_{i - 0.5,j}} \nabla I \cdot n \, ds, \quad
\int_{\Gamma_{i, j - 0.5}} c^{-1} q \cdot n \, ds =
- \int_{\Gamma_{i, j - 0.5}} \nabla I \cdot n \, ds.
\]
Then introduce two-point approximation of the derivative
\[
\frac{1}{|\Gamma_{i - 0.5,j}|}
\int_{\Gamma_{i - 0.5,j}} \nabla I \cdot n \, ds =
\nabla I \cdot n|_{\Gamma_{i - 0.5,j}}
\approx
\frac{I_{i,j} - I_{i-1, j}}{h} ,
\]\[
\frac{1}{|\Gamma_{i, j - 0.5}|}
\int_{\Gamma_{i, j - 0.5}} \nabla I \cdot n \, ds =
\nabla I \cdot n|_{\Gamma_{i, j - 0.5}}
\approx
\frac{I_{i,j} - I_{i, j-1}}{h},
\]
where $|\Gamma_{i - 0.5,j}| = |\Gamma_{i, j - 0.5}| = h$, $h$ is the distance between nodes $x_{i,j}$ and $x_{i-1,j}$ (or $x_{i,j}$ and $x_{i, j-1}$).

Assuming $q \cdot n = const$ on the faces $\Gamma_{i - 0.5, j}$ and $\Gamma_{i, j- 0.5}$, we obtain the following approximation
\[
q_{i - 0.5,j} =
c_{i-0.5,j} (I_{i,j} - I_{i-1,j}),
\quad
q_{i, j - 0.5} =
c_{i, j-0.5} (I_{i,j} - I_{i, j-1} ),
\]
where 
\[
c_{i-0.5,j} = \left( \frac{1}{|\Gamma_{i - 0.5,j}|} \int_{\Gamma_{i - 0.5,j}} c^{-1} ds \right)^{-1}
\quad
c_{i, j-0.5} = \left( \frac{1}{|\Gamma_{i, j - 0.5}|} \int_{\Gamma_{i, j - 0.5}} c^{-1} ds \right)^{-1}.
\]
In the case of heterogeneous coefficients, the following relations (harmonic average) can be used
\[
c_{i-0.5,j} 
= \frac{2}{1/c_{i,j} + 1/c_{i-1,j}}, \quad
c_{i, j-0.5} 
= \frac{2}{1/c_{i,j} + 1/c_{i,j-1}}, .
\]

In the current application for image denoising, we have $h = 1$.  
Then, we have
\[
(I_{i,j})_t 
- c_{i+0.5,j} (I_{i+1,j} - I_{i,j})
+ c_{i-0.5,j} (I_{i,j} - I_{i-1,j} )
- c_{i, j+0.5} (I_{i, j+1} - I_{i,j})
+ c_{i, j-0.5} (I_{i,j} - I_{i, j-1}) = 0,
\]
for $i = 1,..,N_1$, $j = 1,..,N_2$.


Finally by introducing an indexing $l = l(i,j)$ and taking into account free flux boundary conditions, we obtain the following semi-discrete system for $I = (I_1, \ldots , I_{N_h})$
\begin{equation}
\label{fine}
I_t  + L(I) I = 0, 
\end{equation}
with 
\[
L = \{a_{lm}\}_{l,m=1}^{N_h} \quad 
a_{lm} = 
\left\{\begin{matrix}
\sum_j c_{lm}(I_l, I_m) & l = m, \\ 
-c_{lm}(I_l, I_m) & l \neq m
\end{matrix}\right. , 
\]
where $c_{lm}(I_l, I_m)= 2/(1/c(I_l) + 1/c(I_m))$.

Next, we define a uniform time approximation scheme. It is well known that an explicit time approximation to the time-dependent diffusion problem leads to a strong time step size restriction, and implicit schemes are widely used to solve such types of problems. In this work, we consider both explicit and implicit time approximation schemes in order to investigate their performance for image-denoising processes. In this work, we consider constant time step size and leave time adaptive schemes for a future investigation. 
}

\subsection{\rev{Approximation by time}}

\rev{
 Let $\tau = t^n - t^{n-1} = T_{max}/N_t$ be a fixed time step size, $I^n_l = I(x_l, t^n)$ be a solution at time $t^n$ at point $x_l$ and $N_t$ be a total number of time steps. Here, we use superscript $n$ to represent a time layer and subscript for space position $l$. 
To construct a time approximation, we use a general weighted $\theta$-scheme (two-level scheme) and apply linearization from a previous time step \cite{samarskii2001theory, vabishchevich2013additive} 
\begin{equation}
\label{eq:fine}
\frac{I^{n+1} - I^n}{\tau} 
+ L^n (\theta I^{n+1} + (1-\theta) I^n)  = 0,
\end{equation}
where $L^n = L(I^n)$. 
Here, we have the forward Euler (explicit) and backward Euler (implicit) schemes for $\theta=0$ and $1$.
We let $||u|| =  \sqrt{(u,u)}$ and $||u||_L = \sqrt{(Lu,u)}$ be an $L_2$ and energy norms with $(u, v) = v^T u$ and $(L u, v) = v^T L u$.

\begin{theorem}
\label{t:t1}
The solution of the discrete problem \eqref{eq:fine} is stable with $\theta \geq 1/2$ and satisfies the following estimate
\begin{equation}
\label{t1}
||I^{n+1}||_{L^n}^2 \leq ||I^n||_{L^n}^2.
\end{equation}
\end{theorem}
\begin{proof}
We write the equation \eqref{eq:fine} as follows
\[
\left(\mathcal{I} + \tau \theta L^n \right)\frac{I^{n+1} - I^{n}}{\tau} + L^n I^n = 0,
\]
where $\mathcal{I}$ is the identity matrix. 

We multiply to $(I^{n+1} - I^{n})/\tau$ with $I^n = (I^{n+1} + I^n)/2 - (I^{n+1} - I^n)/2$ and obtain 
\begin{equation}
\label{t:stab1}
\begin{split}
& \left( (\mathcal{I} + \tau \theta L^n) \frac{I^{n+1} - I^{n}}{\tau}, \frac{I^{n+1} - I^{n}}{\tau} \right) 
 + \left( L^n I^{n}, \frac{ I^{n+1} - I^{n} }{\tau} \right)\\
&= \left( (\mathcal{I} + \tau \theta L^n) \frac{I^{n+1} - I^{n}}{\tau}, \frac{I^{n+1} - I^{n}}{\tau} \right)  
+ \left( L^n \frac{I^{n+1} + I^n}{2},   \frac{I^{n+1} - I^{n}}{\tau} \right) 
- \left( L^n \frac{I^{n+1} - I^n}{2},   \frac{I^{n+1} - I^{n}}{\tau} \right) \\
&= \left( \left(\mathcal{I} + \tau \left(\theta - \frac{1}{2} \right) L^n \right) \frac{I^{n+1} - I^{n}}{\tau}, \frac{I^{n+1} - I^{n}}{\tau} \right)  
+ \left( L^n \frac{I^{n+1} + I^n}{2},   \frac{I^{n+1} - I^{n}}{\tau} \right) 
= 0.
\end{split}
\end{equation}
We have symmetric semi-positive definite matrix $L^n = (L^n)^T \geq 0$, then 
\[
( L^n (I^{n+1} + I^n),  I^{n+1} - I^n) = L^n I^{n+1}, I^{n+1}) - (L^n I^n, I^n).
\]
Therefore, the two-level  scheme  is stable with $(\theta - 1/2) \geq 0$ and stability estimate holds.
\end{proof}

Based on the theorem \ref{t:t1}, we have an unconditionally stable implicit scheme for $\theta=1$
\begin{equation}
\label{eq:fine-im}
\frac{I^{n+1} - I^n}{\tau} + L^n I^{n+1} = 0.
\end{equation}
We have conditional stability for the explicit scheme with $\theta=0$
\begin{equation}
\label{eq:fine-ex}
\frac{I^{n+1} - I^n}{\tau} + L^n I^n  = 0.
\end{equation}
For operator $L^n$, we have
\[
L^n \geq d \ \bar{c}^n \ \mathcal{I},
\]
with $d = 2$ is the computational domain dimension, $\Omega \subset \mathbb{R}^2$. 
Here $\mathcal{I}$ is the identity matrix and $\bar{c}^n = \max_{K_l} c^n_{l}$ is the maximum diffusion coefficient \cite{vasilyeva_ms_rd, vasilyeva2024decoupled, vasilyeva2023efficient}.

From \eqref{t:stab1}, we have the following condition for stability of the scheme
\[
\mathcal{I} + \tau \left(\theta - \frac{1}{2} \right) L^n 
= \mathcal{I} - \tau \frac{1}{2} L^n 
\geq (1 -  \tau \bar{c}^n) \mathcal{I}
\geq 0.
\]
Therefore, the following theorem holds.

\begin{theorem}
\label{t:t2}
The solution of the explicit scheme \eqref{eq:fine-ex} is conditionally stable with $\tau < 1/\bar{c}^n$, and the a priori estimate holds
\[
||I^{n+1}||_{L^n}^2 \leq ||I^n||_{L^n}^2.
\]
\end{theorem}

We note that the nonlinear coefficient $c_l$ is defined on each cell based on the current gradient of the solution
\[
c_l = \frac{1}{1 + \left(\frac{||\nabla I_l||}{\lambda}\right)^2} \quad \text{(PM)}, 
\quad \quad 
c_l = \frac{1}{\sqrt{||\nabla I_l||^2 + \beta^2}} \quad  \text{(ROF)}.
\] 
Therefore, it would be interesting to investigate the influence of the parameters $\beta$ and $\lambda$ on the stability of the scheme. In the literature, it was shown that the parameters can be adaptively chosen to preserve upper bound and perform calculations with $\tau < 1/4$, \cite{perona1994anisotropic, chambolle2004algorithm}. In this work, we use fixed parameter values and further investigate stability numerically. 

Finally, we have the following two algorithms with a full fine-scale resolution $N_h = N_1 \times N_2$ ($N_1$ and $N_2$ are the number of pixels in vertical and horizontal directions):
\begin{itemize}
\item \textit{Explicit scheme} (Ex-PM and Ex-ROF): 

For a given initial noised image $I_0(x)$, we update the solution as follows at each time iteration ($n=1,2,...,N_t$)
\begin{equation}
\label{eq:alg-ex}
I^{n+1} = I^n -\tau L^n I^{n}.
\end{equation}

\item \textit{Implicit scheme} (Im-PM and Im-ROF):

For a given initial noised image $I^0(x) = I_0(x)$, we solve the following system of linear equations for each time iteration ($n=1,2,...,N_t$)
\begin{equation}
\label{eq:alg-im}
A^{n} I^{n+1} = b^{n},  \quad 
A^{n} = (\mathcal{I} + \tau L^{n}), \quad 
b^{n} = I^{n}.
\end{equation}
\end{itemize}
}
For high-resolution images, the number of unknowns \rev{($N_h$)} in the resulting system is very large and, therefore, computationally expensive to solve. 
Next, we present the construction of an accurate and efficient computational algorithm \rev{for an implicit scheme} based on the local denoising process and multiscale basis functions.

\section{Multiscale method for image denoising}

\begin{figure}[h!]
\centering
\includegraphics[width=0.5\linewidth]{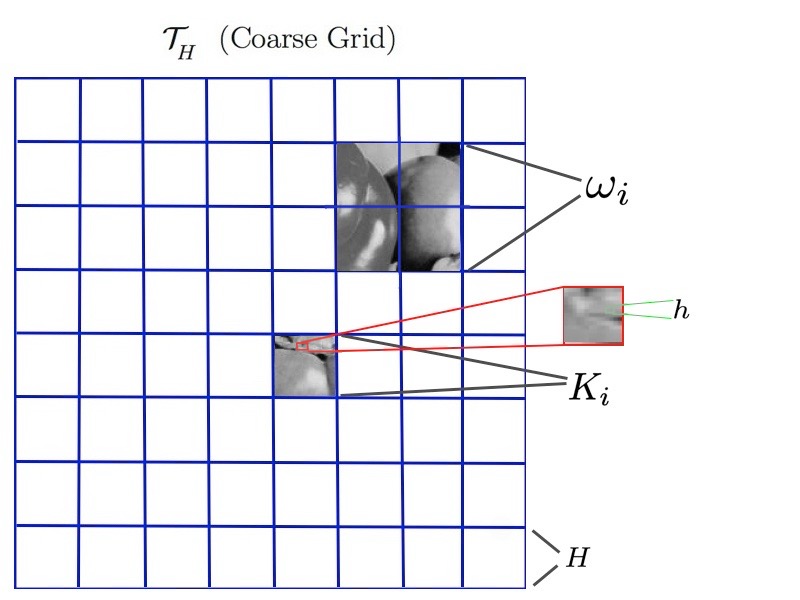}
\caption{Illustration of $8 \times 8$ coarse grid ($\mathcal{T}_H$) with coarse cell $K_i$ and local domain $\omega_i$ corresponded to the local support of multiscale basis function}
\label{fig:ms}
\end{figure}

We start with the one-channel image (greyscale) representation to illustrate basis construction. The computational domain $\Omega \in \mathbb{R}^2$ and $\mathcal{T}_h$  is the structured fine grid with $N_1 \times N_2$ cells (pixels). In order to construct a low-resolution representation and present a multiscale model order reduction technique, we start by defining a coarse grid  $\mathcal{T}_H$ with $H >> h$ (Figure \ref{fig:ms}). 

Let $\mathcal{T}_H$ be the coarse grid defined as follow
\[
\mathcal{T}_H  = \sum_{i=1}^{N_c} K_i, 
\]
where $K_i$ is the coarse grid cell and $N_c$ is the number of coarse grid cells. We consider a conforming approach, where the coarse grid is represented as an agglomeration of the fine grid cells. For example, if we have an input image with resolution $512 \times 512$ pixels, then for $8 \times 8$ coarse grid, we have coarse cell $K_i$ with $64 \times 64$ pixels or for $16 \times 16$ coarse grid we have a coarse cell with $32 \times 32$ pixels. We can relate this coarse grid to an image with low resolution. 

Let, $\omega_i$ be a coarse scale neighborhood related to the coarse grid node $x_i$ and constructed as a combination of the several coarse cells that contain the corresponding coarse grid node. Then, for our case with a structured grid related to the image representation, the local domain $\omega_i$ will contain four coarse cells for interior nodes (see Figure \ref{fig:ms}).  
We construct a multiscale basis function in each local domain $\omega_i$ by performing local denoising and solving local spectral problems.

\subsection{Local image denoising}

In this work, we use the multiscale finite element method concept to create a coarse-scale approximation for a nonlinear time-dependent equation. Specifically, we are developing a multiscale basis function to represent the fine-scale behavior of the solution accurately. Since the problem we are dealing with is nonlinear, the basis function should address the nonlinearity and image dynamics during the multiscale space construction. We also consider the relevant nonlinear coefficients as part of the denoised image representation for our particular denoising application, and we integrate this understanding of the equation's behavior into constructing the multiscale basis functions.


\begin{figure}[h!]
\centering
\includegraphics[width=1\linewidth]{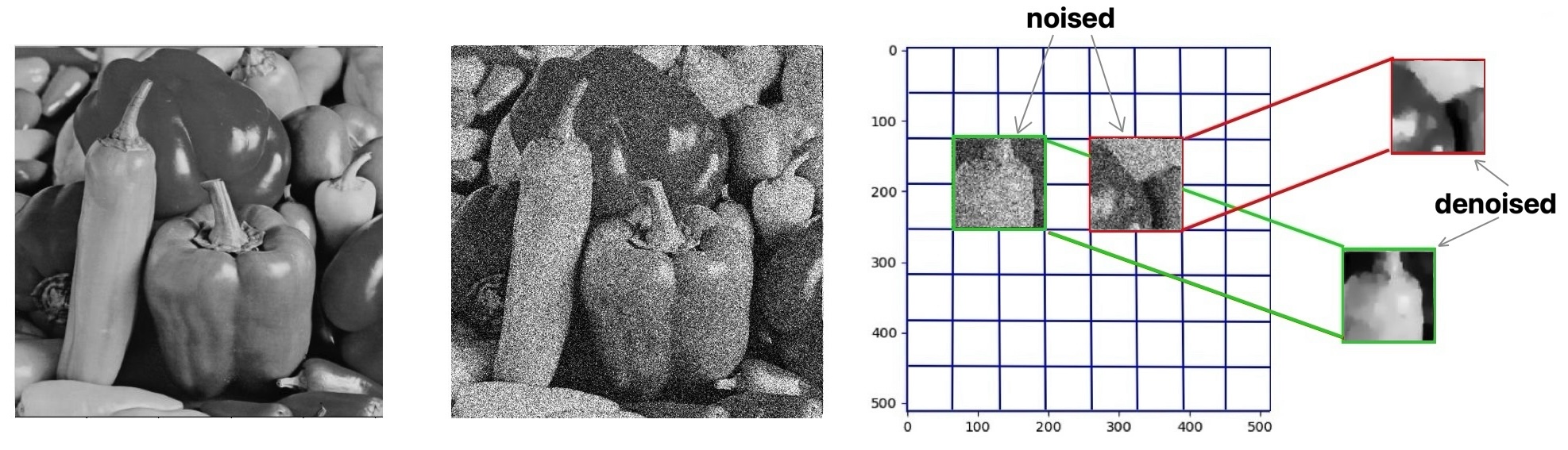}
\caption{Illustration of the original image (first plot), noised image (second plot), coarse grid with two local domains ($\omega_{29}$ and $\omega_{32}$) with local denoising (third plot)}
\label{fig:ms2}
\end{figure}

The considered nonlinear parabolic equation has a heterogeneous diffusion coefficient that depends on the initial noisy image. Therefore, if we use use a denoised image for coefficient calculation, we could obtain better basis functions for multiscale simulations with better feature detection.  Since all basis functions are calculated in local domains and can be done independently in parallel, we can perform a preliminary denoising for each local image and then use them to calculate the coefficient for eigenvalue problems.  The effect of this preliminary local denoising on the constructed eigenvectors will be illustrated in the next section.

Let us illustrate the local denoising process \rev{in Figure \ref{fig:ms2}}. 
The original image is depicted in the first plot. 
Next, we add noise to the original image and show it on the second plot. The illustration of the $8 \times 8$ coarse grid with two local domains of basis support is shown in the third plot. We considered two local domain $\omega_{29}$ (green) and $\omega_{32}$ (red). 
The results of denoising in local domains are depicted in the third plot. It is important to note that this method is lightweight, due to the smaller size of the domain and independent local calculations. 
\rev{
Various techniques for the local denoising process in a small local domain $\omega_i$ might be used. For example, classic denoising techniques such as BM3D and wavelet denoising of total variational minimization might be used. In this work, the goal is to make a local denoising process fast and simple. Furthermore, local denoising is mainly oriented to smoothing the local image, making basis functions smooth and more accurate. For this purpose, we might use a few iterations of linear diffusion denoising with a large time step size and implicit approximation by time for stability. 
Moreover, we may use a local denoised image as a local multiscale basis function after multiplying to the partition unity function and representing a detailed algorithm at the end of this section.
}

Next, we illustrate how this preprocessing step affects the resulting multiscale basis functions and the ability to encapsulate the main features of the image in multiscale approximation.

\subsection{Multiscale space construction}

In order to construct multiscale basis functions in each local domain $\omega_i$, we solve the following generalized eigenvalue problem:
\begin{equation} 
\label{eq:sp}
S^{\omega_i} \psi^{\omega_i}_l = 
\lambda^{\omega_i}_l D^{\omega_i} \psi^{\omega_i}_l,
\end{equation}
with
\[
S^{\omega_i} = \{s_{ij}\}, \quad 
s_{ij} = 
\left\{\begin{matrix}
\sum_j \tilde{W}_{ij} & i = j, \\ 
-\tilde{W}_{ij} & i \neq j
\end{matrix}\right. , \quad 
D^{\omega_i}  = \{d_{ij}\}, \quad 
d_{ij} = 
\left\{\begin{matrix}
a_{ii} & i = j, \\ 
0 & i \neq j
\end{matrix}\right., \quad
i, j = 1,...,N_h^{\omega_i},
\]\[
\tilde{W}_{ij} = \tilde{c}_{ij} \frac{|E_{ij}|}{d_{ij}}, \quad 
\tilde{c}_{ij} = \frac{2}{\frac{1}{c(\tilde{I}_i)} + \frac{1}{c(\tilde{I}_j)}}, \quad 
c(\tilde{I}_i) = \frac{1}{1 + \frac{||\nabla \tilde{I}_i||^2}{\lambda^2} },
\]  
where $\psi^{\omega_i}_l$ and $\lambda^{\omega_i}_l$ are the eigenvectors and eigenvalues,
$\tilde{I}$ is the local denoised image and $N^{\omega_i}_h$ is the number of coarse grid cells in local domain $\omega_i$.


\begin{figure}[h!]
\centering
\begin{subfigure}[b]{1\textwidth}
\includegraphics[width=1\linewidth]{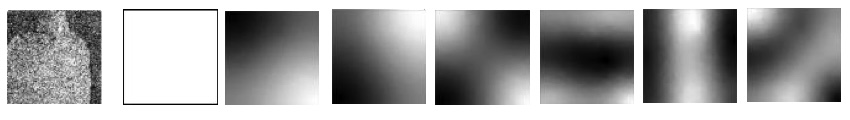}
\caption{Without local denoising}
\end{subfigure}
\begin{subfigure}[b]{1\textwidth}
\includegraphics[width=1\linewidth]{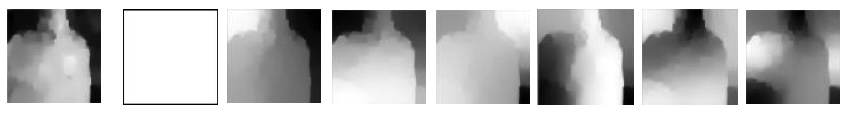}
\caption{With local denoising}
\end{subfigure}
\caption{Illustration of eigenvectors for noised and denoised local images (1st plot: local image, $I(x)$. 2nd-8th plots: eigenvectors corresponded to first seven smallest eigenvalues)}
\label{fig:ms3b}
\end{figure}

In Figure \ref{fig:ms3b}, 
we present the illustration of the local eigenvectors corresponding to the first seven smallest eigenvalues. We considered two cases for image-induced coefficient $c(||\nabla \tilde{I}||)$: (1) an initial noised image, $\tilde{I} = I_0(x)$; and (2) locally denoised image discussed above. 
We plot eigenvectors for two local domains $\omega_{29}$ and  $\omega_{32}$. 
We observe that the accuracy of the local image used for coefficient calculations in generalized eigenvalue problem highly affect to the resulting eigenvectors. We observe that the noised image produces smooth eigenvectors that do not preserve the fine-scale features. However, a local denoising procedure leads to eigenvectors with desired feature-preserving properties.

\rev{
In addition, we investigate the performance of the analytic approach for 'cheap' defining local spectral basis functions without solving local spectral problems. This approach is motivated by a linear denoising process using the Laplace operator. Due to the constant diffusion coefficient, this operator does not preserve edges and multiscale features of the image; however, it may be used to compare with a feature-preserving spectral problem \eqref{eq:sp}. The main advantage of this approach is the existence of an analytic formula that can be used to define local eigenvectors corresponding to the smallest eigenvectors. Similarly to the heterogeneous case with feature-preserving operator \eqref{eq:sp}, we will use the analytic formula for the Laplace operator with free (Neuman) boundary conditions. 
}

\subsection{Coarse-scale approximation and \rev{algorithms}}

In order to construct a multiscale space, we choose eigenvectors that correspond to $M_i$ smallest eigenvalues ($\lambda^{\omega_i}_1 < \lambda^{\omega_i}_2 <... <\lambda^{\omega_i}_{M_i}$). Then, we define a projection matrix as follows
\begin{equation} 
\label{ms-r}
R = \left[ 
\chi^1 \psi^{\omega_1}_1, \ldots, \chi^1 \psi^{\omega_1}_{M_1}
\ldots
\chi^{N_c} \psi^{\omega_{N_c}}_1, \ldots, \chi^{N_c} \psi^{\omega_{N_c}}_{M_{N_c}}
\right]^T.
\end{equation}
where $\chi^i$ is the linear partition of unity functions, and $N_c$ is the number of the local domains (number of coarse grid nodes). 

We use the projection matrix $R$ to project the fine grid system \eqref{eq:fine-im} to the coarse grid
\begin{equation} 
\label{eq:coarsep}
\rev{\mathcal{I}_H}  (I^{n+1}_H - I^{n}_H) + \tau L_H^{n} I_H^{n+1} = 0,
\end{equation}
with  
\begin{equation}
L_H^n = R L_H^n R^T, \quad 
\rev{\mathcal{I}_H} = R R^T,
\end{equation}
\rev{and 
\begin{equation}
L^n = L(I^n_{ms}), \quad I^n_{ms} = R^T I_H^n, \quad I^0_{ms} = R^T (RR^T)^{-1} R \ I^0.
\end{equation}
}
Note that, the size of the system is $DOF_H = \sum_{i=1}^{N_c} M_i$, $M_i$ is the number of local multiscale basis functions in $\omega_i$ and $N_c$ is the number of coarse grid vertices.  For the numerical investigation, we set to take the same number of basis functions in each local domain ($M_i = M$), and therefore, we have $DOF_H = M \cdot N_c$. 

Finally, the multiscale method for image denoising process can be represented as follows \rev{(\textbf{GMs})}:
\begin{enumerate}
\item Construct a coarse grid $\mathcal{T}_H$ and define local domains $\omega_i$, $i = 1,\ldots,N_c$.
\item In each local domain $\omega_i$ and corresponded local noised image $I_0^{\omega_i}$:
\begin{itemize}
\item Perform local denoising process to find $\tilde{I}$ in $\omega_i$.
\item Solve local generalized eigenvalue problems \eqref{eq:sp} to choose eigenvectors $\psi^{\omega_i}_l$ that correspond to $M_i$ smallest eigenvalues, $l=1,\ldots,M_i$.
\item Define multiscale basis functions by multiplying to partition of unity functions,  $\phi^{\omega_{i}}_{l} = \chi^{i} \psi^{\omega_{i}}_{l}$ 
\end{itemize}
\item Form projection matrix $R$ in \eqref{ms-r} and solve coarse-scale system \eqref{eq:coarsep}.
\end{enumerate}

\rev{
In addition to the GMs algorithm, we will study the approach with analytically defined eigenvectors relative to linear denoising with Laplace operator (\textbf{AMs}):
\begin{enumerate}
\item Construct a coarse grid $\mathcal{T}_H$ and define local domains $\omega_i$, $i = 1,\ldots,N_c$.
\item In each local domain $\omega_i$ and corresponded local noised image $I_0^{\omega_i}$:
\begin{itemize}
\item Perform local denoising process to find $\tilde{I}$ in $\omega_i$.
\item Define multiscale basis functions by multiplying to partition of unity functions,  $\phi^{\omega_{i}}_{l} = \chi^{i} \tilde{\psi}^{\omega_{i}}_{l}$, where $\tilde{\psi}^{\omega_{i}}_{l}$ are given by analytic formula. 
\end{itemize}
\item Form projection matrix $R$ and solve coarse-scale system \eqref{eq:coarsep}.
\end{enumerate}

Moreover, we may use a local denoised image as a local multiscale basis function and construct the following algorithm (\textbf{Ms}):
\begin{enumerate}
\item Construct a coarse grid $\mathcal{T}_H$ and define local domains $\omega_i$, $i = 1,\ldots,N_c$.
\item In each local domain $\omega_i$ and corresponded local noised image $I_0^{\omega_i}$:
\begin{itemize}
\item Perform local denoising process to find $\tilde{I}$ in $\omega_i$.
\item Define multiscale basis functions by multiplying to partition of unity functions,  $\phi^{\omega_{i}} = \chi^{i} \tilde{I}^{\omega_{i}}$. 
\end{itemize}
\item Form projection matrix $R$ and solve coarse-scale system \eqref{eq:coarsep}.
\end{enumerate}
In this algorithm, we will also compare results of initial noised image projection to multiscale space ($I^0_{ms} = R^T (RR^T)^{-1} R \ I^0$) and denote method as (\textbf{Ms$^0$}).
}

\section{Numerical results}

We present numerical results for the proposed multiscale method for greyscale and color images. We start with a grayscale image, then we extend a proposed method to a color image and present results for \rev{several datasets}.

Let $u = I^n_{ms}(x)$ be a solution using the multiscale method and $v = I(x)$ be a reference solution (original image without noise). 
To compare the numerical results with a reference solution $v = I(x)$, we use the following metrics
\begin{itemize}
\item RRMSE (relative root mean squared error or relative $L_2$ error)  provides a measure of the differences between the reconstructed image and the reference image 
\[
RRMSE  =  \frac{||u-v||_{L_2}}{||v||_{L_2}}, \quad 
||u||^2_{L_2} = (u, u).
\]
Lower RRMSE values indicate that the reconstructed image is closer to the reference image, while RRMSE = 0 means that the images are identical.

\item SSIM (structural similarity index) is a metric that measures the similarity between two images
\[
SSIM(u, v) = l(u, v) \cdot c(u, v) \cdot s(u, v)
\]
where $l(u, v)$, $c(u, v)$ and $s(u, v)$ are  luminance, contrast, and structure 
\[
   l(u, v) = \frac{2\mu_u \mu_v + C_1}{\mu_u^2 + \mu_v^2 + C_1}, \quad 
   c(u, v) = \frac{2\sigma_u \sigma_v + C_2}{\sigma_u^2 + \sigma_v^2 + C_2}, \quad 
   s(u, v) = \frac{\sigma_{uv} + C_2/2}{\sigma_u \sigma_v + C_2/2},
\]
and
$\mu_u$ and $\mu_v$ are the mean of images $u$ and $v$,  
$\sigma_u^2$  and $\sigma_v^2$ are the variance of $u$ and $v$,  
$\sigma_{uv}$ is the covariance of images  $u$ and $v$, $C_1$ and $C_2$ are constants to stabilize the division with weak denominators \cite{wang2009mean}.  
Higher SSIM  values indicate higher similarity between the two images, where SSIM = 1 means that images are identical.

\item PSNR (peak signal-to-noise ratio) is a metric used to measure the quality of a reconstructed image compared to reference image 
\[
PSNR(u,v) = 10 \cdot \log_{10} \left( \frac{MAX^2(u)}{MSE(u,v)} \right)
\]
where $MAX(u)$ is the maximum possible pixel value of the image (for an 8-bit image, this value is 255) and $MSE(u,v) = ||u-v||_{L_2}^2$ is the mean squared error.  
Higher PSNR values indicate better quality of the reconstructed image.
\end{itemize}

\rev{
The investigation includes the following aspects of the proposed approach: 
\begin{itemize}
\item \textit{Test 1.} Comparison of the fine-scale methods with explicit approximation with traditional techniques. 
\item \textit{Test 2.} Stability and parameters sensitivity analysis of explicit and implicit schemes. 
\item \textit{Test 3.} Performance of the multiscale approximation on the coarse grid.
\item \textit{Test 4.} Extension to color images and experiments on high-resolution datasets.
\end{itemize} 
The proposed PDE-based approach with explicit and implicit time approximation and a multiscale approach were implemented in Python.
All tests are performed with varying noise levels, where we used a Gaussian noise and added it iteratively with ascending noise variance till we reached a given noise level in percentage. We investigate the performance of the denoising algorithms for images with 20 and 40 \% of noise (RRMSE). Simulations are performed on a MacBook Pro with an Apple M2 Max chip and 32 GB memory. }

\subsection{\rev{Comparison with traditional techniques (Test 1)}}

\rev{
We chose five classic greyscale images to perform comparison of the fine-scale methods with explicit approximation with traditional techniques. In Figure \ref{fig:ref}, we represent five images that we use in our first test. All images are in greyscale with $N_1 \times N_2 = 512 \times 512$ pixels. 

\begin{figure}[h!]
\centering
\includegraphics[page=1,width=0.18\linewidth]{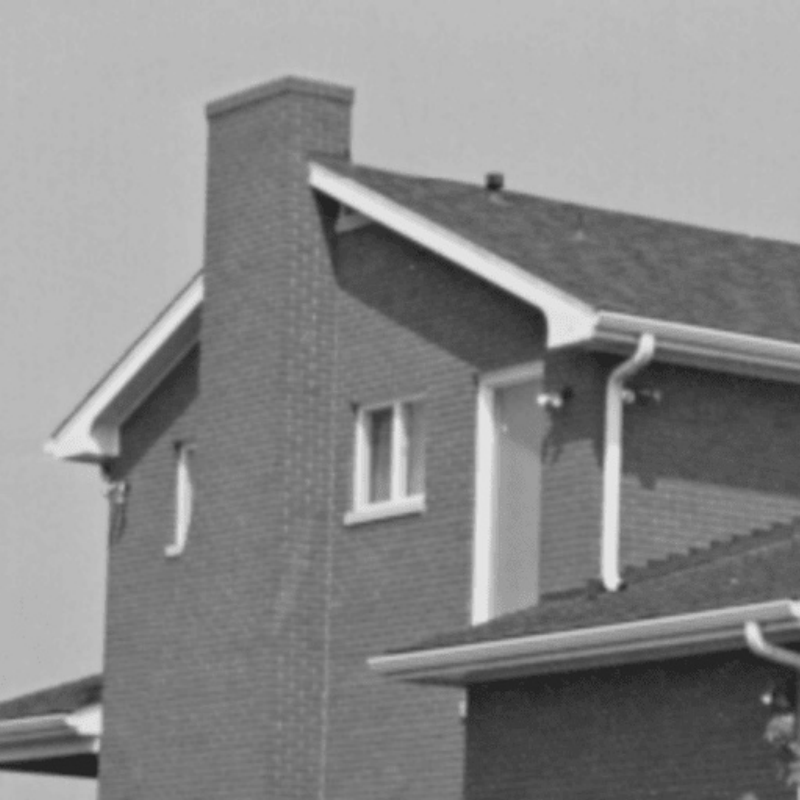}  \ \ 
\includegraphics[page=2,width=0.18\linewidth]{img-ref.pdf}  \ \ 
\includegraphics[page=3,width=0.18\linewidth]{img-ref.pdf}  \ \ 
\includegraphics[page=4,width=0.18\linewidth]{img-ref.pdf}  \ \ 
\includegraphics[page=5,width=0.18\linewidth]{img-ref.pdf}
\caption{\rev{Test 1. Reference image. Image 1 (House), Image 2 (Peppers), Image 3 (Barbara), Image 4 (Baboon) and Image 5 (Boat) (from left to right)}}
\label{fig:ref}
\end{figure}

\begin{table}[h!]
\begin{tabular}{| c | c | c | c | c | c | }
\hline
 & \tiny{RRMSE/SSIM/PSNR} & \tiny{RRMSE/SSIM/PSNR} & \tiny{RRMSE/SSIM/PSNR} & \tiny{RRMSE/SSIM/PSNR} & \tiny{RRMSE/SSIM/PSNR}\\
\hline
& 1.House & 2.Peppers & 3.Barbara & 4.Baboon & 5.Boat\\
\hline
\multicolumn{6}{|c|}{Initial noised image with 20\% of nosie}\\
\hline
$I_0(x)$ & 21.29/0.16/18.38 & 20.33/0.26/19.66 & 20.44/0.41/20.26 & 20.51/0.48/19.52 & 20.77/0.33/19.06 \\ 
\hline
Ex-PM & 4.52/0.84/31.85 & 5.97/0.84/30.31 & 10.75/0.72/25.84 & 11.50/0.67/24.54 & 7.87/0.72/27.49 \\ 
Ex-ROF & 4.19/0.86/32.49 & 5.93/0.85/30.37 & 10.66/0.73/25.92 & 11.27/0.68/24.72 & 7.74/0.73/27.63 \\ 
\hline
Wvlt & 5.79/0.77/29.69 & 8.11/0.71/27.65 & 11.36/0.64/25.37 & 12.70/0.57/23.68 & 9.40/0.63/25.94 \\ 
TV & 5.57/0.73/30.03 & 6.15/0.79/30.06 & 11.00/0.70/25.64 & 11.48/0.62/24.56 & 7.85/0.70/27.50 \\ 
NLM & 4.09/0.85/32.70 & 5.86/0.83/30.47 & 7.56/0.82/28.90 & 12.84/0.56/23.58 & 7.95/0.71/27.39 \\ 
BM3D & 3.32/0.89/34.53 & 4.83/0.86/32.15 & 6.01/0.87/30.89 & 10.19/0.68/25.59 & 6.49/0.77/29.16 \\ 
\hline
\multicolumn{6}{|c|}{Initial noised image with 40\% of noise}\\
\hline
$I_0(x)$ & 40.98/0.06/12.69 & 40.17/0.11/13.75 & 41.04/0.19/14.21 & 41.23/0.22/13.45 & 41.22/0.14/13.10 \\ 
\hline
Ex-PM & 7.47/0.79/27.48 & 9.91/0.76/25.90 & 14.87/0.59/23.03 & 15.84/0.45/21.76 & 11.88/0.58/23.91 \\ 
Ex-ROF & 6.87/0.82/28.20 & 9.59/0.77/26.19 & 14.51/0.60/23.24 & 15.26/0.48/22.09 & 11.44/0.60/24.24 \\ 
\hline
Wvlt & 8.57/0.74/26.28 & 12.25/0.65/24.06 & 15.86/0.54/22.47 & 16.73/0.38/21.29 & 13.32/0.55/22.91 \\ 
TV & 9.38/0.59/25.50 & 10.21/0.66/25.65 & 14.38/0.58/23.32 & 14.98/0.46/22.25 & 12.04/0.56/23.79 \\ 
NLM & 7.74/0.77/27.17 & 10.67/0.72/25.26 & 13.10/0.66/24.13 & 16.53/0.42/21.39 & 12.74/0.58/23.30 \\ 
BM3D & 5.63/0.82/29.93 & 7.99/0.78/27.77 & 9.80/0.75/26.65 & 14.12/0.51/22.76 & 10.01/0.66/25.39 \\ 
\hline
\end{tabular}
\caption{\rev{Test 1. Image denoising using various algorithms}}
\label{table:t1}
\end{table}

\begin{figure}[h!]
\centering
\begin{subfigure}[b]{0.14\textwidth}
Noised, $I_0(x)$
\end{subfigure}
\begin{subfigure}[b]{0.14\textwidth}
Ex-PM
\end{subfigure}
\begin{subfigure}[b]{0.14\textwidth}
Ex-ROF
\end{subfigure}
\begin{subfigure}[b]{0.14\textwidth}
Wvlt
\end{subfigure}
\begin{subfigure}[b]{0.14\textwidth}
TV
\end{subfigure}
\begin{subfigure}[b]{0.14\textwidth}
NLM
\end{subfigure}
\begin{subfigure}[b]{0.1\textwidth}
BM3D
\end{subfigure}

\vspace{5pt} 

\begin{subfigure}[b]{1\textwidth}
\includegraphics[page=1,width=0.13\linewidth]{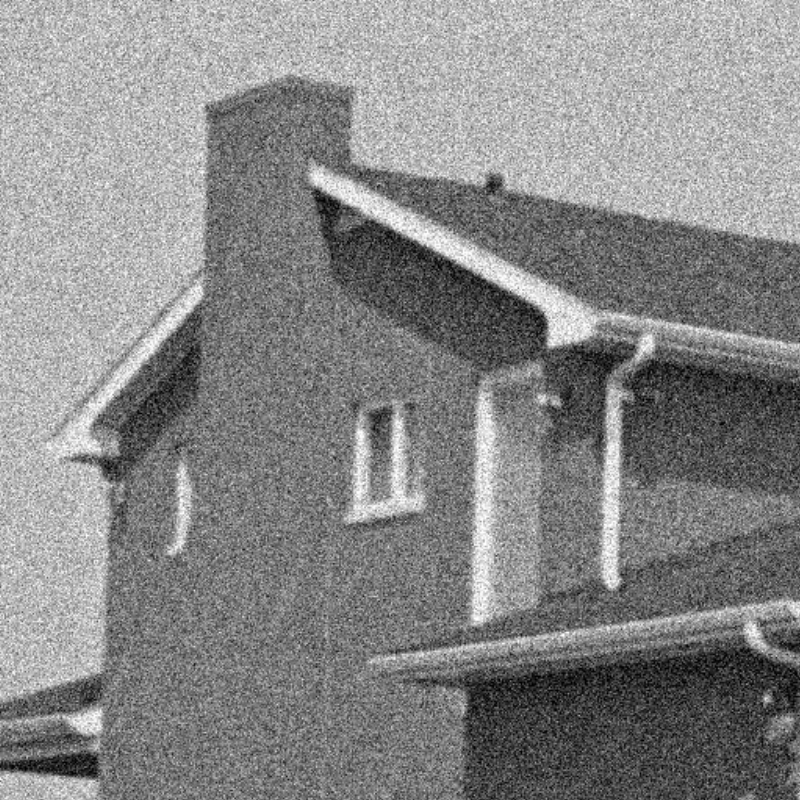} \ \ 
\includegraphics[page=1,width=0.13\linewidth]{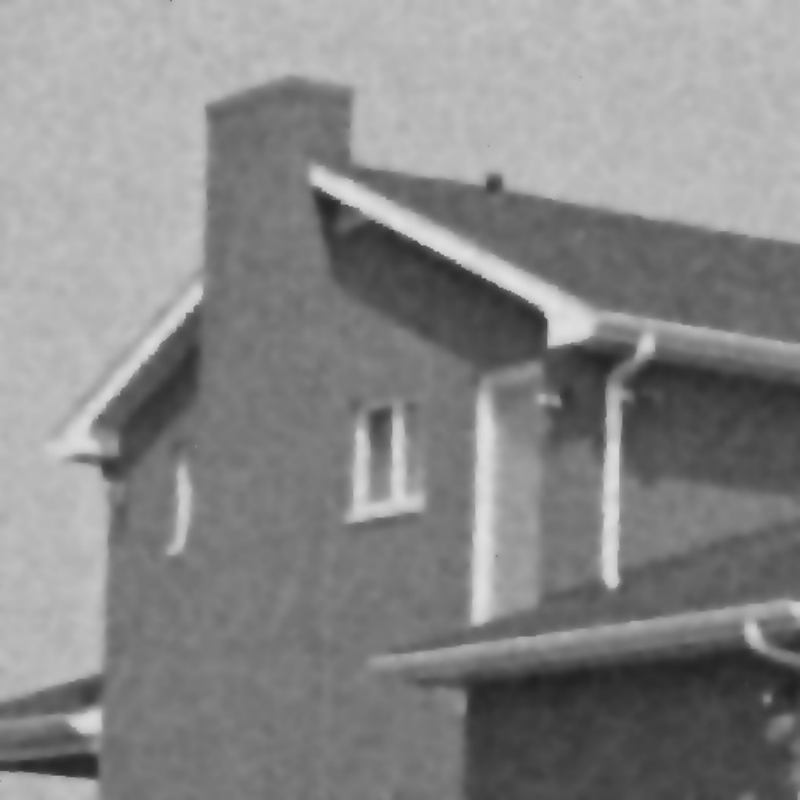} \ \
\includegraphics[page=1,width=0.13\linewidth]{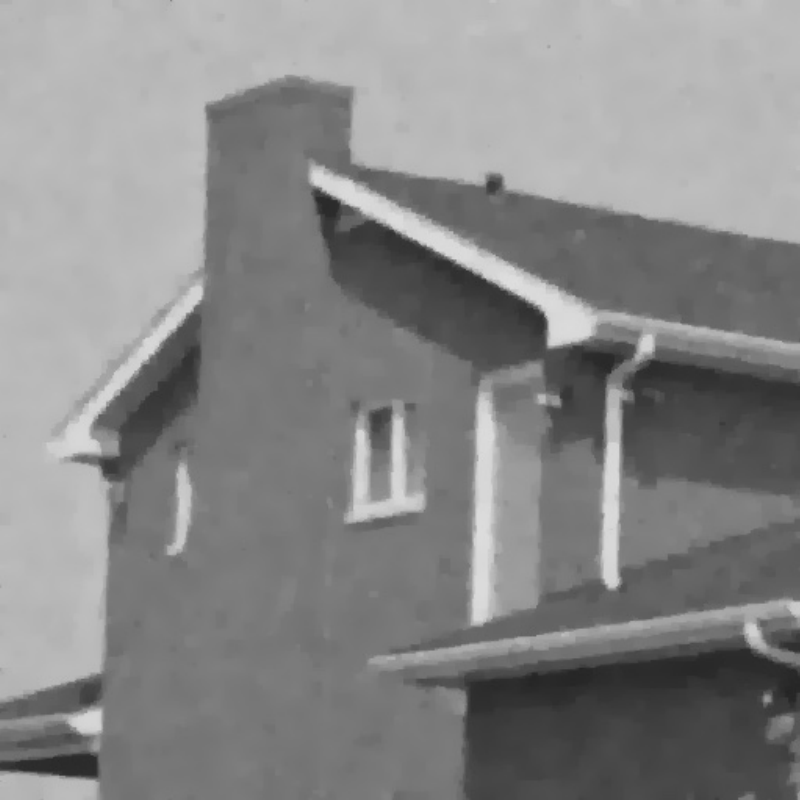}  \ \
\includegraphics[page=1,width=0.13\linewidth]{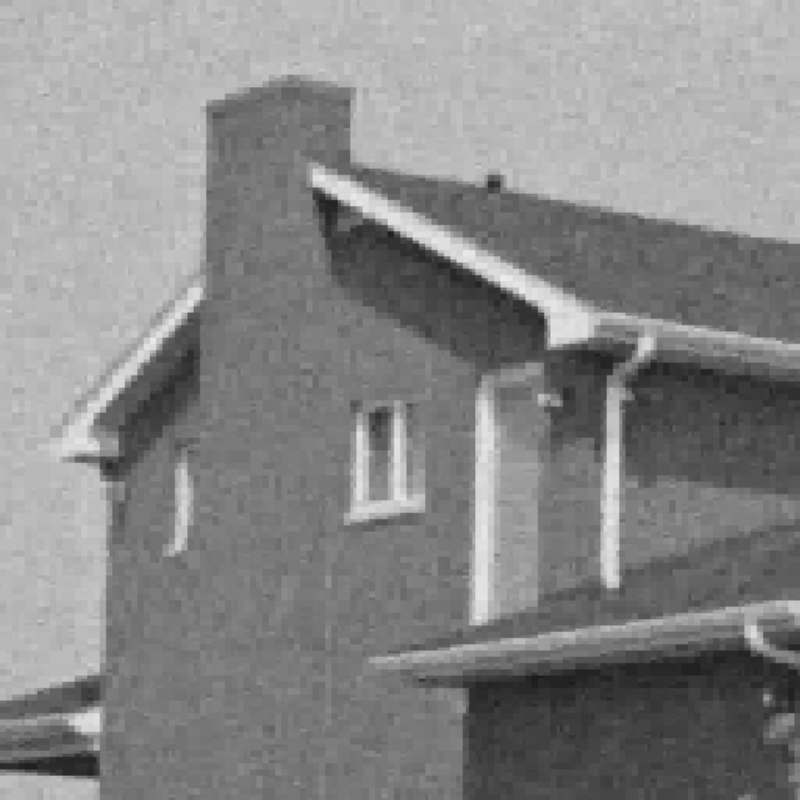} \ \
\includegraphics[page=1,width=0.13\linewidth]{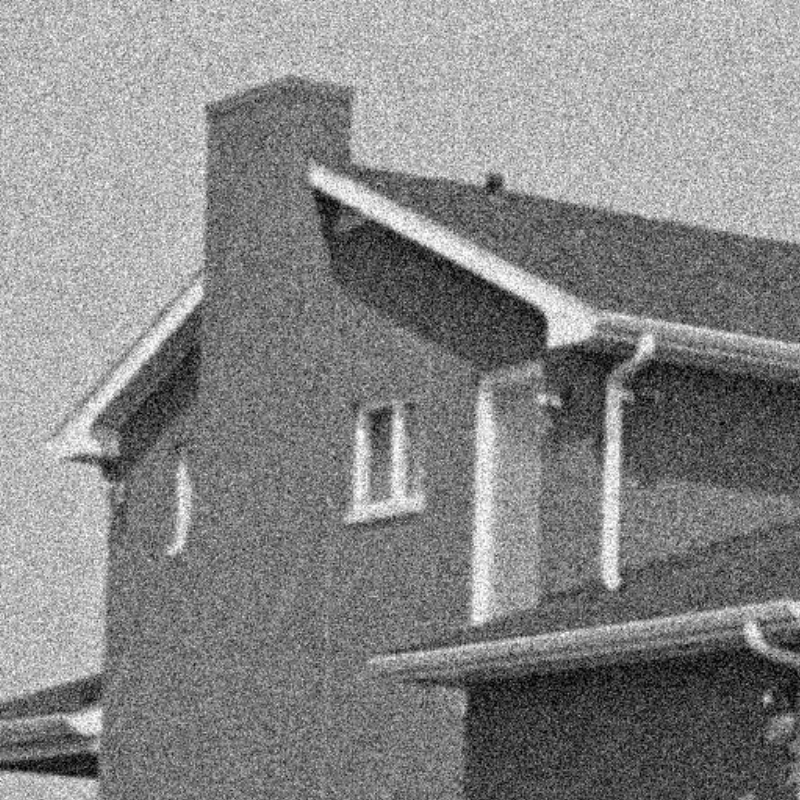} \ \
\includegraphics[page=1,width=0.13\linewidth]{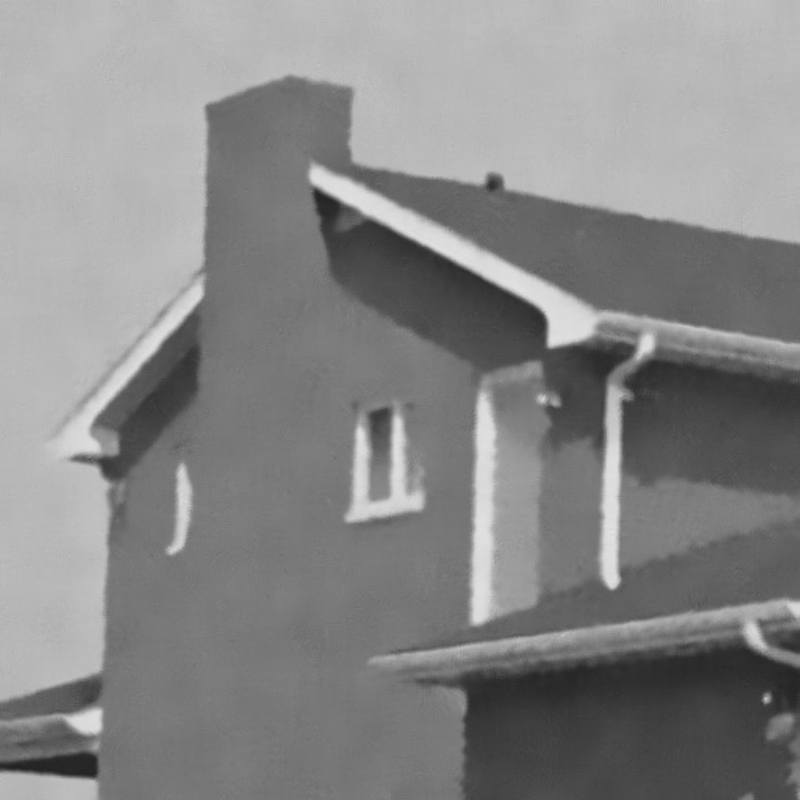}\ \
\includegraphics[page=1,width=0.13\linewidth]{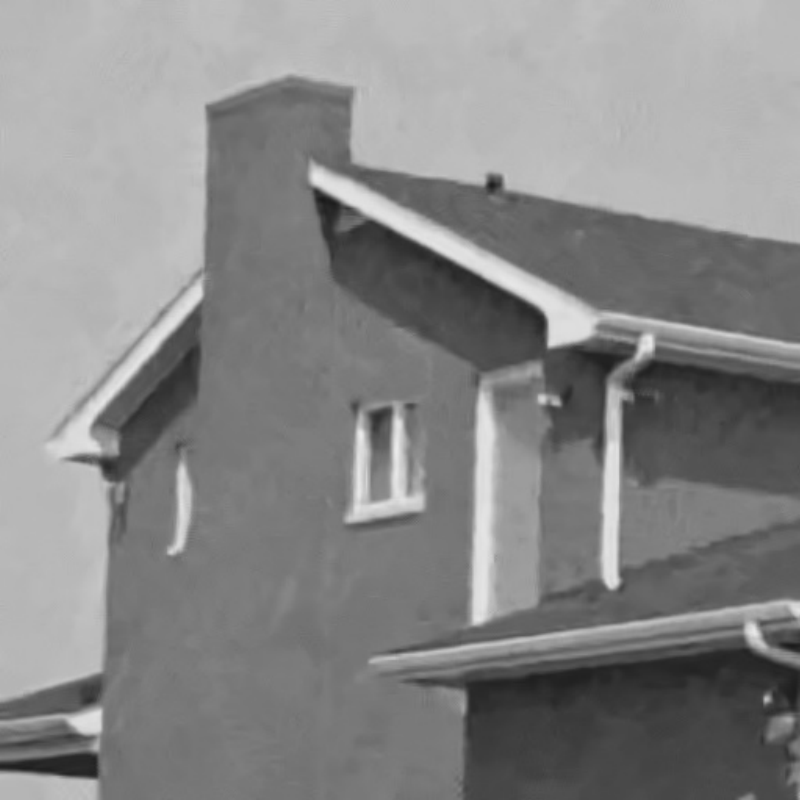}
\caption{Image 1 (House)}
\end{subfigure}
\begin{subfigure}[b]{1\textwidth}
\includegraphics[page=2,width=0.13\linewidth]{img-noi20.pdf} \ \ 
\includegraphics[page=2,width=0.13\linewidth]{img-10-nl2-20.pdf} \ \
\includegraphics[page=2,width=0.13\linewidth]{img-10-nl3-20.pdf}  \ \
\includegraphics[page=2,width=0.13\linewidth]{img-4-20.pdf} \ \
\includegraphics[page=2,width=0.13\linewidth]{img-5-20.pdf} \ \
\includegraphics[page=2,width=0.13\linewidth]{img-6-20.pdf}\ \
\includegraphics[page=2,width=0.13\linewidth]{img-3-20.pdf}
\caption{Image 2 (Peppers)}
\end{subfigure}
\begin{subfigure}[b]{1\textwidth}
\includegraphics[page=3,width=0.13\linewidth]{img-noi20.pdf} \ \ 
\includegraphics[page=3,width=0.13\linewidth]{img-10-nl2-20.pdf} \ \
\includegraphics[page=3,width=0.13\linewidth]{img-10-nl3-20.pdf}  \ \
\includegraphics[page=3,width=0.13\linewidth]{img-4-20.pdf} \ \
\includegraphics[page=3,width=0.13\linewidth]{img-5-20.pdf} \ \
\includegraphics[page=3,width=0.13\linewidth]{img-6-20.pdf}\ \
\includegraphics[page=3,width=0.13\linewidth]{img-3-20.pdf}
\caption{Image 3 (Barbara)}
\end{subfigure}
\begin{subfigure}[b]{1\textwidth}
\includegraphics[page=4,width=0.13\linewidth]{img-noi20.pdf} \ \ 
\includegraphics[page=4,width=0.13\linewidth]{img-10-nl2-20.pdf} \ \
\includegraphics[page=4,width=0.13\linewidth]{img-10-nl3-20.pdf}  \ \
\includegraphics[page=4,width=0.13\linewidth]{img-4-20.pdf} \ \
\includegraphics[page=4,width=0.13\linewidth]{img-5-20.pdf} \ \
\includegraphics[page=4,width=0.13\linewidth]{img-6-20.pdf}\ \
\includegraphics[page=4,width=0.13\linewidth]{img-3-20.pdf}
\caption{Image 4 (Baboon)}
\end{subfigure}
\begin{subfigure}[b]{1\textwidth}
\includegraphics[page=5,width=0.13\linewidth]{img-noi20.pdf} \ \ 
\includegraphics[page=5,width=0.13\linewidth]{img-10-nl2-20.pdf} \ \
\includegraphics[page=5,width=0.13\linewidth]{img-10-nl3-20.pdf}  \ \
\includegraphics[page=5,width=0.13\linewidth]{img-4-20.pdf} \ \
\includegraphics[page=5,width=0.13\linewidth]{img-5-20.pdf} \ \
\includegraphics[page=5,width=0.13\linewidth]{img-6-20.pdf}\ \
\includegraphics[page=5,width=0.13\linewidth]{img-3-20.pdf}
\caption{Image 5 (Boat)}
\end{subfigure}
\caption{\rev{Test 1. Image denoising using various algorithms. 20 \% of noise.  Image 1 (House), Image 2 (Peppers), Image 3 (Barbara), Image 4 (Baboon) and Image 5 (Boat) (from top to bottom)}}
\label{fig:t1-noi20}
\end{figure}

\begin{figure}[h!]
\centering
\begin{subfigure}[b]{0.14\textwidth}
Noised, $I_0(x)$
\end{subfigure}
\begin{subfigure}[b]{0.14\textwidth}
Ex-PM
\end{subfigure}
\begin{subfigure}[b]{0.14\textwidth}
Ex-ROF
\end{subfigure}
\begin{subfigure}[b]{0.14\textwidth}
Wvlt
\end{subfigure}
\begin{subfigure}[b]{0.14\textwidth}
TV
\end{subfigure}
\begin{subfigure}[b]{0.14\textwidth}
NLM
\end{subfigure}
\begin{subfigure}[b]{0.1\textwidth}
BM3D
\end{subfigure}

\vspace{5pt} 

\begin{subfigure}[b]{1\textwidth}
\includegraphics[page=1,width=0.13\linewidth]{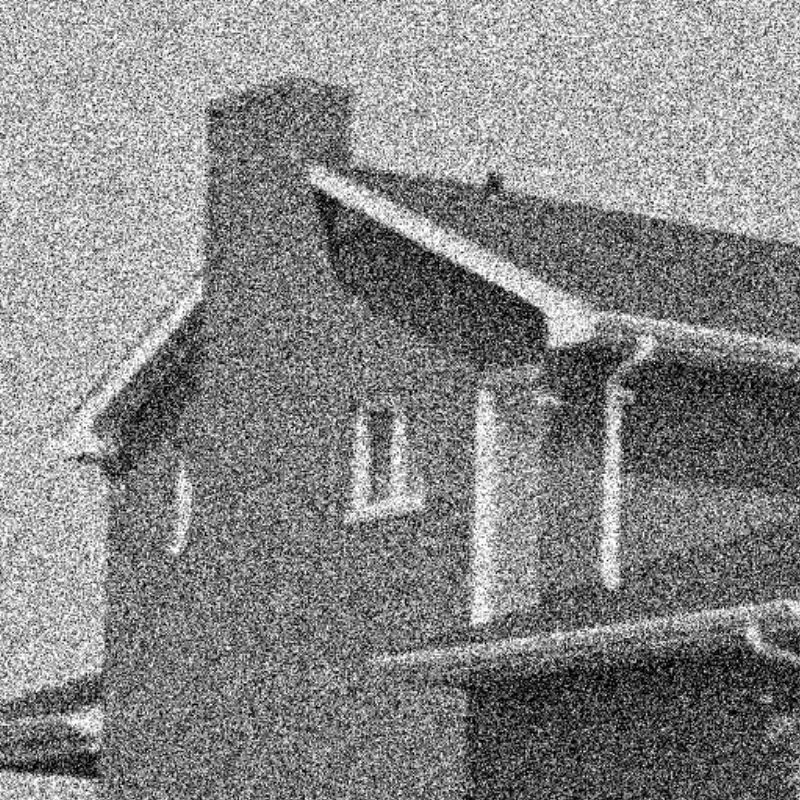} \ \ 
\includegraphics[page=1,width=0.13\linewidth]{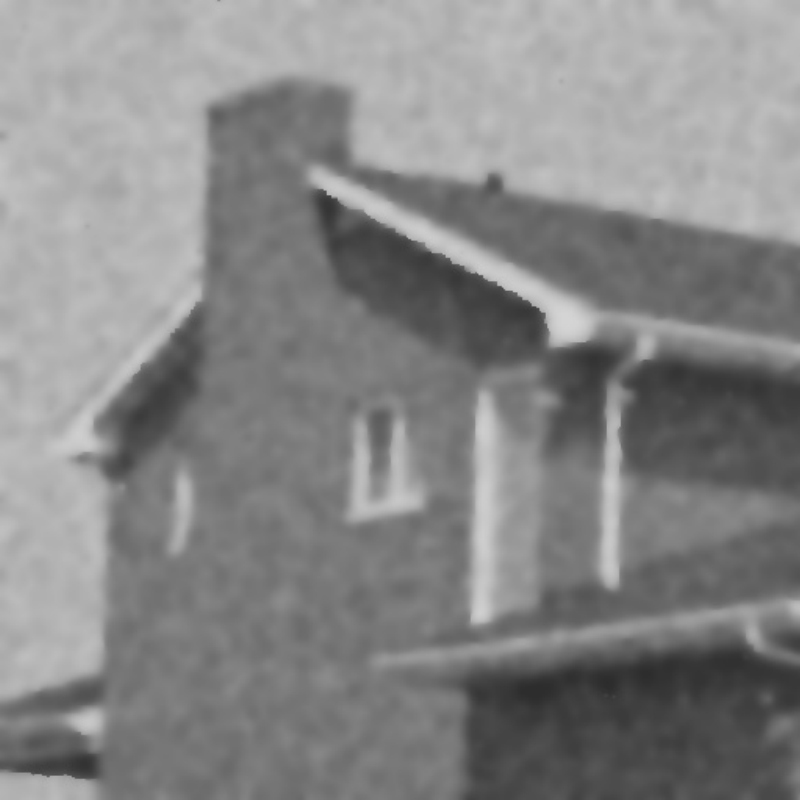} \ \
\includegraphics[page=1,width=0.13\linewidth]{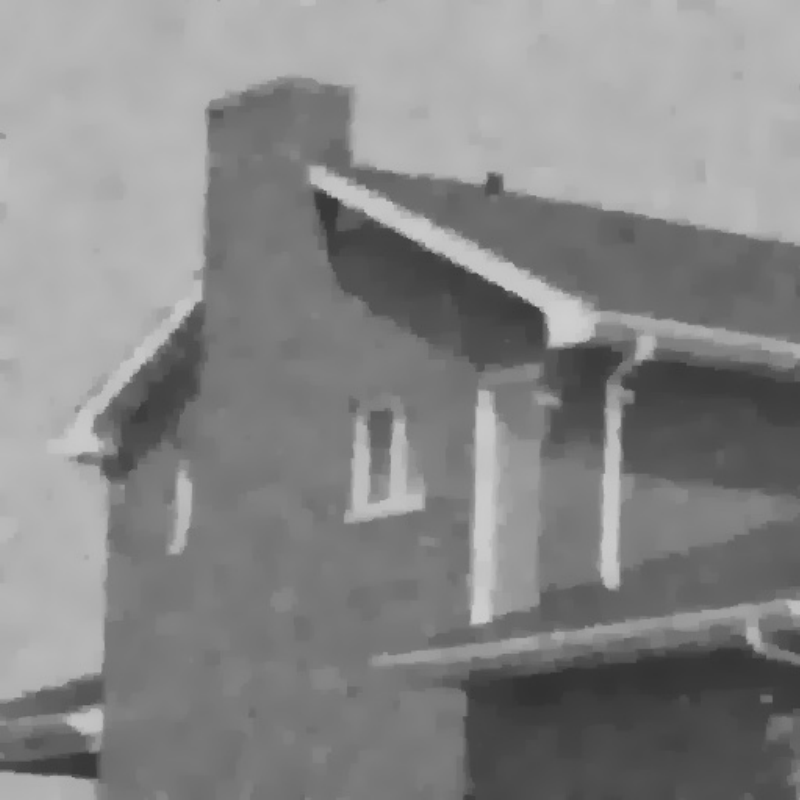}  \ \
\includegraphics[page=1,width=0.13\linewidth]{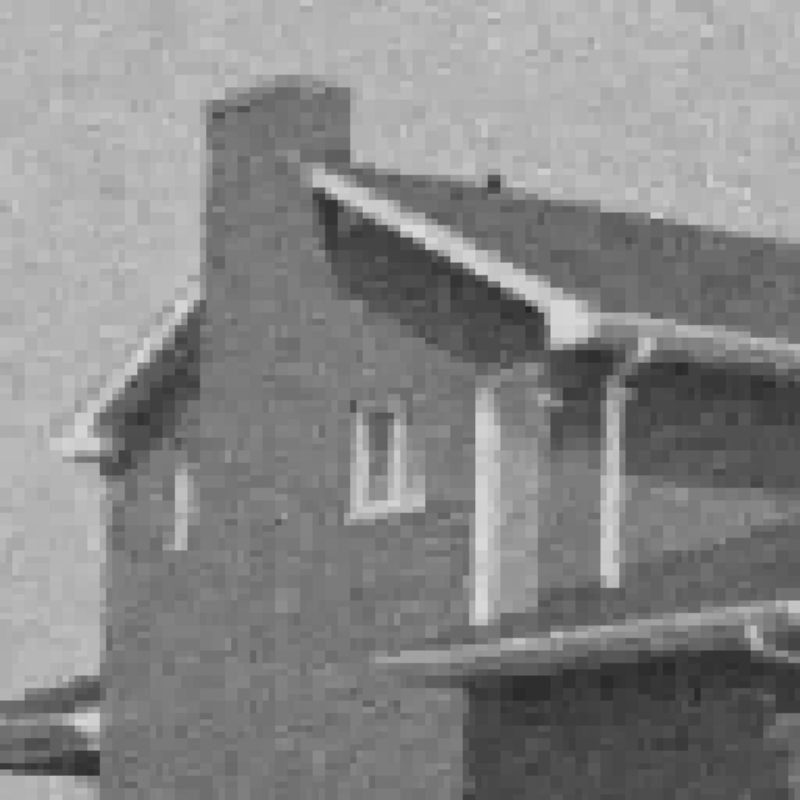} \ \
\includegraphics[page=1,width=0.13\linewidth]{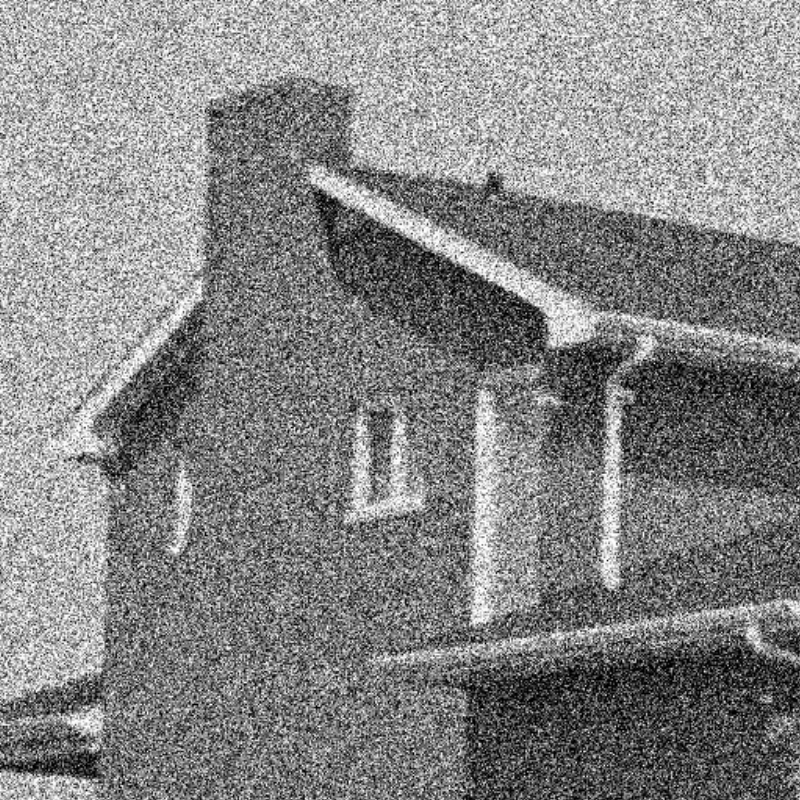} \ \
\includegraphics[page=1,width=0.13\linewidth]{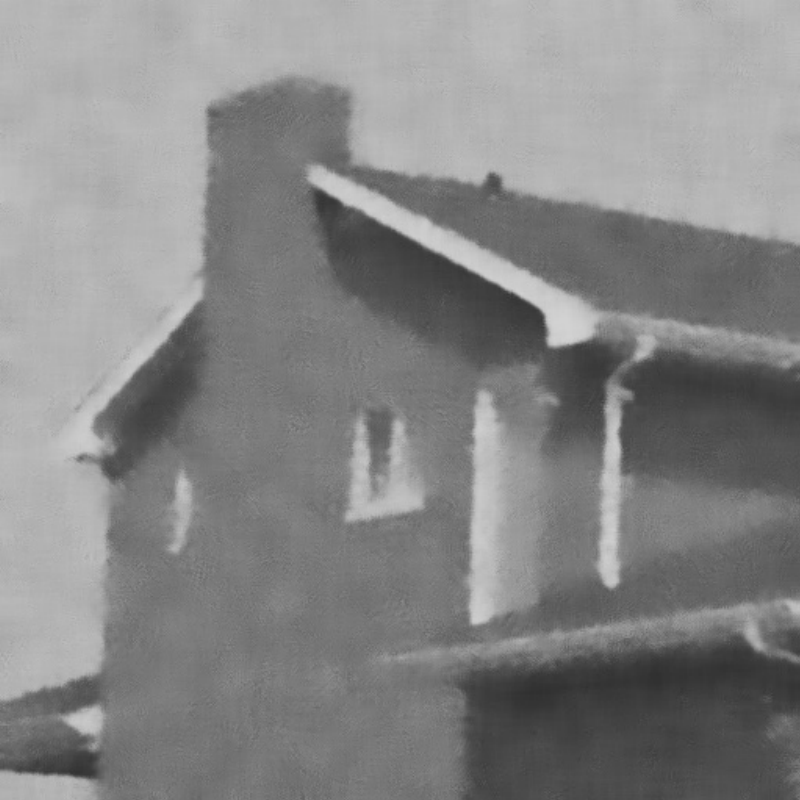}\ \
\includegraphics[page=1,width=0.13\linewidth]{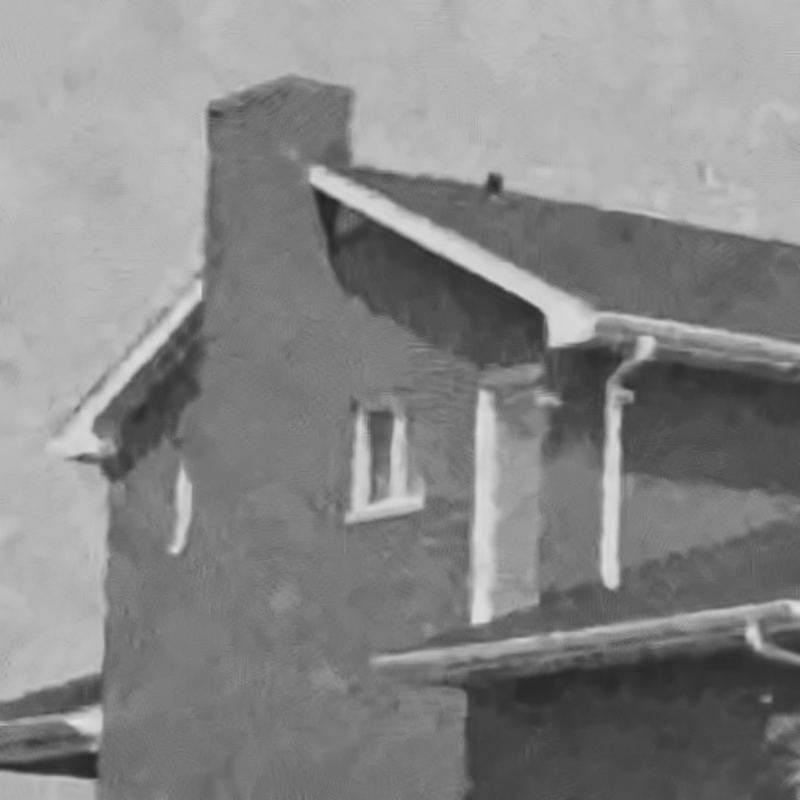}
\caption{Image 1 (House)}
\end{subfigure}
\begin{subfigure}[b]{1\textwidth}
\includegraphics[page=2,width=0.13\linewidth]{img-noi40.pdf} \ \ 
\includegraphics[page=2,width=0.13\linewidth]{img-10-nl2-40.pdf} \ \
\includegraphics[page=2,width=0.13\linewidth]{img-10-nl3-40.pdf}  \ \
\includegraphics[page=2,width=0.13\linewidth]{img-4-40.pdf} \ \
\includegraphics[page=2,width=0.13\linewidth]{img-5-40.pdf} \ \
\includegraphics[page=2,width=0.13\linewidth]{img-6-40.pdf}\ \
\includegraphics[page=2,width=0.13\linewidth]{img-3-40.pdf}
\caption{Image 2 (Peppers)}
\end{subfigure}
\begin{subfigure}[b]{1\textwidth}
\includegraphics[page=3,width=0.13\linewidth]{img-noi40.pdf} \ \ 
\includegraphics[page=3,width=0.13\linewidth]{img-10-nl2-40.pdf} \ \
\includegraphics[page=3,width=0.13\linewidth]{img-10-nl3-40.pdf}  \ \
\includegraphics[page=3,width=0.13\linewidth]{img-4-40.pdf} \ \
\includegraphics[page=3,width=0.13\linewidth]{img-5-40.pdf} \ \
\includegraphics[page=3,width=0.13\linewidth]{img-6-40.pdf}\ \
\includegraphics[page=3,width=0.13\linewidth]{img-3-40.pdf}
\caption{Image 3 (Barbara)}
\end{subfigure}
\begin{subfigure}[b]{1\textwidth}
\includegraphics[page=4,width=0.13\linewidth]{img-noi40.pdf} \ \ 
\includegraphics[page=4,width=0.13\linewidth]{img-10-nl2-40.pdf} \ \
\includegraphics[page=4,width=0.13\linewidth]{img-10-nl3-40.pdf}  \ \
\includegraphics[page=4,width=0.13\linewidth]{img-4-40.pdf} \ \
\includegraphics[page=4,width=0.13\linewidth]{img-5-40.pdf} \ \
\includegraphics[page=4,width=0.13\linewidth]{img-6-40.pdf}\ \
\includegraphics[page=4,width=0.13\linewidth]{img-3-40.pdf}
\caption{Image 4 (Baboon)}
\end{subfigure}
\begin{subfigure}[b]{1\textwidth}
\includegraphics[page=5,width=0.13\linewidth]{img-noi40.pdf} \ \ 
\includegraphics[page=5,width=0.13\linewidth]{img-10-nl2-40.pdf} \ \
\includegraphics[page=5,width=0.13\linewidth]{img-10-nl3-40.pdf}  \ \
\includegraphics[page=5,width=0.13\linewidth]{img-4-40.pdf} \ \
\includegraphics[page=5,width=0.13\linewidth]{img-5-40.pdf} \ \
\includegraphics[page=5,width=0.13\linewidth]{img-6-40.pdf}\ \
\includegraphics[page=5,width=0.13\linewidth]{img-3-40.pdf}
\caption{Image 5 (Boat)}
\end{subfigure}
\caption{\rev{Test 1. Image denoising using various algorithms. 40 \% of noise. Image 1 (House), Image 2 (Peppers), Image 3 (Barbara), Image 4 (Baboon) and Image 5 (Boat) (from top to bottom)}}
\label{fig:t1-noi40}
\end{figure}

We consider explicit scheme as a reference implementation of the proposed approach and investigate choice of nonlionearity to image denoising performance: 
\begin{itemize}
\item Ex-ROF: we solve nonlinear PDE \eqref{eq:1} with 
$c(||\nabla I||) = 1/(\sqrt{||\nabla I||^2 + \beta^2})$ using explicit scheme. We set $\beta = 10^{-5}$ and $T_{max} = 45$ and $90$ for 20 and 40 \% of noise
\item Ex-PM: we solve nonlinear PDE \eqref{eq:1} with $c(||\nabla I||) = 1/(1 + \left(\frac{||\nabla I||}{\lambda}\right)^2)$ using explicit scheme. We set $\lambda = 20$ and $T_{max} = 9$ and $18$ for 20 and 40 \% of noise
\end{itemize}
In both methods, we set a small time step size $\tau = 0.1$ to have a stability of the explicit time approximation scheme. 

We consider the following four classic methods: 
\begin{itemize}
\item Wvlt: wavelet denoising provides a sparse representation of images by coefficient thresholding those below a certain value to zero.  An adaptive thresholding method (BayesShrink) is used from the skimage restoration library \cite{chang2000adaptive}.

\item TV: total variational minimization given as follows \cite{chambolle2004algorithm}
\[
\min_I  \int_{\Omega} \left( ||\nabla I(x)|| + \frac{\mu}{2}  (I(x) - I_0(x))^2 \right)  dx,
\]
where $\mu > 0$ is the parameter with $\mu = 50$ in simulations.

\item NLM: non-local means algorithm denoises images by averaging a pixels value with those of similar neighboring pixels \cite{darbon2008fast, buades2011non, froment2014parameter}.  In numerical simulations, we estimate noise variance by estimate\_sigma function from the skimage library; the size of patches is 5, and the patch distance in pixels is 6.

\item BM3D: block-matching and 3D filtering is a 3-D block-matching algorithm \cite{dabov2007image, makinen2020collaborative}. BM3D has two main stages: hard-thresholding and Wiener filtering. Each stage involves grouping similar patches, applying collaborative filtering, and aggregating the results.  In numerical simulations, the noise variance parameter is estimated by skimage library. 
\end{itemize}
An implementation of the first three algorithms (Wvlt, TV, NLM, BM3D) is performed using skimage restoration algorithms on Python (skimage.restoration) \cite{van2014scikit}.  BM3D is bm3d python library \cite{makinen2020collaborative}. 

In Table \ref{table:t1}, we present numerical results for six considered image denoising algorithms for images with 20 \% and 40 \% of noise. The first row for a given noise level represents three initial characteristics (RRMSE, SSIM, PSNR) of the noised image. We see that initial condition $I_0(x)$ is given with SSIM $\approx$ 0.16-0.48 and PSNR $\approx$ 18-19 for 20 \% of noise. Much lower SSIM and PSNR of the initial image given for 40 \% of noise: SSIM $\approx$ 0.06-0.22 and PSNR $\approx$ 12-14. We see a smaller SSIM and PSNR for smoother images with fewer details, for example, Images 1 and 2 compared with Images 3 and 4. We compare a classical denoising technique with the proposed PDE-based denoising with two types of nonlinear coefficients: PM  (Perona-Malik model) and ROF (Rudin-Osher-Fatemi model). In this test, we consider only explicit approximation on small images to compare it with other conceptually different techniques. By comparison of the classic techniques, we see that BM3D algorithms outperform all other considered classic techniques (Wvlt, TV, NLM) for both noise levels. By comparing results for the considered PDE-based approach, we observe that the ROF nonlinear model performs better than the PM nonlinear model. However, both models strongly depend on the parameters $\lambda$ and $\beta$; therefore, we will investigate both methods in detail in the next test with varying parameter values.  By comparing Ex-ROF method with BM3D, we see that we obtain similar denoising  preformance but slightly worse in the given set of images. We will investigate it in details for other datasets with different resolution of images. 

In Figure \ref{fig:t1-noi20} and \ref{fig:t1-noi40}, we present a visual comparison of the image denoising results for various algorithms for 20 \% and 40 \% of noise, respectively. We see that Ex-ROF gives better results compared with Ex-PM with better noise removal, keeping details, and smoothing an image where needed. The best denoising performance is obtained for the BM3D method. Ex-ROF and BM3D methods are comparable and give the best results. For the classic techniques, we observe rough results for wavelets (Wvlt), bad denoising performance for the total variational minimization (TV) method that keeps noise at almost the same level, and over-smoothing for the non-local means (NLM) method on all images. 
} 

\subsection{\rev{Stability and parameters sensitivity analysis (Test 2)}}

\rev{
Similarly to the results presented above,  we take five classic greyscale images to discuss and illustrate the stability and parameter sensitivity analysis of explicit and implicit schemes. 

We compare the results of the image denoising process using the proposed PDE-based approach with varying values of the nonlinear coefficients:
\begin{itemize}
\item Ex-ROF and Im-ROF: $\beta = 10^{-5}, 10^{-3}, 10^{-1}, 1$ and $10$.
\item Ex-PM and Im-PM: $\lambda = 10, 20, 40$ and $80$.
\end{itemize}

We compare implicit and explicit time approximations with varying parameter values and time step sizes.
\begin{itemize}
\item Ex-ROF and Im-ROF: We set $T_{max} = 45$ and $90$ for 20 and 40 \% of noise. We consider $\tau = 0.1, 0.25, 1, 3, 15$ for explicit scheme and $\tau = 3, 15$ for implicit scheme. 
\item Ex-PM and Im-PM: We set $T_{max} = 9$ and $18$ for 20 and 40 \% of noise. We consider $\tau = 0.1, 0.25, 1, 3$ for explicit scheme and $\tau = 1, 3$ for implicit scheme. 
\end{itemize}

Here, we used a small time step size for the explicit scheme to preserve stability, and a large time step size is used for the unconditionally stable implicit scheme. We note that solving using an implicit scheme leads to the inversion of the problem matrix and can be computationally expensive for larger matrices; however, implicit approximation allows for much smaller time iterations.

\begin{table}[h!]
\begin{tabular}{|c | c | c | c | c | c |}
\hline
 & \tiny{RRMSE/SSIM/PSNR} & \tiny{RRMSE/SSIM/PSNR} & \tiny{RRMSE/SSIM/PSNR} & \tiny{RRMSE/SSIM/PSNR} & \tiny{RRMSE/SSIM/PSNR}\\
\hline
& 1.House & 2.Peppers & 3.Barbara & 4.Baboon & 5.Boat\\
\hline
\multicolumn{6}{|c|}{Initial noised image with 40\% of nosie}\\
\hline
 & 40.98/0.06/12.69 & 40.17/0.11/13.75 & 41.04/0.19/14.21 & 41.23/0.22/13.45 & 41.22/0.14/13.10 \\ 
\hline
\multicolumn{6}{|c|}{Ex-ROF, $\tau = 0.1$}\\
\hline
$10^{-5}$ & 6.87/0.82/28.20 & 9.59/0.77/26.19 & 14.51/0.60/23.24 & 15.26/0.48/22.09 & 11.44/0.60/24.24 \\ 
$10^{-3}$ & 6.87/0.82/28.20 & 9.59/0.77/26.19 & 14.51/0.60/23.24 & 15.26/0.48/22.08 & 11.44/0.60/24.23 \\ 
$10^{-1}$ & 6.86/0.82/28.22 & 9.57/0.77/26.21 & 14.51/0.60/23.24 & 15.27/0.48/22.08 & 11.43/0.60/24.24 \\ 
$10^{0}$ & 6.85/0.82/28.23 & 9.47/0.78/26.30 & 14.51/0.61/23.24 & 15.25/0.48/22.09 & 11.41/0.60/24.26 \\ 
$10^{1}$ & 7.21/0.79/27.79 & 9.67/0.76/26.12 & 14.73/0.59/23.11 & 15.29/0.48/22.07 & 11.54/0.59/24.16 \\ 
\hline
\multicolumn{6}{|c|}{Ex-ROF, $\tau = 0.25$}\\
\hline
$10^{-5}$ & 6.91/0.81/28.16 & 9.56/0.77/26.22 & 14.50/0.61/23.25 & 15.24/0.48/22.10 & 11.43/0.60/24.24 \\ 
$10^{-3}$ & 6.91/0.81/28.16 & 9.56/0.77/26.22 & 14.51/0.60/23.24 & 15.26/0.48/22.08 & 11.44/0.60/24.24 \\ 
$10^{-1}$ & 6.90/0.81/28.17 & 9.55/0.77/26.22 & 14.51/0.60/23.24 & 15.25/0.48/22.09 & 11.43/0.60/24.24 \\ 
$10^{0}$ & 6.85/0.82/28.23 & 9.47/0.78/26.30 & 14.51/0.61/23.24 & 15.25/0.48/22.09 & 11.41/0.60/24.26 \\ 
$10^{1}$ & 7.21/0.79/27.79 & 9.67/0.76/26.12 & 14.73/0.59/23.11 & 15.28/0.48/22.08 & 11.54/0.59/24.16 \\ 
\hline
\multicolumn{6}{|c|}{Ex-ROF, $\tau = 1$}\\
\hline
$10^{-5}$ & 7.28/0.73/27.70 & 9.65/0.73/26.14 & 14.60/0.58/23.18 & 15.21/0.48/22.12 & 11.52/0.59/24.17 \\ 
$10^{-3}$ & 7.28/0.73/27.70 & 9.65/0.73/26.14 & 14.62/0.58/23.17 & 15.22/0.48/22.11 & 11.53/0.58/24.17 \\ 
$10^{-1}$ & 7.28/0.73/27.70 & 9.65/0.73/26.14 & 14.63/0.58/23.17 & 15.28/0.48/22.07 & 11.53/0.58/24.17 \\ 
$10^{0}$ & 7.21/0.75/27.79 & 9.59/0.74/26.19 & 14.61/0.59/23.18 & 15.27/0.48/22.08 & 11.49/0.59/24.20 \\ 
$10^{1}$ & 7.20/0.79/27.79 & 9.67/0.76/26.12 & 14.74/0.59/23.10 & 15.27/0.48/22.08 & 11.54/0.59/24.16 \\ 
\hline
\multicolumn{6}{|c|}{Ex-ROF, $\tau = 3$}\\
\hline
$10^{-5}$ & 8.81/0.51/26.05 & 10.79/0.56/25.17 & 15.43/0.48/22.71 & 15.78/0.46/21.79 & 12.31/0.50/23.60 \\ 
$10^{-3}$ & 8.80/0.51/26.05 & 10.79/0.56/25.17 & 15.44/0.48/22.70 & 15.78/0.46/21.79 & 12.31/0.50/23.60 \\ 
$10^{-1}$ & 8.79/0.51/26.06 & 10.78/0.56/25.18 & 15.45/0.48/22.70 & 15.79/0.46/21.79 & 12.30/0.50/23.60 \\ 
$10^{0}$ & 8.72/0.52/26.13 & 10.70/0.57/25.24 & 15.38/0.49/22.73 & 15.74/0.46/21.81 & 12.25/0.50/23.64 \\ 
$10^{1}$ & 7.20/0.77/27.79 & 9.66/0.75/26.12 & 14.78/0.59/23.08 & 15.22/0.48/22.11 & 11.55/0.59/24.15 \\ 
\hline
\multicolumn{6}{|c|}{Ex-ROF, $\tau = 15$}\\
\hline
$10^{-5}$ & 20.21/0.16/18.83 & 22.69/0.21/18.71 & 25.72/0.28/18.27 & 23.96/0.32/18.17 & 21.97/0.25/18.57 \\ 
$10^{-3}$ & 20.21/0.16/18.83 & 22.68/0.21/18.71 & 25.72/0.28/18.27 & 23.95/0.32/18.17 & 21.96/0.25/18.57 \\ 
$10^{-1}$ & 20.06/0.16/18.90 & 22.51/0.21/18.78 & 25.72/0.28/18.27 & 23.82/0.32/18.22 & 21.81/0.25/18.63 \\ 
$10^{0}$ & 19.87/0.16/18.98 & 22.30/0.21/18.86 & 25.59/0.28/18.31 & 23.63/0.33/18.29 & 21.62/0.26/18.70 \\ 
$10^{1}$ & 17.31/0.20/20.18 & 19.26/0.26/20.14 & 22.78/0.30/19.32 & 20.99/0.34/19.32 & 18.99/0.30/19.83 \\ 
\hline
\multicolumn{6}{|c|}{Im-ROF, $\tau = 3$}\\
\hline
$10^{-5}$ & 7.05/0.82/27.99 & 9.80/0.76/26.00 & 14.43/0.60/23.29 & 15.02/0.48/22.22 & 11.60/0.60/24.11 \\ 
$10^{-3}$ & 7.04/0.82/27.99 & 9.80/0.76/26.01 & 14.43/0.60/23.29 & 15.02/0.48/22.22 & 11.60/0.60/24.11 \\ 
$10^{-1}$ & 7.01/0.82/28.03 & 9.75/0.77/26.05 & 14.42/0.61/23.29 & 15.02/0.48/22.22 & 11.58/0.60/24.13 \\ 
$10^{0}$ & 7.00/0.82/28.05 & 9.65/0.77/26.14 & 14.43/0.61/23.29 & 15.02/0.48/22.22 & 11.53/0.60/24.16 \\ 
$10^{1}$ & 7.25/0.79/27.73 & 9.76/0.76/26.04 & 14.65/0.59/23.16 & 15.12/0.49/22.17 & 11.59/0.59/24.12 \\ 
\hline
\multicolumn{6}{|c|}{Im-ROF $\tau = 15$}\\
\hline
$10^{-5}$ & 7.45/0.81/27.50 & 10.27/0.76/25.59 & 14.59/0.60/23.19 & 15.14/0.47/22.15 & 11.96/0.58/23.85 \\ 
$10^{-3}$ & 7.45/0.81/27.50 & 10.27/0.76/25.59 & 14.59/0.60/23.19 & 15.14/0.47/22.15 & 11.96/0.58/23.85 \\ 
$10^{-1}$ & 7.45/0.81/27.50 & 10.26/0.76/25.60 & 14.58/0.60/23.20 & 15.14/0.47/22.16 & 11.96/0.59/23.85 \\ 
$10^{0}$ & 7.44/0.81/27.51 & 10.22/0.76/25.64 & 14.56/0.60/23.21 & 15.12/0.47/22.16 & 11.93/0.59/23.87 \\ 
$10^{1}$ & 7.41/0.79/27.55 & 10.10/0.75/25.74 & 14.56/0.59/23.21 & 15.04/0.48/22.21 & 11.87/0.58/23.91 \\ 
\hline
\end{tabular}
\caption{\rev{Test 2. Image denoising using various algorithms. 40 \% of noise. Ex-ROF and Im-ROF methods}}
\label{table:t2-rof}
\end{table}

\begin{table}[h!]
\begin{tabular}{|c | c | c | c | c | c |}
\hline
 & \tiny{RRMSE/SSIM/PSNR} & \tiny{RRMSE/SSIM/PSNR} & \tiny{RRMSE/SSIM/PSNR} & \tiny{RRMSE/SSIM/PSNR} & \tiny{RRMSE/SSIM/PSNR}\\
\hline
& 1.House & 2.Peppers & 3.Barbara & 4.Baboon & 5.Boat\\
\hline
\multicolumn{6}{|c|}{Initial noised image with 40\% of nosie}\\
\hline
 & 40.98/0.06/12.69 & 40.17/0.11/13.75 & 41.04/0.19/14.21 & 41.23/0.22/13.45 & 41.22/0.14/13.10 \\ 
\hline
\multicolumn{6}{|c|}{Ex-PM, $\tau = 0.1$}\\
\hline
$10$ & 8.26/0.78/26.61 & 10.88/0.73/25.09 & 15.42/0.57/22.71 & 16.87/0.41/21.21 & 12.69/0.54/23.33 \\ 
$20$ & 7.47/0.79/27.48 & 9.91/0.76/25.90 & 14.87/0.59/23.03 & 15.84/0.45/21.76 & 11.88/0.58/23.91 \\ 
$40$ & 7.46/0.78/27.49 & 9.91/0.74/25.91 & 14.90/0.58/23.01 & 15.39/0.48/22.01 & 11.71/0.58/24.03 \\ 
$80$ & 7.57/0.77/27.36 & 10.37/0.73/25.51 & 15.09/0.57/22.90 & 15.45/0.48/21.98 & 11.98/0.57/23.83 \\ 
\hline
\multicolumn{6}{|c|}{Ex-PM, $\tau = 0.25$}\\
\hline
$10$ & 8.26/0.78/26.61 & 10.88/0.73/25.09 & 15.42/0.57/22.71 & 16.83/0.41/21.24 & 12.69/0.54/23.33 \\ 
$20$ & 7.48/0.79/27.47 & 9.92/0.76/25.89 & 14.89/0.59/23.02 & 15.81/0.45/21.78 & 11.89/0.58/23.90 \\ 
$40$ & 7.45/0.78/27.50 & 9.90/0.74/25.91 & 14.95/0.58/22.98 & 15.37/0.47/22.02 & 11.73/0.58/24.02 \\ 
$80$ & 7.56/0.77/27.37 & 10.37/0.73/25.51 & 15.13/0.57/22.88 & 15.49/0.47/21.96 & 12.00/0.57/23.82 \\ 
\hline
\multicolumn{6}{|c|}{Ex-PM, $\tau = 1$}\\
\hline
$10$ & 10.37/0.41/24.63 & 12.55/0.46/23.86 & 16.98/0.41/21.88 & 17.82/0.39/20.74 & 13.85/0.41/22.57 \\ 
$20$ & 14.71/0.24/21.59 & 16.58/0.30/21.43 & 20.24/0.33/20.35 & 20.04/0.36/19.72 & 17.08/0.32/20.75 \\ 
$40$ & 21.60/0.15/18.26 & 24.02/0.20/18.22 & 27.44/0.26/17.71 & 24.90/0.30/17.83 & 23.14/0.24/18.11 \\ 
$80$ & 38.39/0.06/13.26 & 40.17/0.11/13.75 & 41.04/0.19/14.21 & 41.23/0.22/13.45 & 41.01/0.14/13.15 \\ 
\hline
\multicolumn{6}{|c|}{Ex-PM, $\tau = 3$}\\
\hline
$10$ & 15.69/0.21/21.03 & 18.12/0.26/20.67 & 22.48/0.28/19.44 & 22.16/0.32/18.85 & 18.63/0.28/20.00 \\ 
$20$ & 28.80/0.10/15.76 & 32.43/0.14/15.61 & 34.58/0.20/15.70 & 32.60/0.24/15.49 & 31.23/0.18/15.51 \\ 
$40$ & 40.98/0.06/12.69 & 40.17/0.11/13.75 & 41.04/0.19/14.21 & 41.23/0.22/13.45 & 41.22/0.14/13.10 \\ 
$80$ & 40.98/0.06/12.69 & 40.17/0.11/13.75 & 41.04/0.19/14.21 & 41.23/0.22/13.45 & 41.22/0.14/13.10 \\ 
\hline
\multicolumn{6}{|c|}{Im-PM, $\tau = 1$}\\
\hline
$10$ & 8.29/0.78/26.57 & 10.89/0.73/25.09 & 15.32/0.57/22.77 & 16.31/0.42/21.51 & 12.68/0.54/23.34 \\ 
$20$ & 7.59/0.79/27.34 & 10.14/0.75/25.71 & 14.75/0.59/23.10 & 15.41/0.46/22.00 & 12.02/0.57/23.80 \\ 
$40$ & 7.57/0.78/27.37 & 10.26/0.74/25.61 & 14.66/0.58/23.15 & 15.13/0.47/22.16 & 12.00/0.58/23.82 \\ 
$80$ & 7.68/0.77/27.24 & 10.49/0.73/25.41 & 14.82/0.57/23.06 & 15.31/0.48/22.06 & 12.04/0.57/23.79 \\ 
\hline
\multicolumn{6}{|c|}{Im-PM, $\tau = 3$}\\
\hline
$10$ & 8.31/0.77/26.55 & 10.92/0.73/25.06 & 15.23/0.57/22.82 & 15.90/0.43/21.73 & 12.69/0.55/23.33 \\ 
$20$ & 7.79/0.78/27.12 & 10.46/0.75/25.43 & 14.87/0.58/23.03 & 15.41/0.45/22.00 & 12.23/0.57/23.65 \\ 
$40$ & 7.76/0.78/27.14 & 10.84/0.73/25.12 & 14.64/0.58/23.16 & 15.14/0.47/22.16 & 12.42/0.57/23.52 \\ 
$80$ & 7.84/0.77/27.06 & 11.11/0.73/24.91 & 15.21/0.58/22.83 & 15.77/0.46/21.80 & 12.63/0.57/23.37 \\ 
\hline
\end{tabular}
\caption{\rev{Test 2. Image denoising using proposed approach with implicit and explicit time approximation and varying time step size. 40 \% of noise. Ex-PM and Im-PM methods}}
\label{table:t2-pm}
\end{table}

In Tables \ref{table:t2-rof} and \ref{table:t2-pm}, we present numerical results for varying coefficient and time step size values for explicit and implicit time approximation schemes of the PDE-based approach. In this test, we consider 40\% of noise in the initial image. 
For the Ex-ROF method with a stable time step size, the smaller parameter $\beta$ is preferable. However, we obtain almost similar results except for very large $\beta = 10$. For the Ex-PM model, we observe that $\lambda=20$ gives better denoising results. We also mention that ROF and PM coefficient nonlinearity also affect the final time that we use in simulation, where faster denoising with smaller $T_{max}$ is observed for the PM model. Moreover it is effects on the choice of the time step size inimplicit approximation. For the case with a minimum 3 time iterations in the implicit scheme, we obtain maximum time step size $\tau = 3$ in Im-PM and  $\tau = 15$ in Im-ROF. 
For both nonlinear coefficients in explicit time approximation, we observe stable solutions for small time step size $\tau = 0.1$ and $0.25$ with the same denoising performance. 
We should mention that we use first-order time approximation schemes, so the time step size may affect the final solution. 
For the bigger time step size, we observe a high influence on the final denoising characteristics that make all three characteristics, RRMSE, SSIM, and PSNT, perform worse due to the instability of the explicit scheme. However, we see a little better stability for the ROF nonlinear diffusion coefficient for $\tau = 1$ compared with the PM nonlinearity with the same $\tau = 1$ that becomes unstable and does not remove the noise (even produces more oscillations due to the instability of the scheme). We observe poor results for bigger time step size $\tau = 3$ or $15$. 
In contrast with the unstable solution of the explicit scheme, we observe the unconditional stability of implicit time approximation. We see that with a large time step size $\tau = 1$ for Im-PM or $3$ for Im-ROF, we may perform the same denoising performance with a huge reduction of time step iterations. We see that we can perform a simulation with extremely large time step size compared with a final time $T_{max}$, which gives only 3-time iterations. Once again, we should consider the influence of time scheme accuracy, which is the first order in a given time approximation scheme. It would be interesting to consider the applicability of the high-order time approximation to the image-denoising process in future works. 
}

\subsection{\rev{Multiscale approximation on the coarse grid (Test 3)}}

\rev{
Next, we consider the performance of the multiscale approximation on the coarse grid. We consider ROF nonlinear coefficient with $\beta = 10^{-5}$ and implicit time approximation scheme. We note that the proposed multiscale methods are implemented only for implicit time approximation schemes in order to reduce the size of the matrix that is inverted in each iteration to perform faster calculations.

\begin{table}[h!]
\begin{tabular}{| c | c | c | c | c | c | }
\hline
 & \tiny{RRMSE/SSIM/PSNR} & \tiny{RRMSE/SSIM/PSNR} & \tiny{RRMSE/SSIM/PSNR} & \tiny{RRMSE/SSIM/PSNR} & \tiny{RRMSE/SSIM/PSNR}\\
\hline
& 1.House & 2.Peppers & 3.Barbara & 4.Baboon & 5.Boat\\
\hline
\multicolumn{6}{|c|}{Initial noised image with 20\% of nosie}\\
\hline
$I_0(x)$ & 21.29/0.16/18.38 & 20.33/0.26/19.66 & 20.44/0.41/20.26 & 20.51/0.48/19.52 & 20.77/0.33/19.06 \\ 
\hline
\multicolumn{6}{|c|}{Fine-scale method}\\
\hline
Im-ROF & 4.37/0.86/32.13 & 6.13/0.84/30.08 & 10.72/0.72/25.87 & 11.29/0.67/24.71 & 7.95/0.72/27.40 \\ 
Ex-ROF & 4.23/0.86/32.42 & 5.91/0.85/30.39 & 10.82/0.73/25.79 & 11.45/0.67/24.58 & 7.80/0.73/27.56 \\ 
\hline
\multicolumn{6}{|c|}{Multiscale method, $8 \times 8$ coarse grid, $N_h^{K_i} = 16^2$}\\
\hline
Ms$^0$ & 6.49/0.81/28.70 & 11.10/0.76/24.92 & 15.49/0.58/22.67 & 16.47/0.48/21.42 & 12.98/0.56/23.13 \\ 
Ms & 6.49/0.81/28.70 & 11.10/0.76/24.92 & 15.49/0.58/22.67 & 16.46/0.48/21.43 & 12.98/0.56/23.13 \\ 
\hline
AMs1 & 10.98/0.72/24.13 & 15.92/0.63/21.79 & 18.65/0.47/21.06 & 18.26/0.48/20.53 & 16.60/0.46/21.00 \\ 
AMs2 & 8.91/0.76/25.94 & 13.05/0.69/23.52 & 16.53/0.53/22.11 & 17.48/0.48/20.90 & 15.56/0.48/21.57 \\ 
AMs4 & 6.03/0.80/29.34 & 10.52/0.73/25.38 & 15.30/0.56/22.78 & 16.62/0.48/21.35 & 13.15/0.53/23.02 \\ 
AMs8 & 5.15/0.83/30.72 & 8.92/0.77/26.82 & 14.54/0.60/23.22 & 15.81/0.48/21.78 & 11.89/0.57/23.90 \\ 
AMs12 & 4.40/0.85/32.08 & 7.71/0.80/28.09 & 14.00/0.63/23.55 & 15.11/0.48/22.17 & 10.78/0.63/24.75 \\ 
AMs16 & 4.27/0.86/32.33 & 7.26/0.82/28.61 & 13.74/0.64/23.71 & 14.77/0.49/22.37 & 10.14/0.65/25.28 \\ 
\hline
GMs1 & 10.98/0.72/24.13 & 15.92/0.63/21.79 & 18.65/0.47/21.06 & 18.26/0.48/20.53 & 16.60/0.46/21.00 \\ 
GMs2 & 6.50/0.80/28.68 & 11.16/0.73/24.87 & 15.70/0.56/22.55 & 16.92/0.48/21.19 & 13.74/0.52/22.64 \\ 
GMs4 & 5.14/0.84/30.73 & 9.16/0.78/26.58 & 14.61/0.60/23.18 & 16.04/0.48/21.65 & 12.10/0.58/23.74 \\ 
GMs8 & 4.35/0.86/32.17 & 7.87/0.81/27.90 & 13.92/0.63/23.60 & 15.13/0.48/22.16 & 10.63/0.63/24.87 \\ 
GMs12 & 4.14/0.86/32.61 & 7.25/0.82/28.61 & 13.65/0.65/23.77 & 14.69/0.49/22.42 & 9.96/0.66/25.44 \\ 
GMs16 & 4.03/0.86/32.83 & 6.90/0.83/29.05 & 13.48/0.65/23.88 & 14.35/0.51/22.62 & 9.56/0.67/25.80 \\ 
\hline
\multicolumn{6}{|c|}{Multiscale method, $16 \times 16$ coarse grid, $N_h^{K_i} = 32^2$}\\
\hline
Ms$^0$ & 6.49/0.81/28.70 & 11.10/0.76/24.92 & 15.49/0.58/22.67 & 16.47/0.48/21.42 & 12.98/0.56/23.13 \\ 
Ms & 6.49/0.81/28.70 & 11.10/0.76/24.92 & 15.49/0.58/22.67 & 16.46/0.48/21.43 & 12.98/0.56/23.13 \\ 
\hline
AMs1 & 15.24/0.69/21.29 & 20.33/0.56/19.66 & 20.44/0.41/20.26 & 19.56/0.48/19.93 & 20.77/0.41/19.06 \\ 
AMs2 & 13.99/0.70/22.03 & 18.86/0.60/20.32 & 20.41/0.44/20.28 & 18.85/0.48/20.25 & 19.96/0.42/19.40 \\ 
AMs4 & 10.45/0.72/24.56 & 15.26/0.64/22.16 & 18.28/0.47/21.24 & 18.18/0.48/20.56 & 16.30/0.46/21.16 \\ 
AMs8 & 8.56/0.75/26.29 & 13.14/0.67/23.46 & 16.83/0.51/21.95 & 17.43/0.48/20.93 & 14.82/0.48/21.99 \\ 
AMs12 & 6.49/0.79/28.69 & 11.24/0.71/24.81 & 15.80/0.54/22.50 & 16.80/0.48/21.25 & 13.60/0.52/22.73 \\ 
AMs16 & 5.93/0.81/29.49 & 10.26/0.74/25.61 & 15.24/0.56/22.81 & 16.48/0.48/21.42 & 13.09/0.53/23.07 \\ 
\hline
GMs1 & 15.24/0.69/21.29 & 20.33/0.56/19.66 & 20.44/0.41/20.26 & 19.56/0.48/19.93 & 20.77/0.41/19.06 \\ 
GMs2 & 10.65/0.75/24.39 & 15.79/0.65/21.86 & 18.77/0.49/21.01 & 18.32/0.48/20.50 & 16.97/0.46/20.81 \\ 
GMs4 & 7.57/0.79/27.37 & 12.60/0.71/23.82 & 16.36/0.54/22.20 & 17.46/0.48/20.92 & 14.94/0.50/21.92 \\ 
GMs8 & 5.73/0.82/29.78 & 10.27/0.75/25.59 & 15.12/0.58/22.88 & 16.51/0.48/21.40 & 12.98/0.55/23.14 \\ 
GMs12 & 5.14/0.84/30.72 & 9.40/0.77/26.36 & 14.68/0.60/23.14 & 16.03/0.48/21.66 & 12.21/0.57/23.67 \\ 
GMs16 & 4.83/0.84/31.27 & 8.78/0.79/26.96 & 14.40/0.61/23.30 & 15.73/0.48/21.82 & 11.71/0.59/24.03 \\ 
\hline
\end{tabular}
\caption{\rev{Test 3. Image denoising using multiscale method. 20\% of noise}}
\label{table:t3-1}
\end{table}

\begin{table}[h!]
\begin{tabular}{| c | c | c | c | c | c | }
\hline
 & \tiny{RRMSE/SSIM/PSNR} & \tiny{RRMSE/SSIM/PSNR} & \tiny{RRMSE/SSIM/PSNR} & \tiny{RRMSE/SSIM/PSNR} & \tiny{RRMSE/SSIM/PSNR}\\
\hline
& 1.House & 2.Peppers & 3.Barbara & 4.Baboon & 5.Boat\\
\hline
\multicolumn{6}{|c|}{Initial noised image with 40\% of noise}\\
\hline
$I_0(x)$ & 40.98/0.06/12.69 & 40.17/0.11/13.75 & 41.04/0.19/14.21 & 41.23/0.22/13.45 & 41.22/0.14/13.10 \\ 
\hline
\multicolumn{6}{|c|}{Fine-scale method}\\
\hline
Ex-ROF & 6.91/0.81/28.15 & 9.57/0.77/26.21 & 14.50/0.60/23.25 & 15.11/0.48/22.17 & 11.45/0.60/24.23 \\ 
Im-ROF & 7.05/0.82/27.99 & 9.80/0.76/26.00 & 14.43/0.60/23.29 & 15.02/0.48/22.22 & 11.60/0.60/24.11 \\ 
\hline
\multicolumn{6}{|c|}{Multiscale method, $8 \times 8$ coarse grid, $N_h^{K_i} = 16^2$}\\
\hline
Ms$^0$ & 7.98/0.77/26.91 & 12.36/0.72/23.98 & 16.11/0.56/22.33 & 16.87/0.35/21.22 & 14.01/0.53/22.48 \\ 
Ms & 7.97/0.78/26.91 & 12.36/0.73/23.98 & 16.11/0.56/22.33 & 16.86/0.35/21.22 & 14.01/0.53/22.48 \\ 
\hline
AMs1 & 11.54/0.72/23.70 & 16.43/0.63/21.51 & 18.93/0.47/20.93 & 18.44/0.28/20.44 & 17.07/0.45/20.76 \\ 
AMs2 & 9.80/0.75/25.12 & 13.86/0.68/22.99 & 16.99/0.52/21.87 & 17.79/0.31/20.75 & 16.22/0.46/21.20 \\ 
AMs4 & 7.60/0.79/27.33 & 11.88/0.71/24.33 & 16.01/0.54/22.38 & 17.12/0.33/21.09 & 14.22/0.50/22.35 \\ 
AMs8 & 7.14/0.81/27.87 & 10.86/0.73/25.11 & 15.52/0.56/22.66 & 16.60/0.37/21.35 & 13.39/0.53/22.87 \\ 
AMs12 & 6.70/0.82/28.42 & 10.08/0.76/25.75 & 15.16/0.58/22.86 & 16.18/0.39/21.58 & 12.62/0.56/23.38 \\ 
AMs16 & 6.68/0.82/28.44 & 9.86/0.77/25.95 & 15.00/0.59/22.95 & 16.03/0.41/21.66 & 12.24/0.58/23.65 \\ 
\hline
GMs1 & 11.54/0.72/23.70 & 16.43/0.63/21.51 & 18.93/0.47/20.93 & 18.44/0.28/20.44 & 17.07/0.45/20.76 \\ 
GMs2 & 7.82/0.79/27.08 & 12.43/0.71/23.94 & 16.37/0.54/22.19 & 17.41/0.32/20.94 & 14.58/0.50/22.13 \\ 
GMs4 & 7.10/0.81/27.92 & 11.15/0.74/24.88 & 15.52/0.57/22.66 & 16.83/0.36/21.24 & 13.41/0.53/22.85 \\ 
GMs8 & 6.72/0.82/28.39 & 10.30/0.76/25.57 & 15.07/0.59/22.91 & 16.27/0.40/21.53 & 12.46/0.57/23.49 \\ 
GMs12 & 6.68/0.83/28.45 & 9.96/0.77/25.86 & 14.93/0.60/22.99 & 16.04/0.41/21.66 & 12.12/0.58/23.73 \\ 
GMs16 & 6.65/0.83/28.48 & 9.78/0.77/26.02 & 14.86/0.60/23.03 & 15.85/0.42/21.75 & 11.95/0.59/23.85 \\ 
\hline
\multicolumn{6}{|c|}{Multiscale method, $16 \times 16$ coarse grid, $N_h^{K_i} = 32^2$}\\
\hline
Ms$^0$ & 7.98/0.77/26.91 & 12.36/0.72/23.98 & 16.11/0.56/22.33 & 16.87/0.35/21.22 & 14.01/0.53/22.48 \\ 
Ms & 7.97/0.78/26.91 & 12.36/0.73/23.98 & 16.11/0.56/22.33 & 16.86/0.35/21.22 & 14.01/0.53/22.48 \\ 
\hline
AMs1 & 15.49/0.70/21.14 & 23.24/0.56/18.50 & 23.75/0.41/18.96 & 19.63/0.26/19.90 & 21.10/0.41/18.92 \\ 
AMs2 & 14.32/0.70/21.82 & 19.21/0.60/20.16 & 20.60/0.44/20.20 & 18.97/0.27/20.20 & 20.27/0.42/19.27 \\ 
AMs4 & 11.06/0.72/24.07 & 15.81/0.63/21.85 & 18.57/0.47/21.10 & 18.37/0.28/20.47 & 16.79/0.45/20.90 \\ 
AMs8 & 9.56/0.74/25.34 & 13.97/0.66/22.92 & 17.30/0.50/21.71 & 17.75/0.30/20.77 & 15.54/0.47/21.57 \\ 
AMs12 & 7.93/0.78/26.96 & 12.40/0.69/23.96 & 16.43/0.53/22.16 & 17.25/0.32/21.02 & 14.54/0.49/22.15 \\ 
AMs16 & 7.54/0.79/27.39 & 11.68/0.71/24.48 & 15.99/0.54/22.40 & 17.01/0.34/21.14 & 14.16/0.50/22.38 \\ 
\hline
GMs1 & 15.49/0.70/21.14 & 23.24/0.56/18.50 & 23.75/0.41/18.96 & 19.63/0.26/19.90 & 21.10/0.41/18.92 \\ 
GMs2 & 11.18/0.74/23.97 & 16.63/0.64/21.41 & 19.25/0.48/20.79 & 18.59/0.28/20.37 & 17.28/0.45/20.65 \\ 
GMs4 & 8.75/0.78/26.11 & 13.68/0.68/23.10 & 16.94/0.53/21.90 & 17.80/0.31/20.75 & 15.43/0.48/21.64 \\ 
GMs8 & 7.37/0.80/27.60 & 11.80/0.73/24.39 & 15.89/0.56/22.45 & 17.15/0.34/21.07 & 14.08/0.52/22.43 \\ 
GMs12 & 7.06/0.81/27.97 & 11.22/0.74/24.83 & 15.55/0.57/22.64 & 16.84/0.36/21.23 & 13.58/0.53/22.75 \\ 
GMs16 & 6.94/0.82/28.11 & 10.88/0.75/25.10 & 15.36/0.58/22.75 & 16.60/0.37/21.35 & 13.21/0.54/22.98 \\ 
\hline
\end{tabular}
\caption{\rev{Test 3. Image denoising using multiscale method. 40\% of noise}}
\label{table:t3-2}
\end{table}

\begin{figure}[h!]
\centering
\begin{subfigure}[b]{0.18\textwidth}
\includegraphics[page=1,width=1\linewidth]{img-noi40.pdf} 
\caption{Noised, $I^0(x)$}
\end{subfigure} \ 
\begin{subfigure}[b]{0.18\textwidth}
\includegraphics[page=1,width=1\linewidth]{img-10-nl3-40.pdf} 
\caption{Ex-ROF}
\end{subfigure} \ 
\begin{subfigure}[b]{0.18\textwidth}
\includegraphics[page=1,width=1\linewidth]{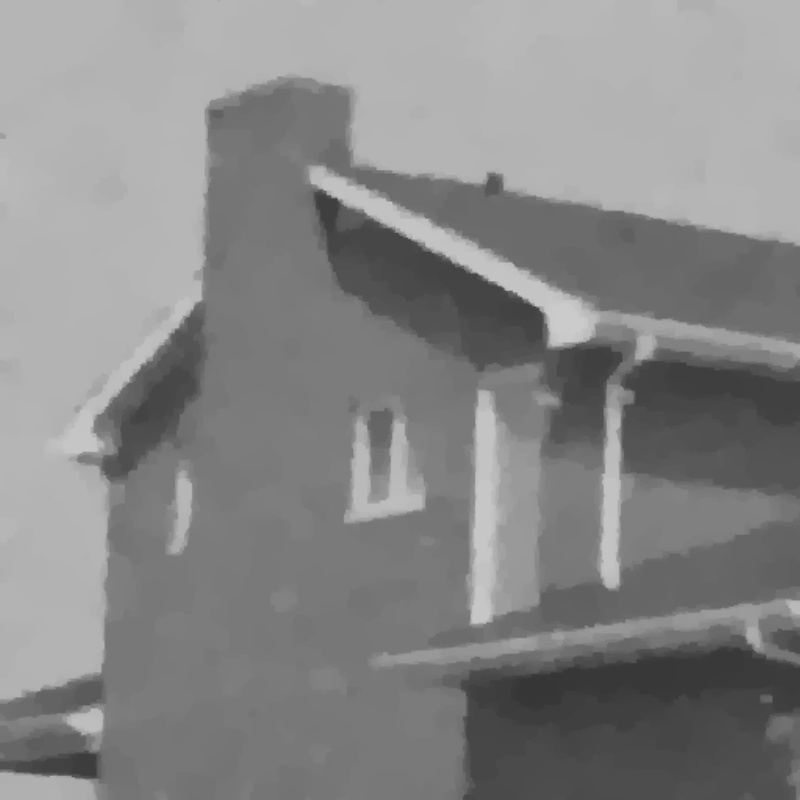} 
\caption{Im-ROF}
\end{subfigure} \ 
\begin{subfigure}[b]{0.18\textwidth}
\includegraphics[page=1,width=1\linewidth]{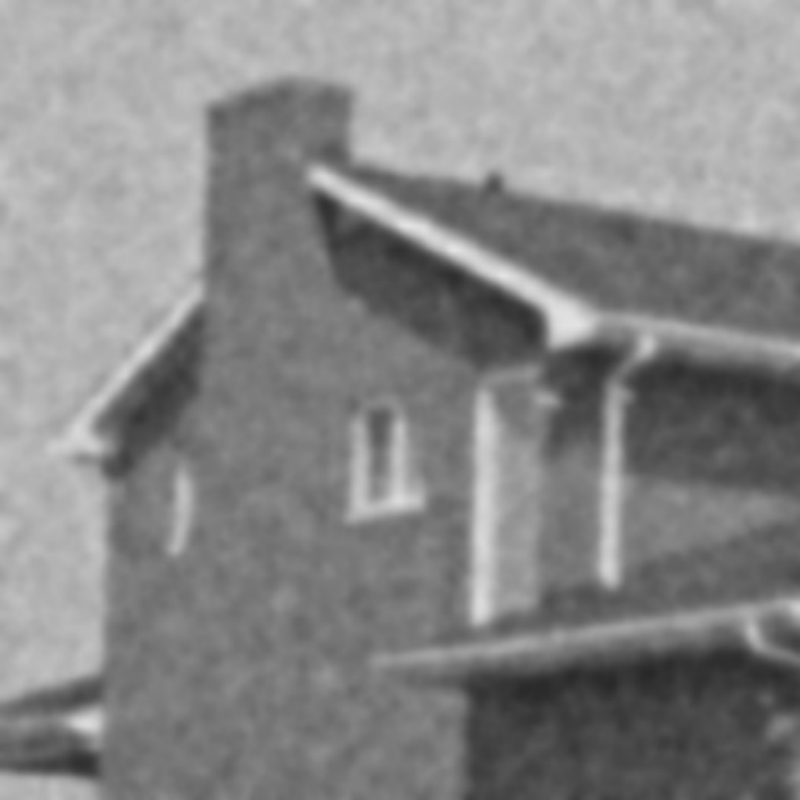}
\caption{Ms$^0$, $16^2$}
\end{subfigure} \
\begin{subfigure}[b]{0.18\textwidth}
\includegraphics[page=1,width=1\linewidth]{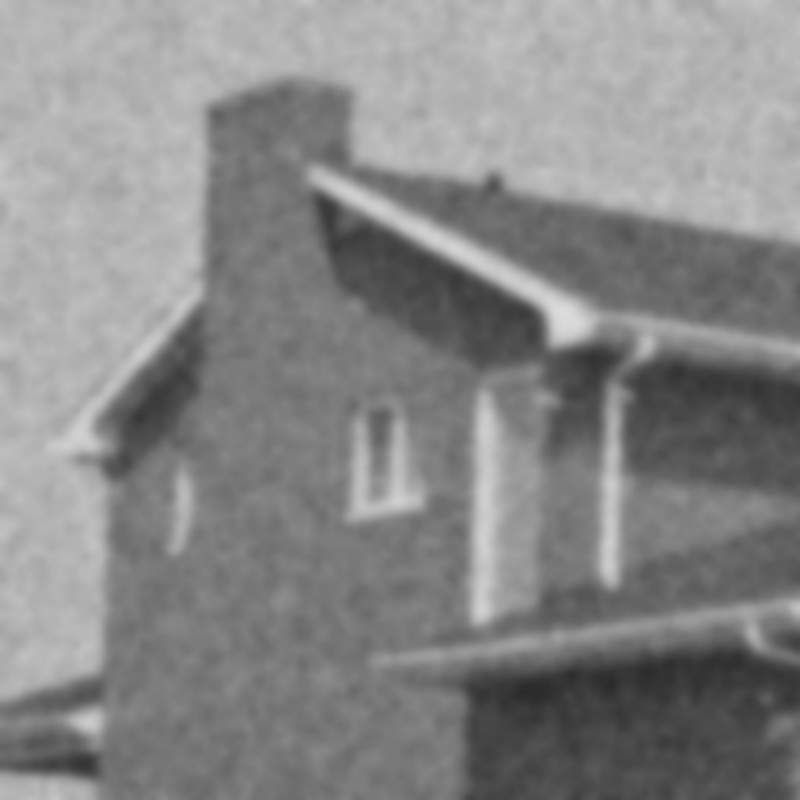}
\caption{Ms$^0$, $32^2$}
\end{subfigure} 
\\
\begin{subfigure}[b]{0.18\textwidth}
\includegraphics[page=1,width=1\linewidth]{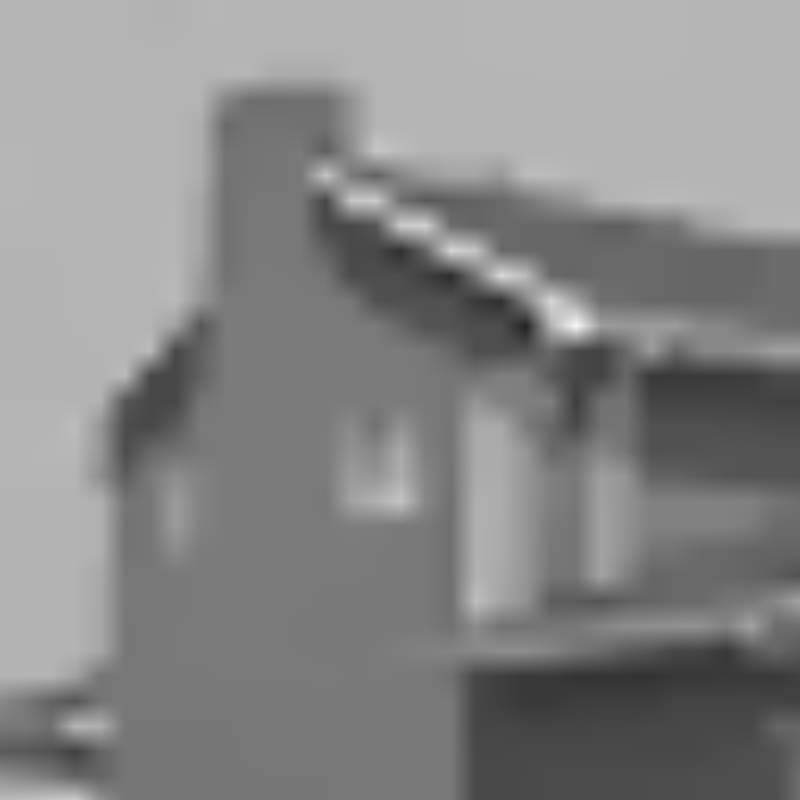}
\caption{AMs1, $16^2$}
\end{subfigure} \
\begin{subfigure}[b]{0.18\textwidth}
\includegraphics[page=1,width=1\linewidth]{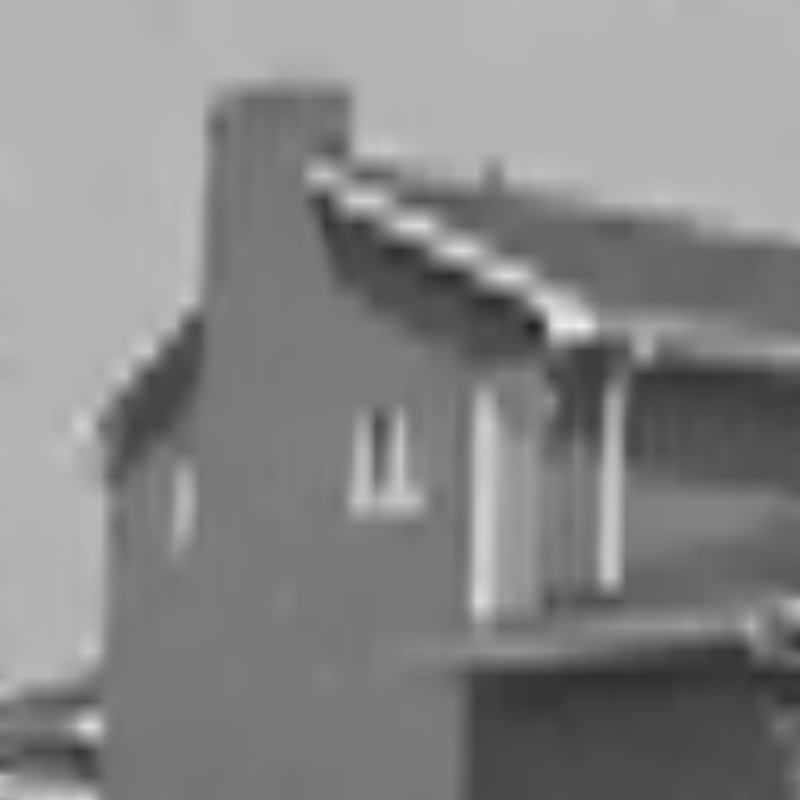}
\caption{AMs2, $16^2$}
\end{subfigure} \
\begin{subfigure}[b]{0.18\textwidth}
\includegraphics[page=1,width=1\linewidth]{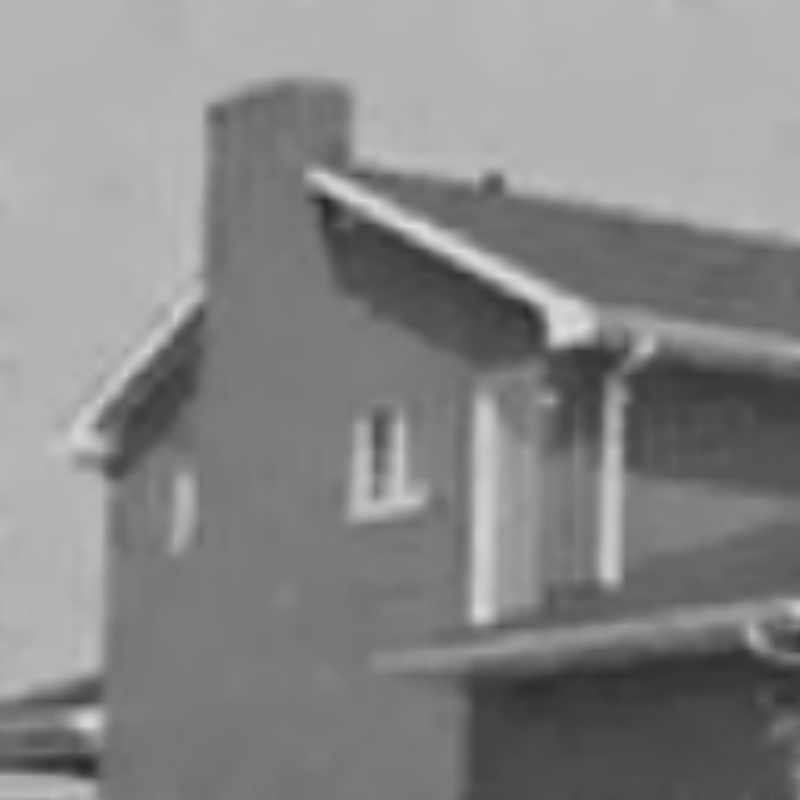}
\caption{AMs4, $16^2$}
\end{subfigure} \
\begin{subfigure}[b]{0.18\textwidth}
\includegraphics[page=1,width=1\linewidth]{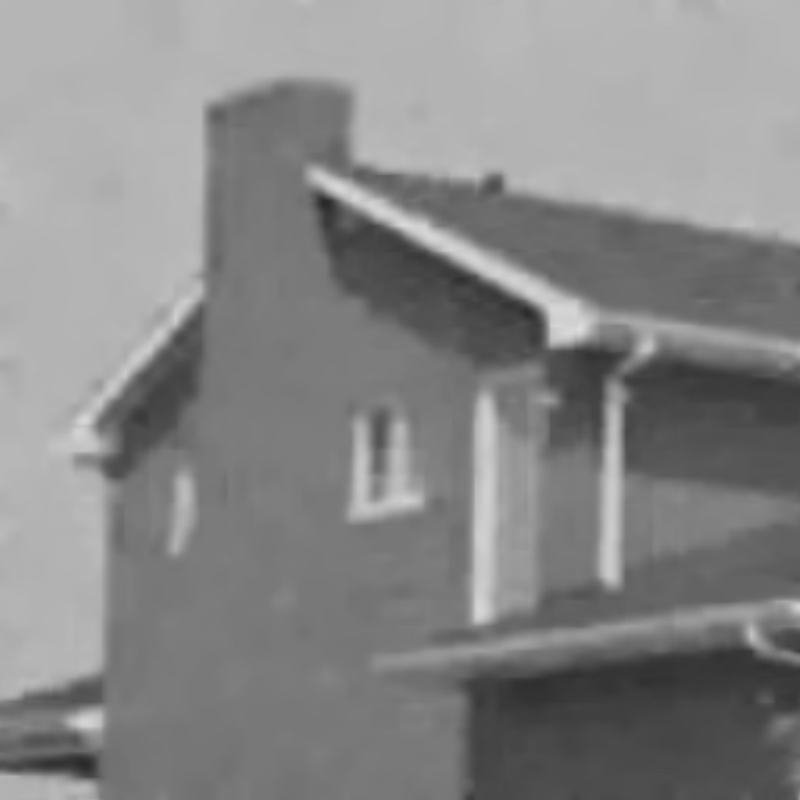}
\caption{AMs8, $16^2$}
\end{subfigure} \
\begin{subfigure}[b]{0.18\textwidth}
\includegraphics[page=1,width=1\linewidth]{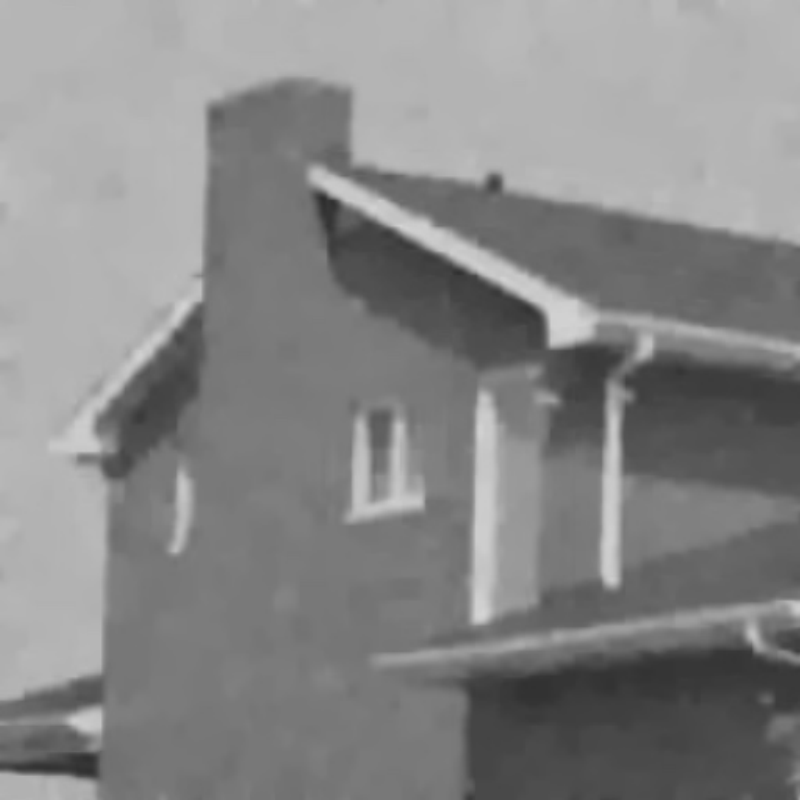} 
\caption{AMs16, $16^2$}
\end{subfigure}
\\
\begin{subfigure}[b]{0.18\textwidth}
\includegraphics[page=1,width=1\linewidth]{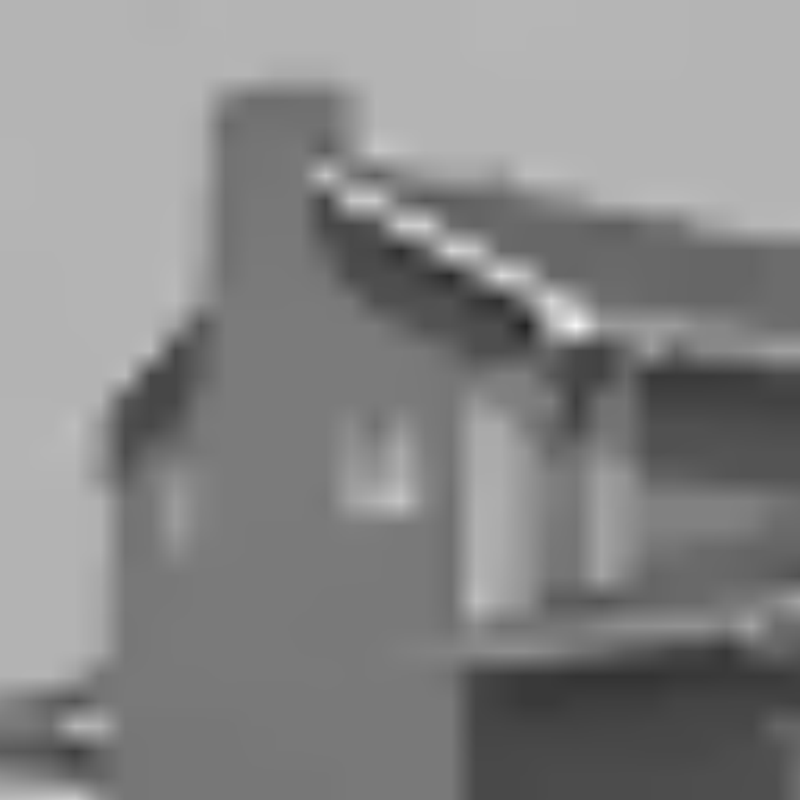}
\caption{GMs1, $16^2$}
\end{subfigure} \ 
\begin{subfigure}[b]{0.18\textwidth}
\includegraphics[page=1,width=1\linewidth]{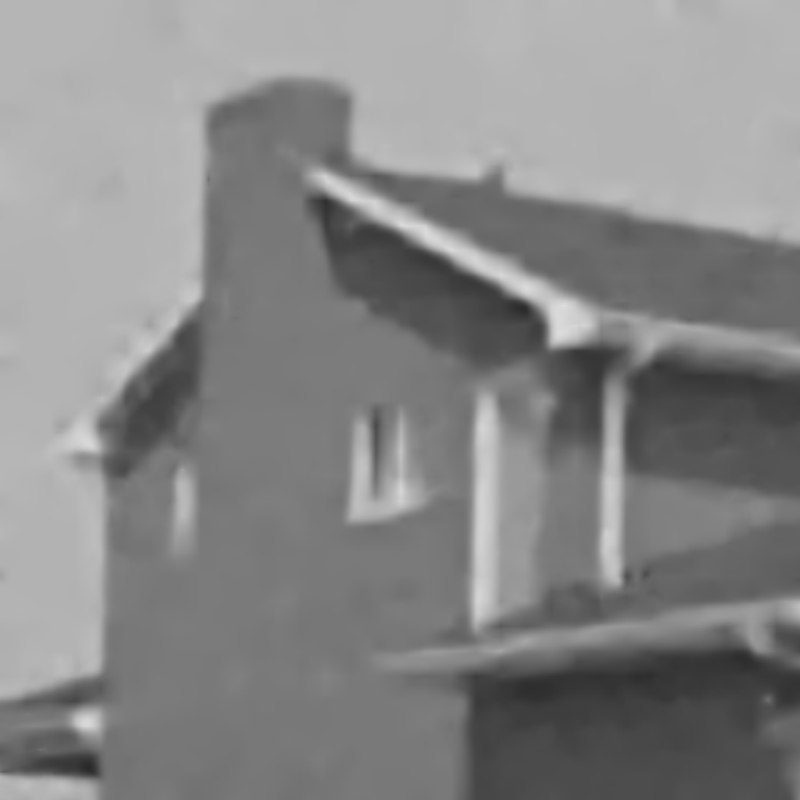}
\caption{GMs2, $16^2$}
\end{subfigure} \ 
\begin{subfigure}[b]{0.18\textwidth}
\includegraphics[page=1,width=1\linewidth]{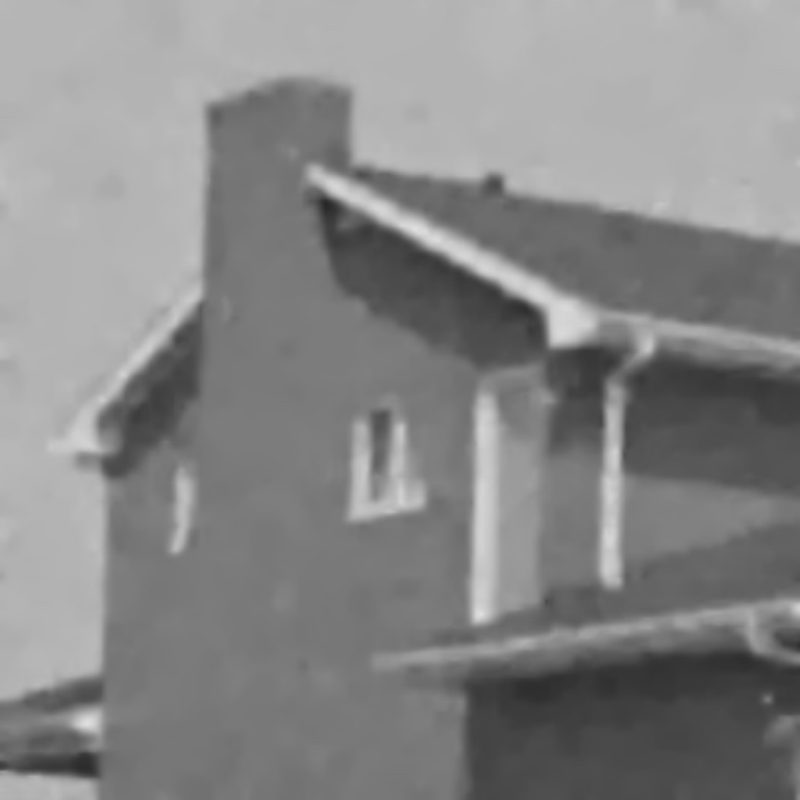}
\caption{GMs4, $16^2$}
\end{subfigure} \ 
\begin{subfigure}[b]{0.18\textwidth}
\includegraphics[page=1,width=1\linewidth]{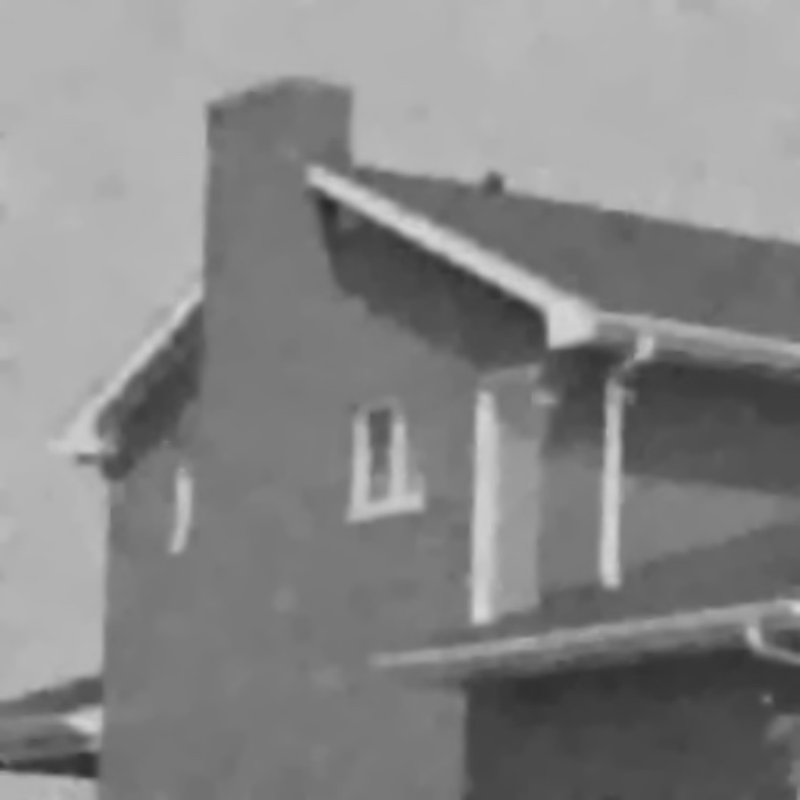}
\caption{GMs8, $16^2$}
\end{subfigure} \ 
\begin{subfigure}[b]{0.18\textwidth}
\includegraphics[page=1,width=1\linewidth]{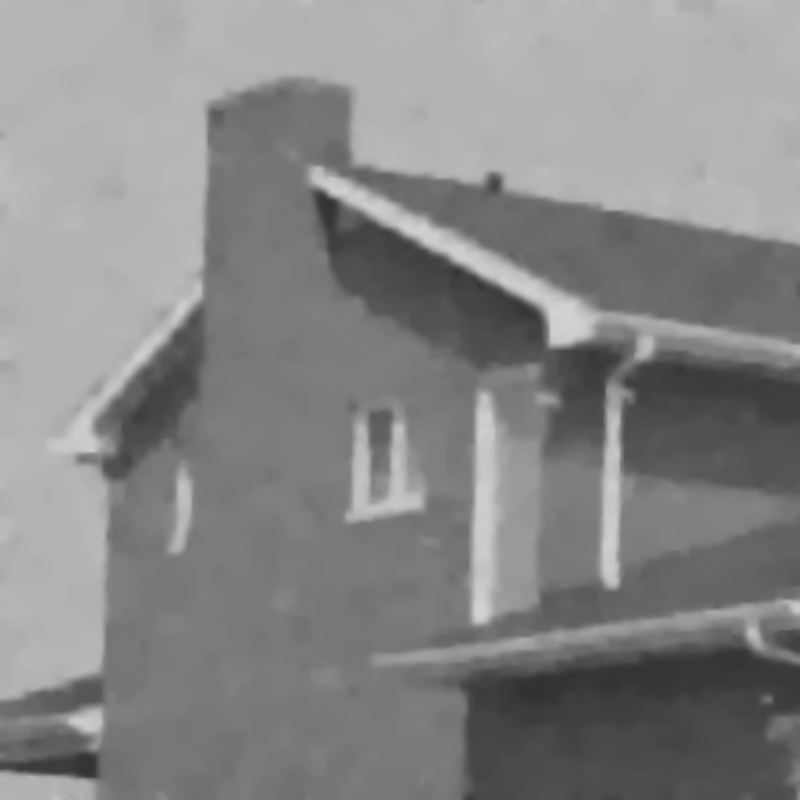}  
\caption{GMs16, $16^2$}
\end{subfigure}
\\
\begin{subfigure}[b]{0.18\textwidth}
\includegraphics[page=1,width=1\linewidth]{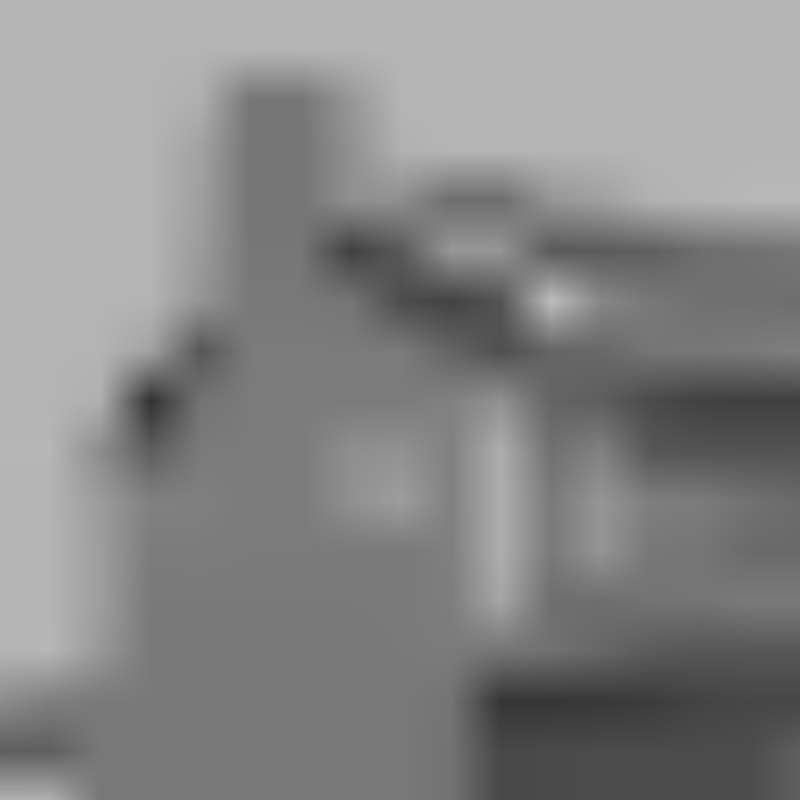}
\caption{AMs1, $32^2$}
\end{subfigure} \ 
\begin{subfigure}[b]{0.18\textwidth}
\includegraphics[page=1,width=1\linewidth]{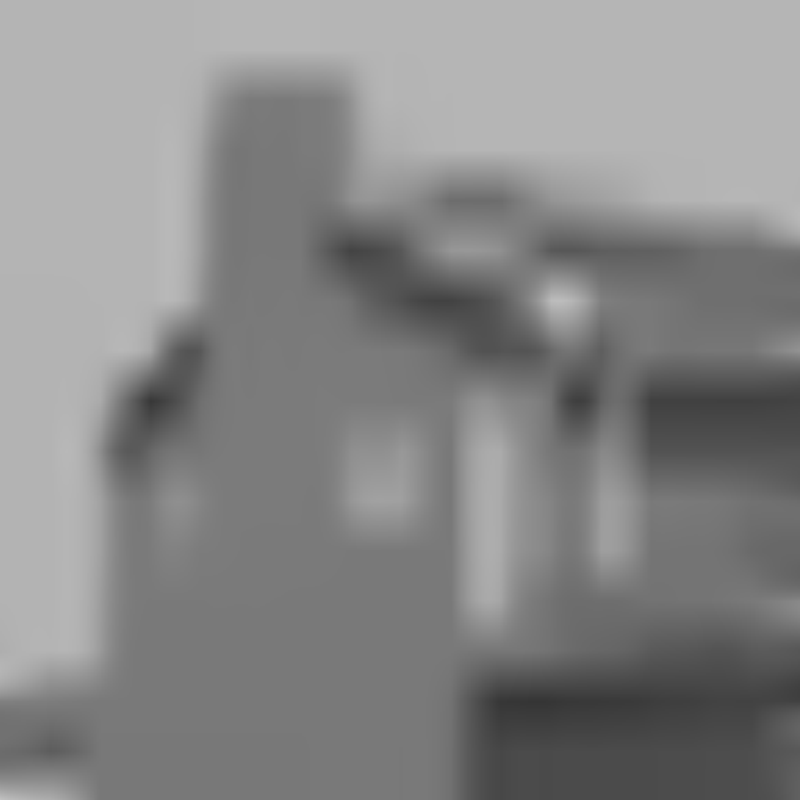}
\caption{AMs2, $32^2$}
\end{subfigure} \ 
\begin{subfigure}[b]{0.18\textwidth}
\includegraphics[page=1,width=1\linewidth]{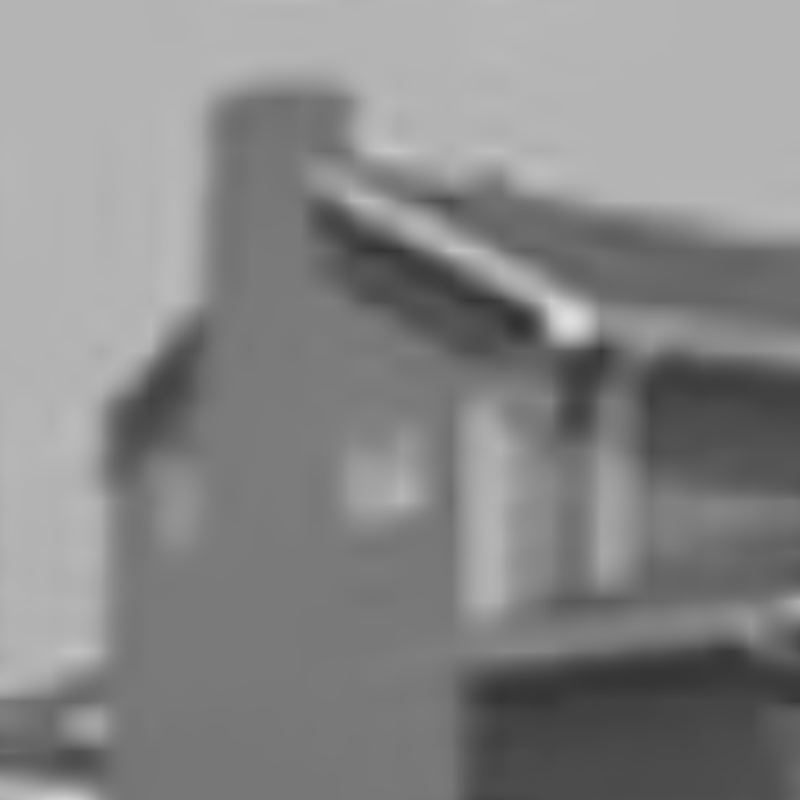}
\caption{AMs4, $32^2$}
\end{subfigure} \ 
\begin{subfigure}[b]{0.18\textwidth}
\includegraphics[page=1,width=1\linewidth]{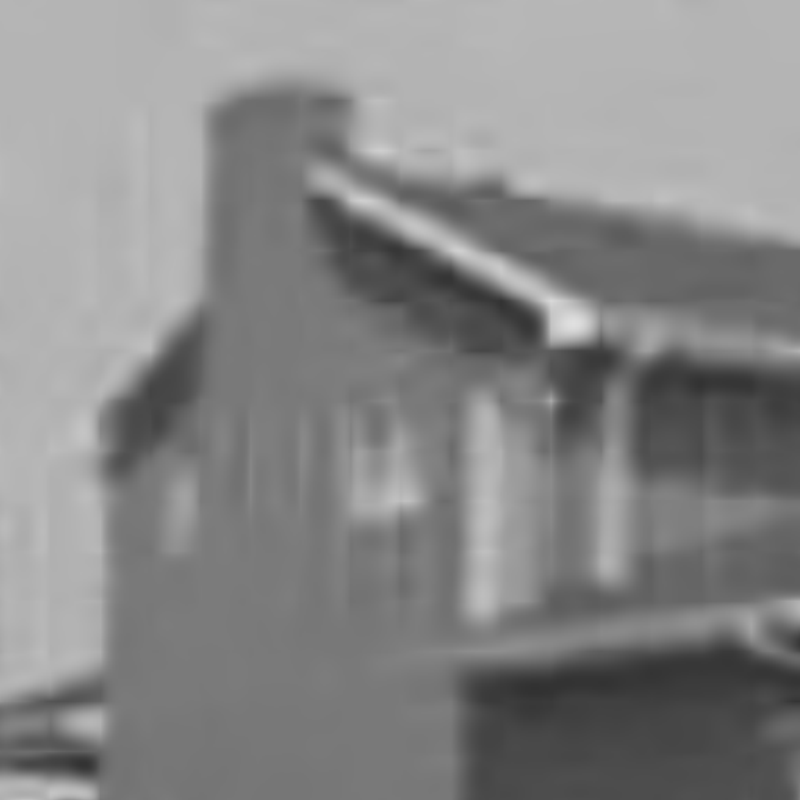}
\caption{AMs8, $32^2$}
\end{subfigure} \ 
\begin{subfigure}[b]{0.18\textwidth}
\includegraphics[page=1,width=1\linewidth]{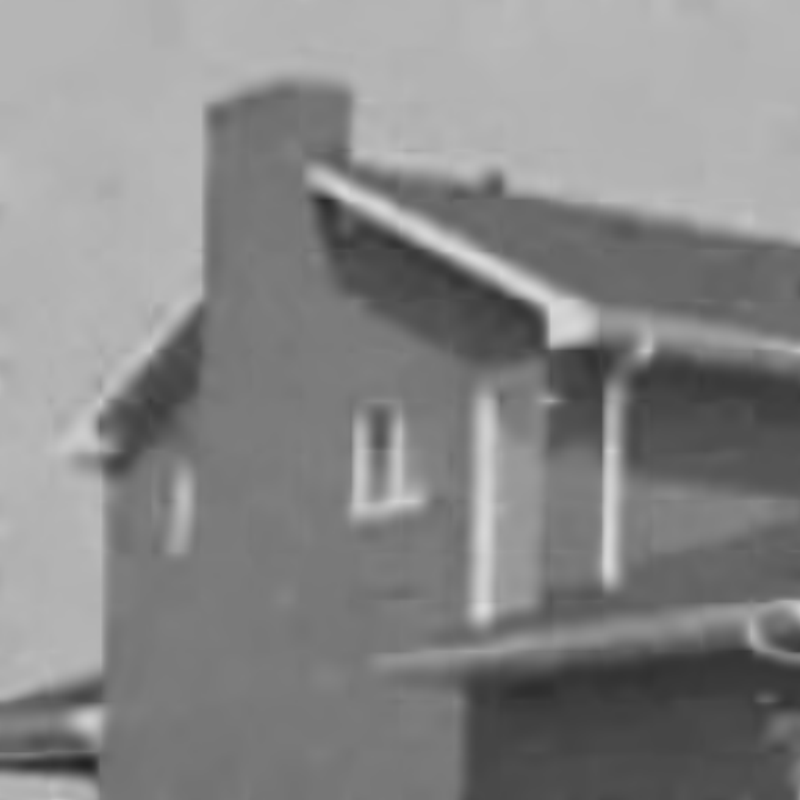} 
\caption{AMs16, $32^2$}
\end{subfigure}
\\
\begin{subfigure}[b]{0.18\textwidth}
\includegraphics[page=1,width=1\linewidth]{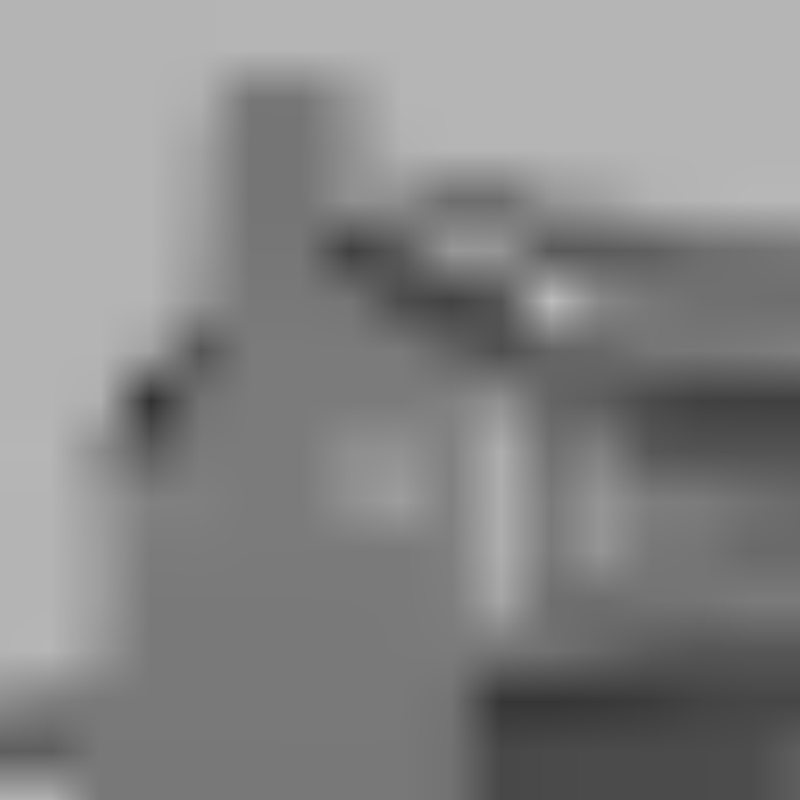}
\caption{GMs1, $32^2$}
\end{subfigure} \ 
\begin{subfigure}[b]{0.18\textwidth}
\includegraphics[page=1,width=1\linewidth]{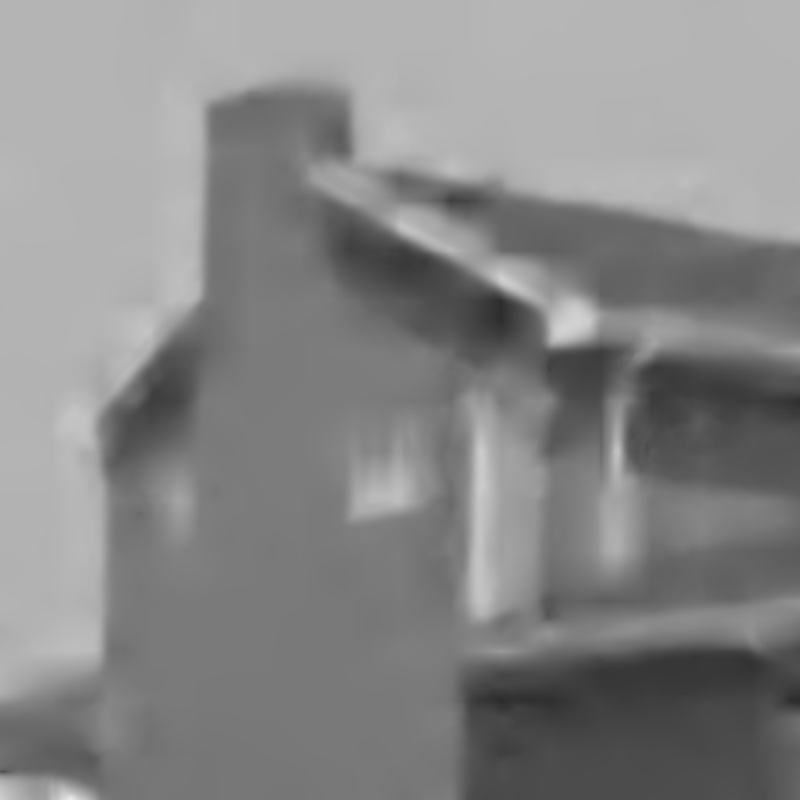}
\caption{GMs2, $32^2$}
\end{subfigure} \ 
\begin{subfigure}[b]{0.18\textwidth}
\includegraphics[page=1,width=1\linewidth]{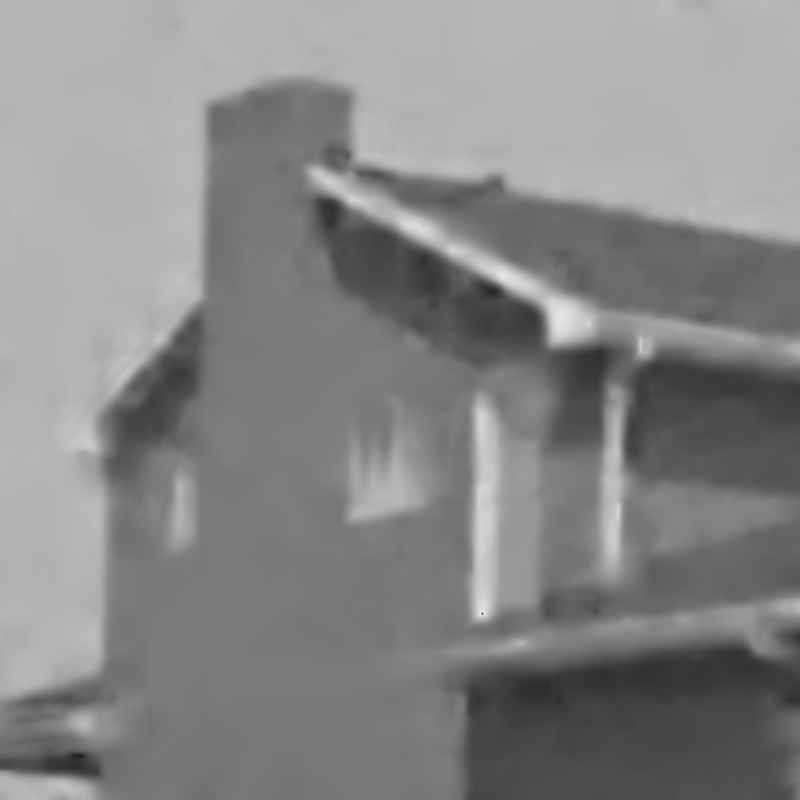}
\caption{GMs4, $32^2$}
\end{subfigure} \ 
\begin{subfigure}[b]{0.18\textwidth}
\includegraphics[page=1,width=1\linewidth]{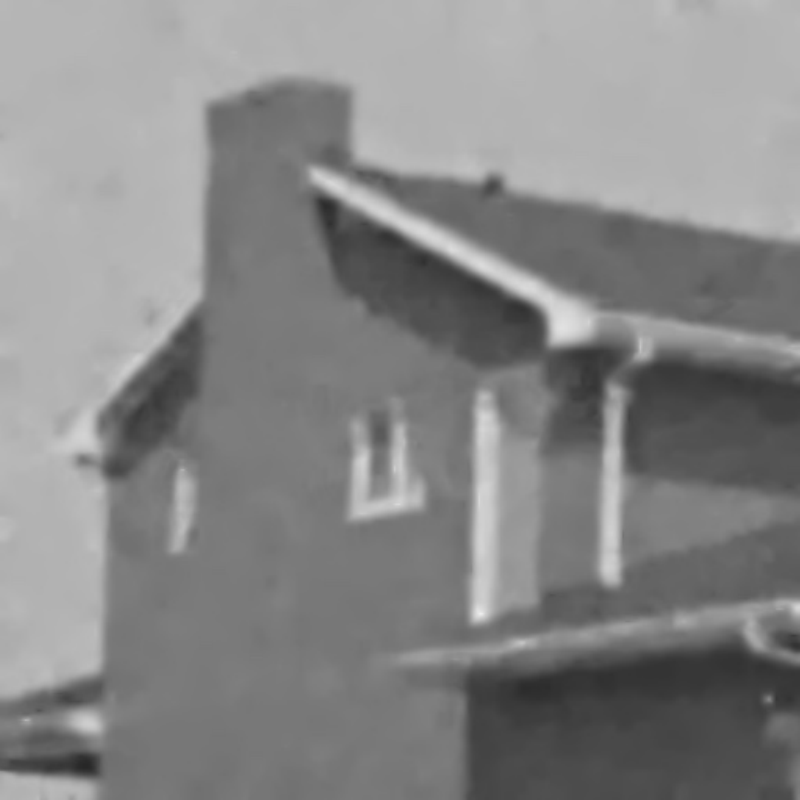}
\caption{GMs8, $32^2$}
\end{subfigure} \ 
\begin{subfigure}[b]{0.18\textwidth}
\includegraphics[page=1,width=1\linewidth]{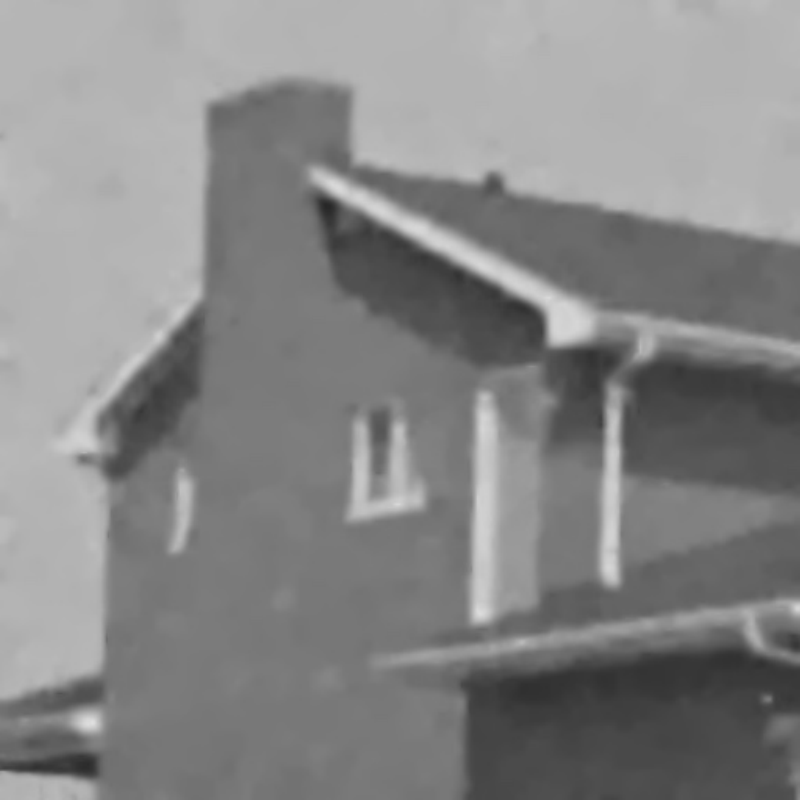}  
\caption{GMs16, $32^2$}
\end{subfigure}
\caption{\rev{Test 3. Image 1 (House). Image denoising using multiscale method. 40 \% of noise}}
\label{fig:t3-noi40-1}
\end{figure}

\begin{figure}[h!]
\centering
\begin{subfigure}[b]{0.18\textwidth}
\includegraphics[page=2,width=1\linewidth]{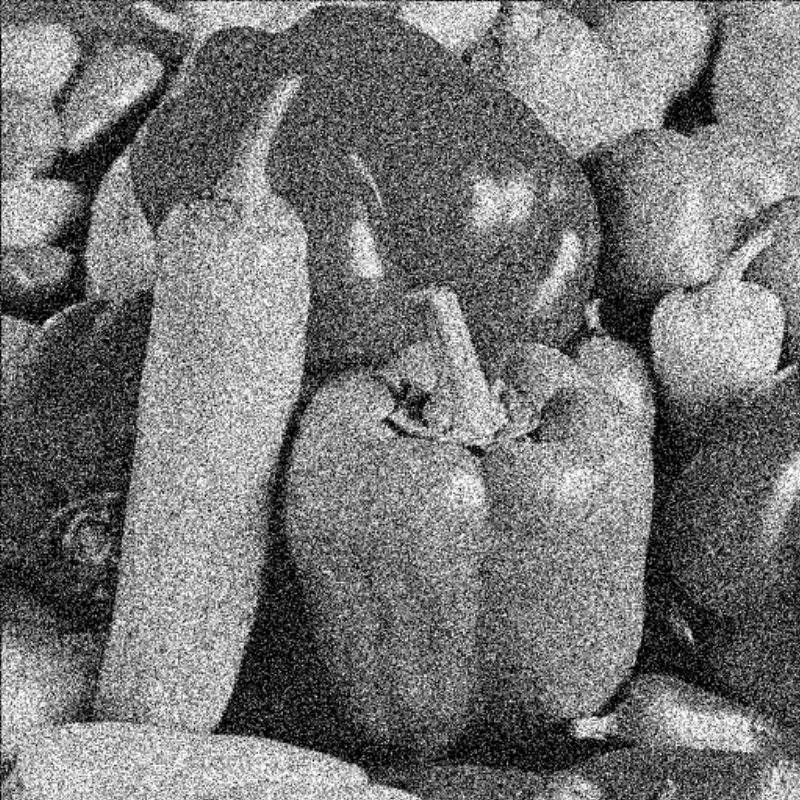} 
\caption{Noised, $I^0(x)$}
\end{subfigure} \ 
\begin{subfigure}[b]{0.18\textwidth}
\includegraphics[page=2,width=1\linewidth]{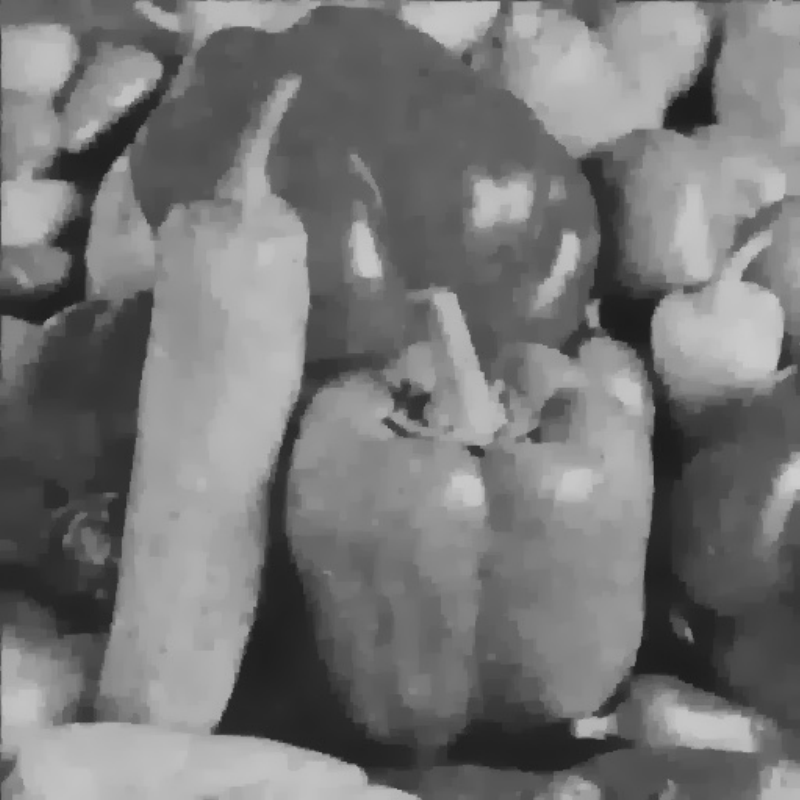} 
\caption{Ex-ROF}
\end{subfigure} \ 
\begin{subfigure}[b]{0.18\textwidth}
\includegraphics[page=2,width=1\linewidth]{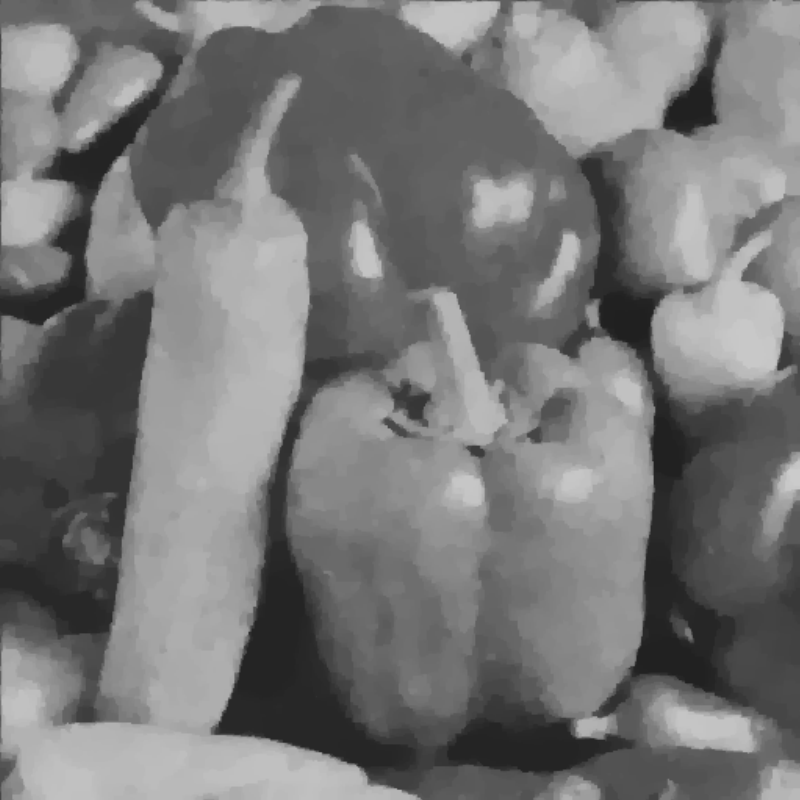} 
\caption{Im-ROF}
\end{subfigure} \ 
\begin{subfigure}[b]{0.18\textwidth}
\includegraphics[page=2,width=1\linewidth]{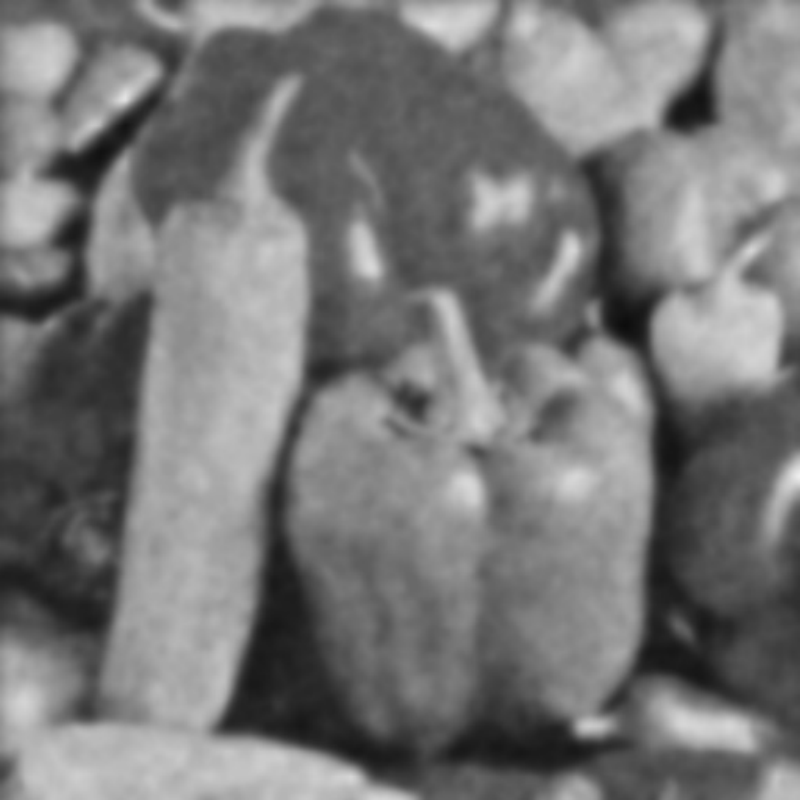}
\caption{Ms$^0$, $16^2$}
\end{subfigure} \
\begin{subfigure}[b]{0.18\textwidth}
\includegraphics[page=2,width=1\linewidth]{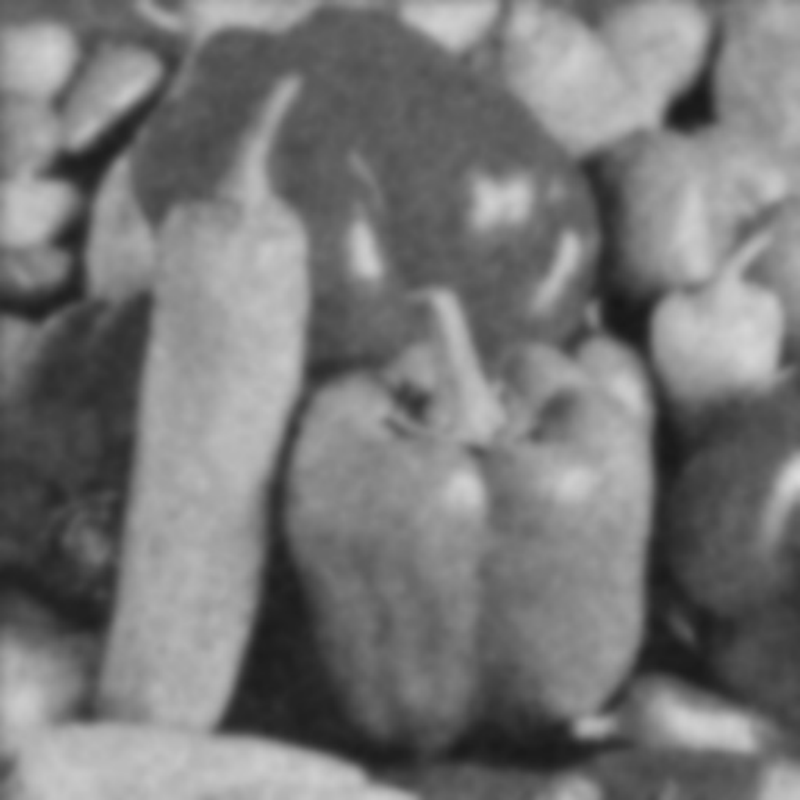}
\caption{Ms$^0$, $32^2$}
\end{subfigure} 
\\
\begin{subfigure}[b]{0.18\textwidth}
\includegraphics[page=2,width=1\linewidth]{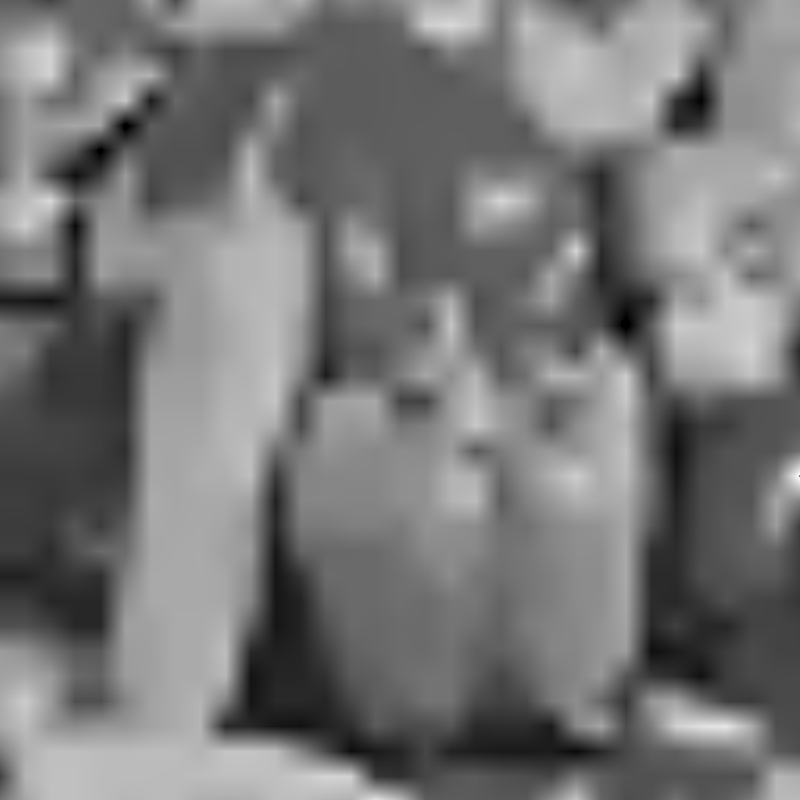}
\caption{AMs1, $16^2$}
\end{subfigure} \
\begin{subfigure}[b]{0.18\textwidth}
\includegraphics[page=2,width=1\linewidth]{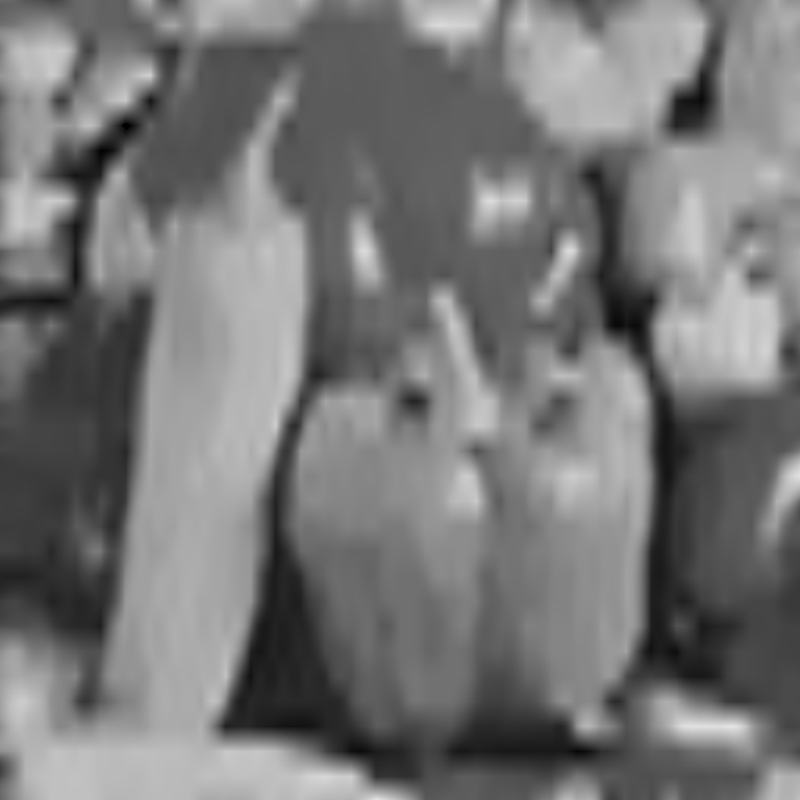}
\caption{AMs2, $16^2$}
\end{subfigure} \
\begin{subfigure}[b]{0.18\textwidth}
\includegraphics[page=2,width=1\linewidth]{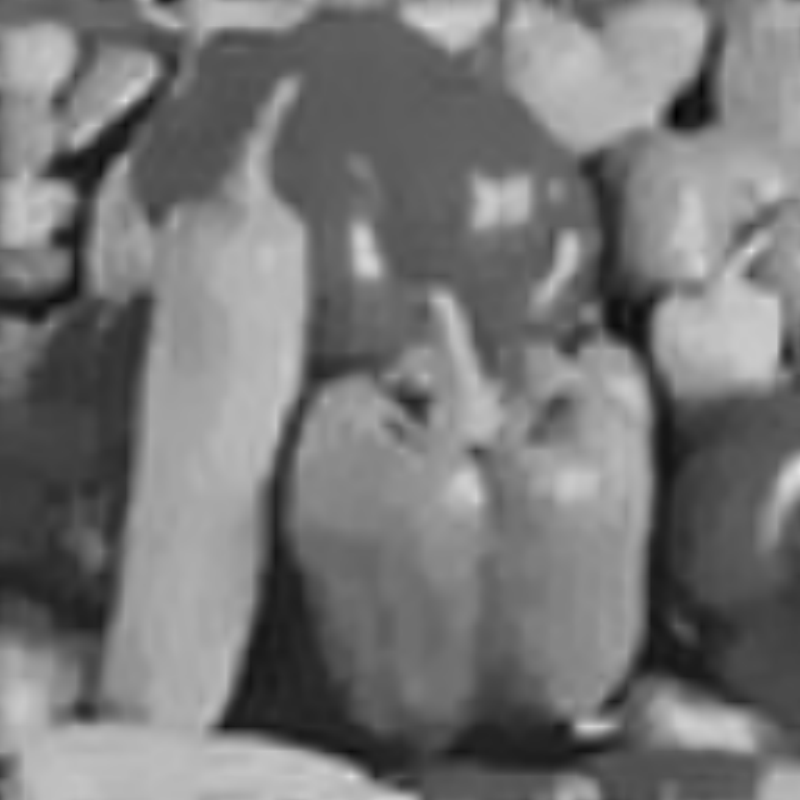}
\caption{AMs4, $16^2$}
\end{subfigure} \
\begin{subfigure}[b]{0.18\textwidth}
\includegraphics[page=2,width=1\linewidth]{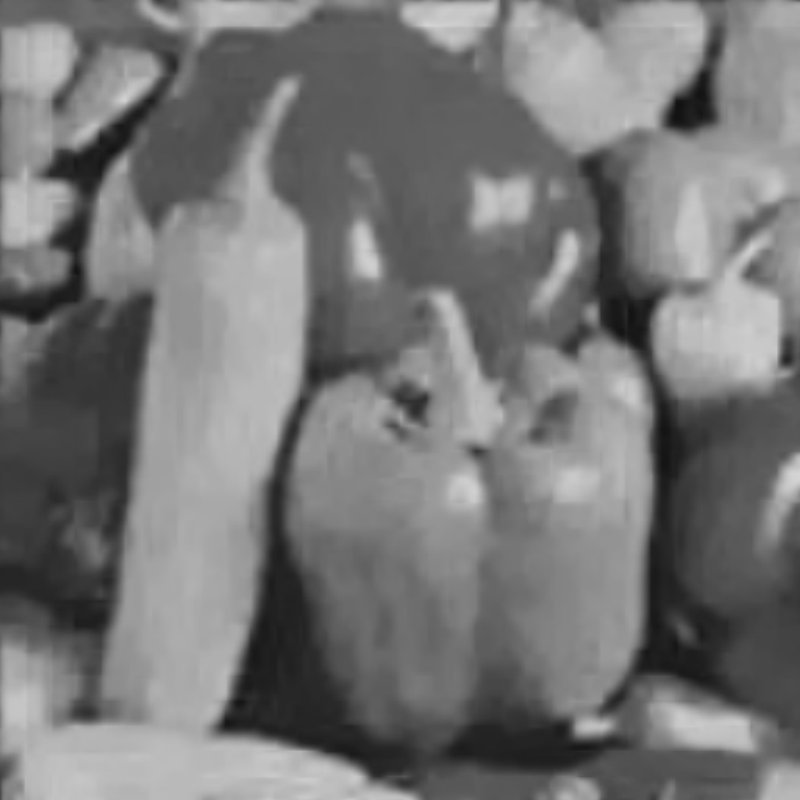}
\caption{AMs8, $16^2$}
\end{subfigure} \
\begin{subfigure}[b]{0.18\textwidth}
\includegraphics[page=2,width=1\linewidth]{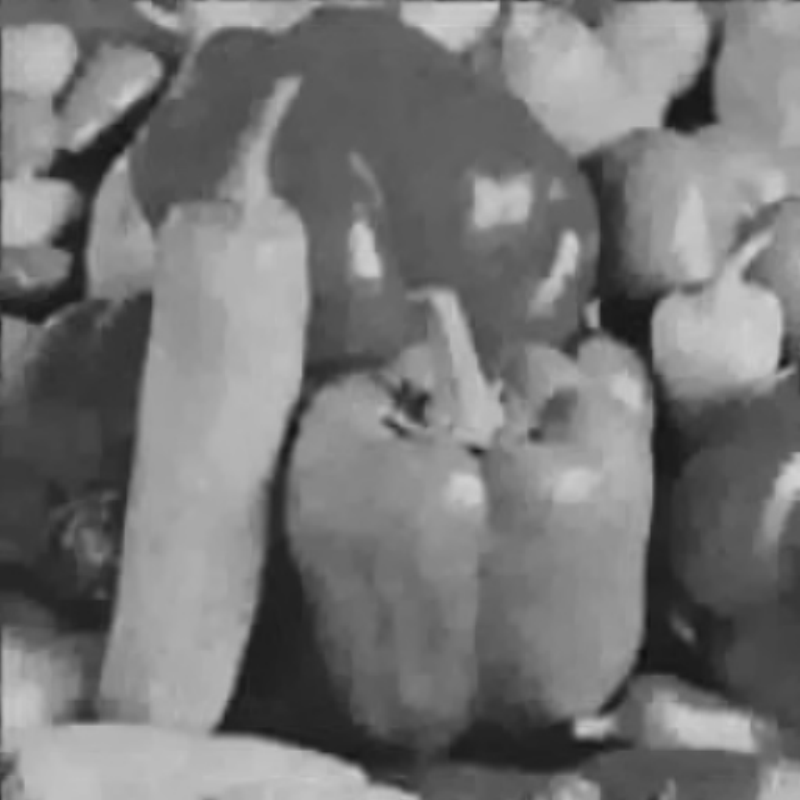} 
\caption{AMs16, $16^2$}
\end{subfigure}
\\
\begin{subfigure}[b]{0.18\textwidth}
\includegraphics[page=2,width=1\linewidth]{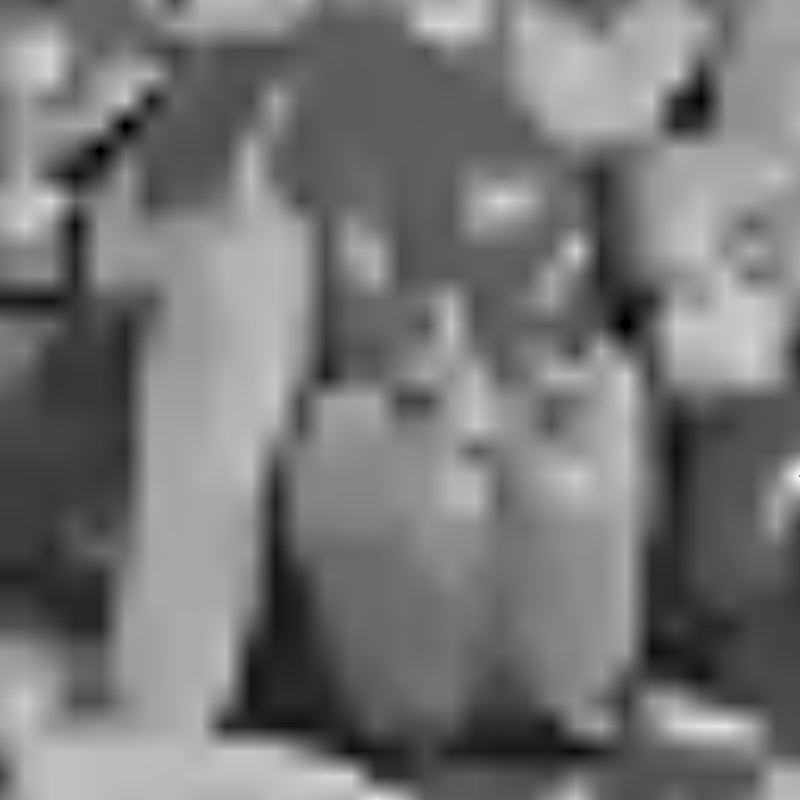}
\caption{GMs1, $16^2$}
\end{subfigure} \ 
\begin{subfigure}[b]{0.18\textwidth}
\includegraphics[page=2,width=1\linewidth]{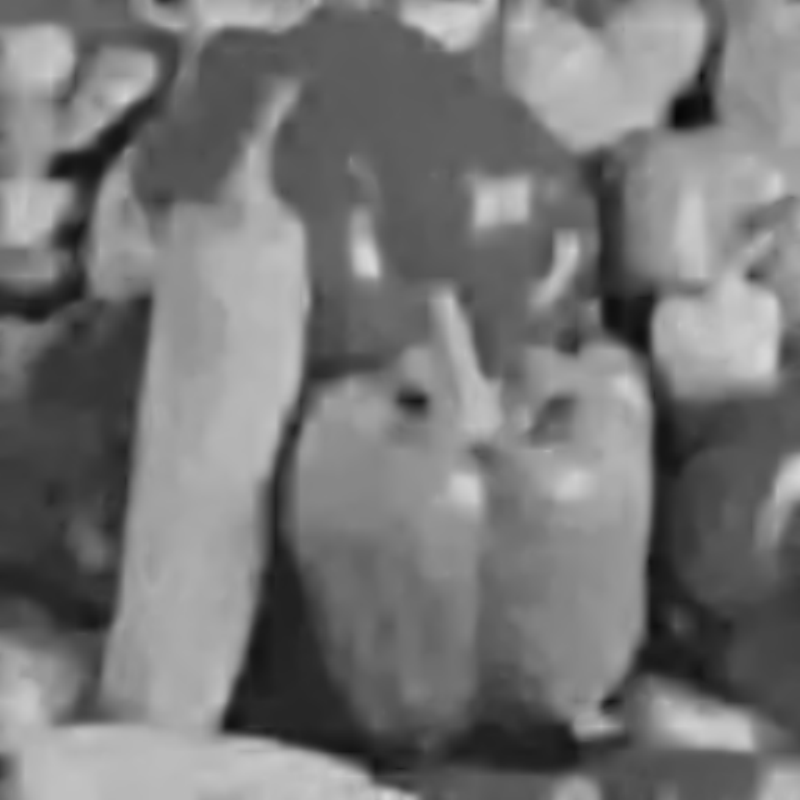}
\caption{GMs2, $16^2$}
\end{subfigure} \ 
\begin{subfigure}[b]{0.18\textwidth}
\includegraphics[page=2,width=1\linewidth]{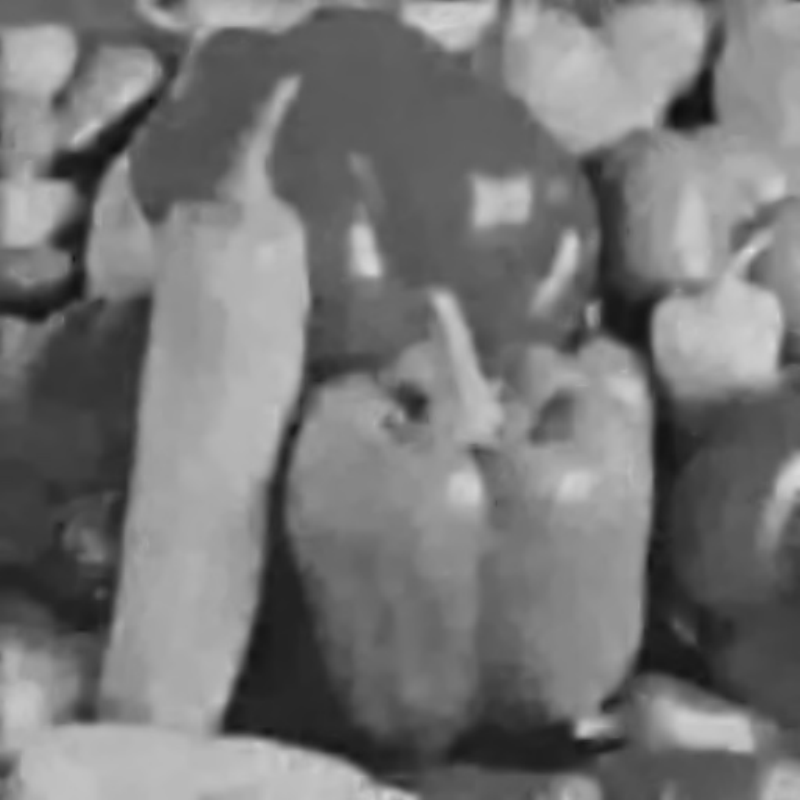}
\caption{GMs4, $16^2$}
\end{subfigure} \ 
\begin{subfigure}[b]{0.18\textwidth}
\includegraphics[page=2,width=1\linewidth]{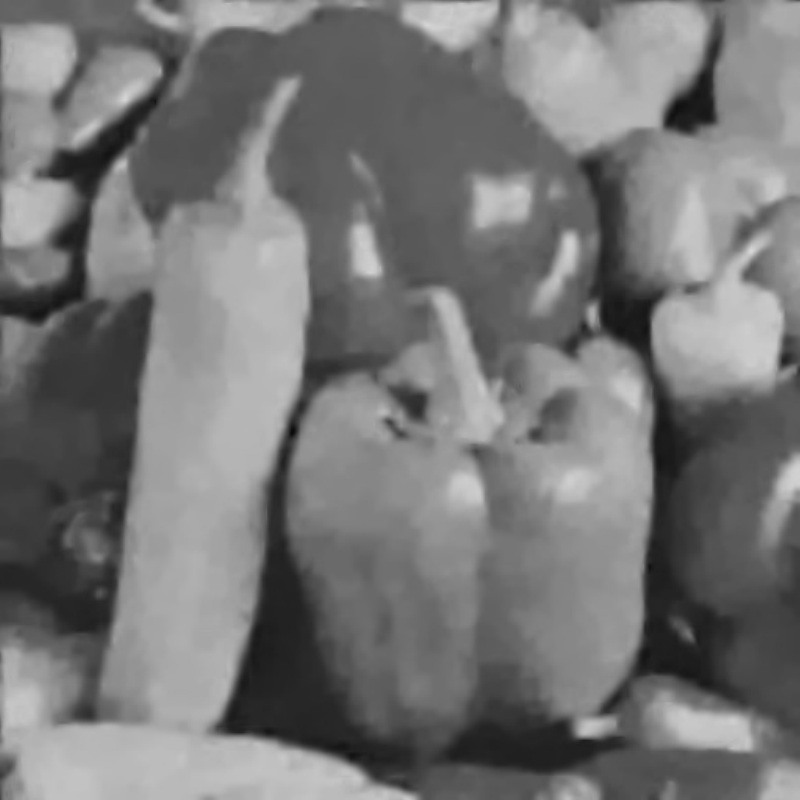}
\caption{GMs8, $16^2$}
\end{subfigure} \ 
\begin{subfigure}[b]{0.18\textwidth}
\includegraphics[page=2,width=1\linewidth]{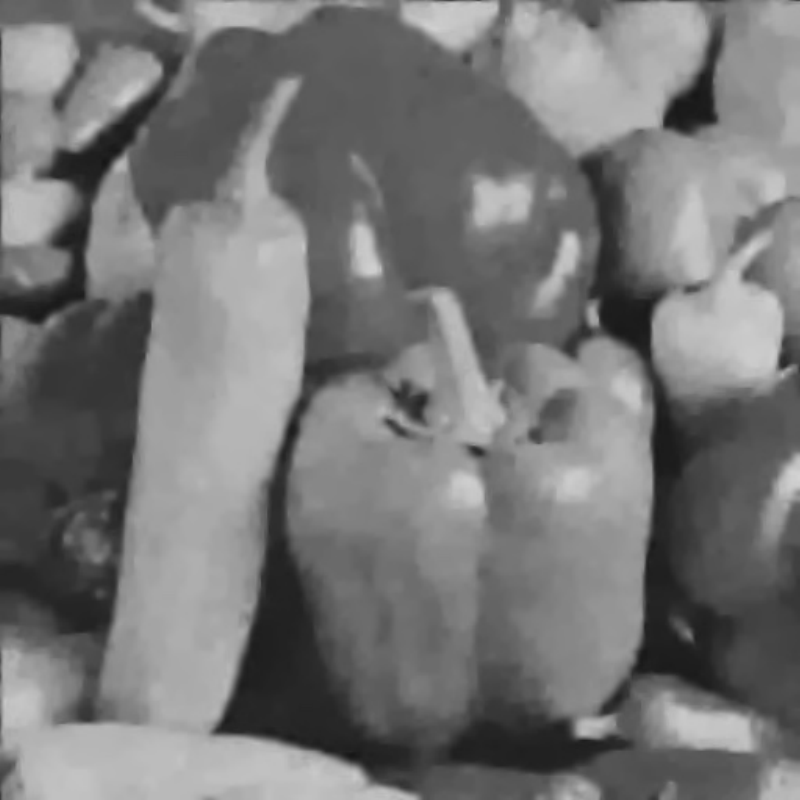}  
\caption{GMs16, $16^2$}
\end{subfigure}
\\
\begin{subfigure}[b]{0.18\textwidth}
\includegraphics[page=2,width=1\linewidth]{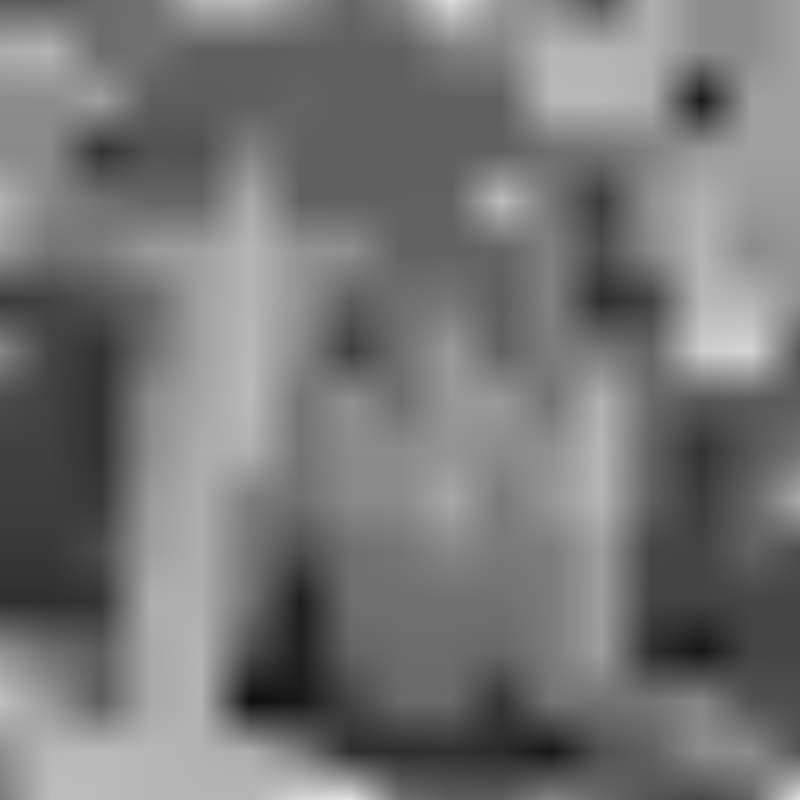}
\caption{AMs1, $32^2$}
\end{subfigure} \ 
\begin{subfigure}[b]{0.18\textwidth}
\includegraphics[page=2,width=1\linewidth]{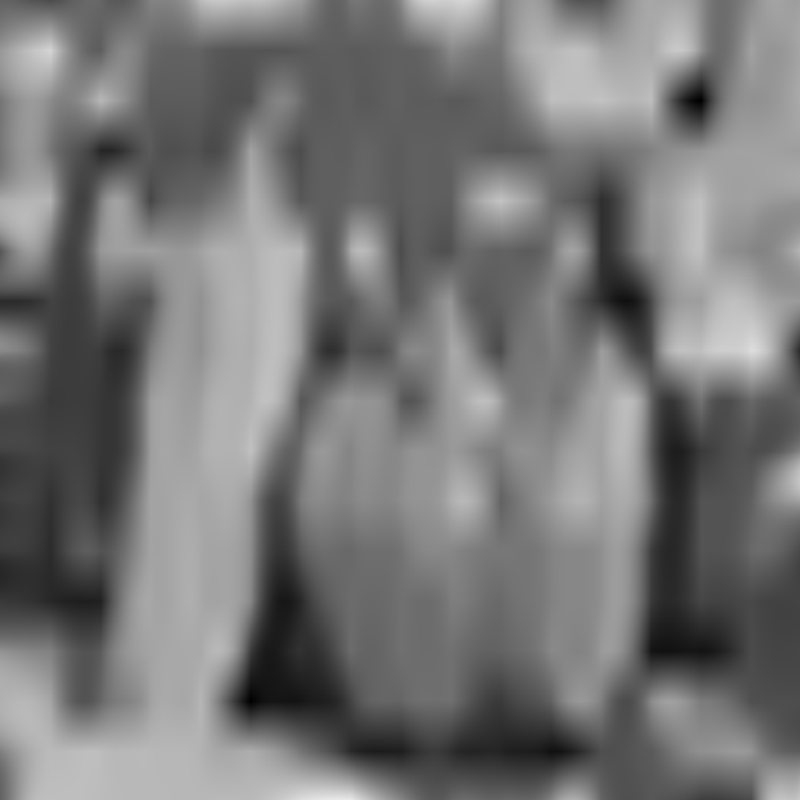}
\caption{AMs2, $32^2$}
\end{subfigure} \ 
\begin{subfigure}[b]{0.18\textwidth}
\includegraphics[page=2,width=1\linewidth]{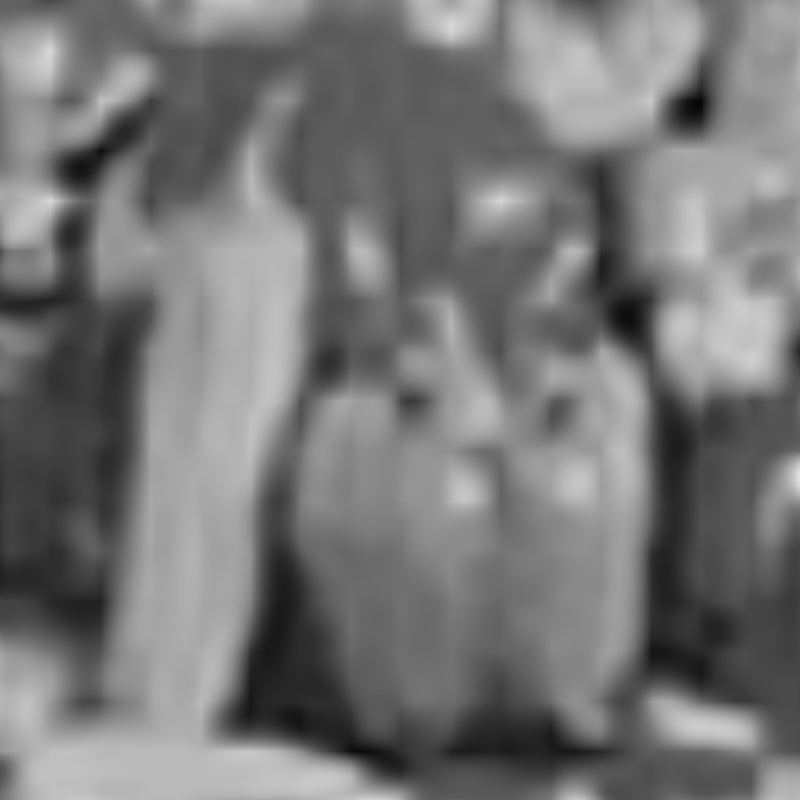}
\caption{AMs4, $32^2$}
\end{subfigure} \ 
\begin{subfigure}[b]{0.18\textwidth}
\includegraphics[page=2,width=1\linewidth]{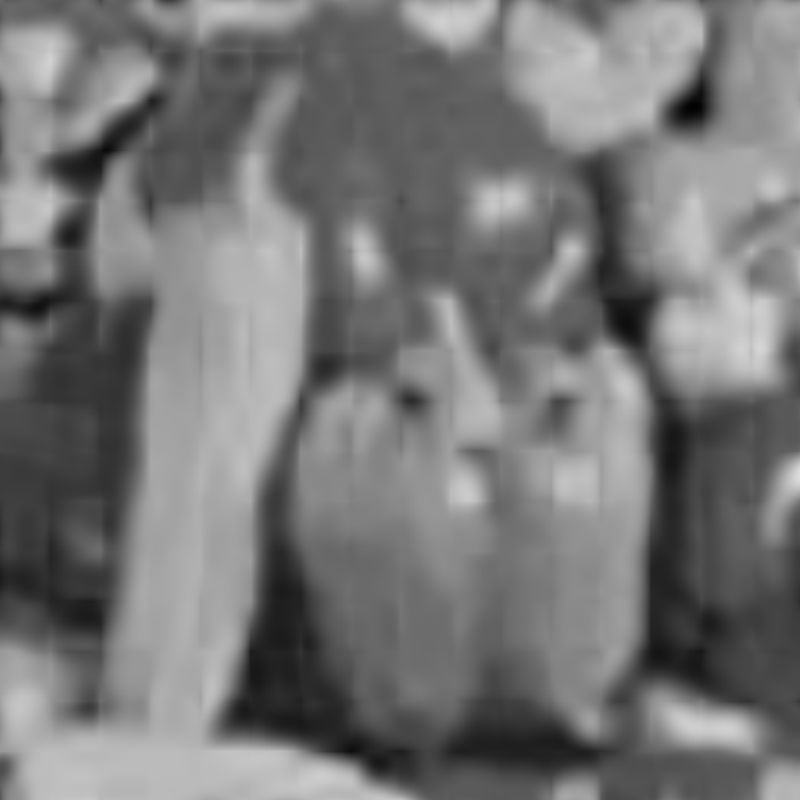}
\caption{AMs8, $32^2$}
\end{subfigure} \ 
\begin{subfigure}[b]{0.18\textwidth}
\includegraphics[page=2,width=1\linewidth]{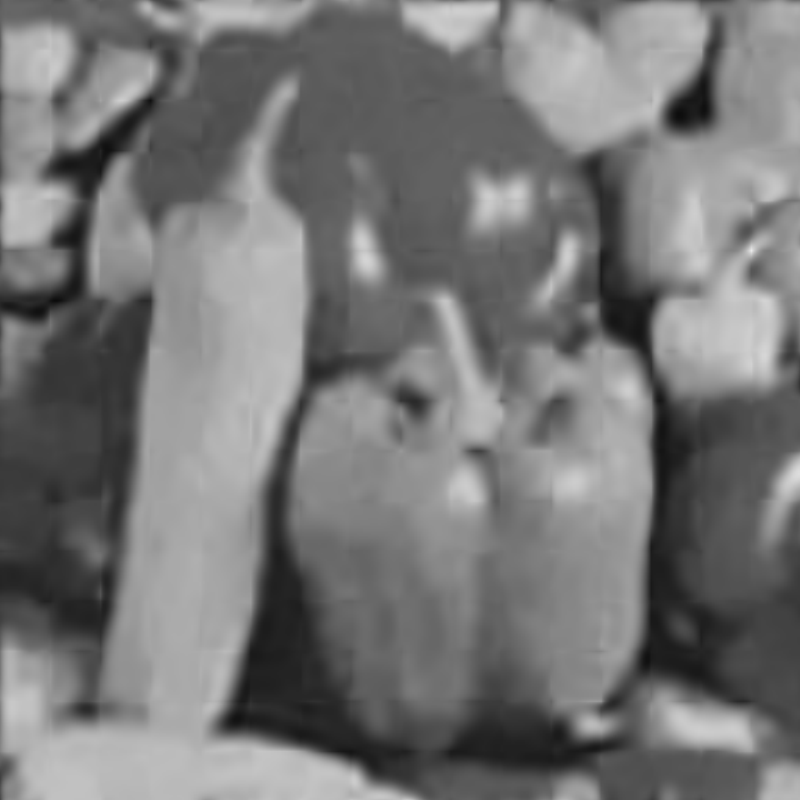} 
\caption{AMs16, $32^2$}
\end{subfigure}
\\
\begin{subfigure}[b]{0.18\textwidth}
\includegraphics[page=2,width=1\linewidth]{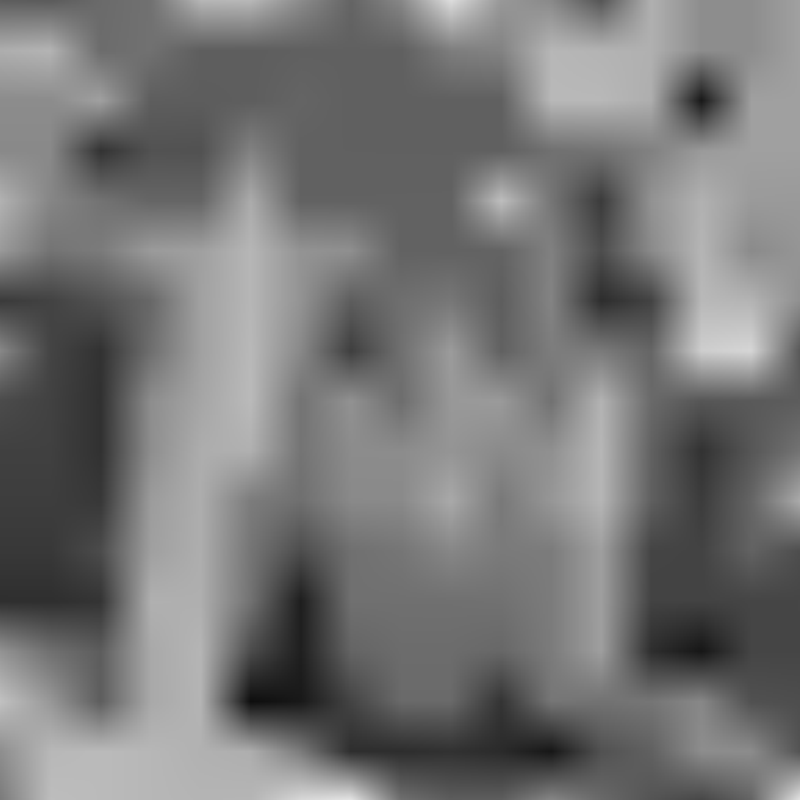}
\caption{GMs1, $32^2$}
\end{subfigure} \ 
\begin{subfigure}[b]{0.18\textwidth}
\includegraphics[page=2,width=1\linewidth]{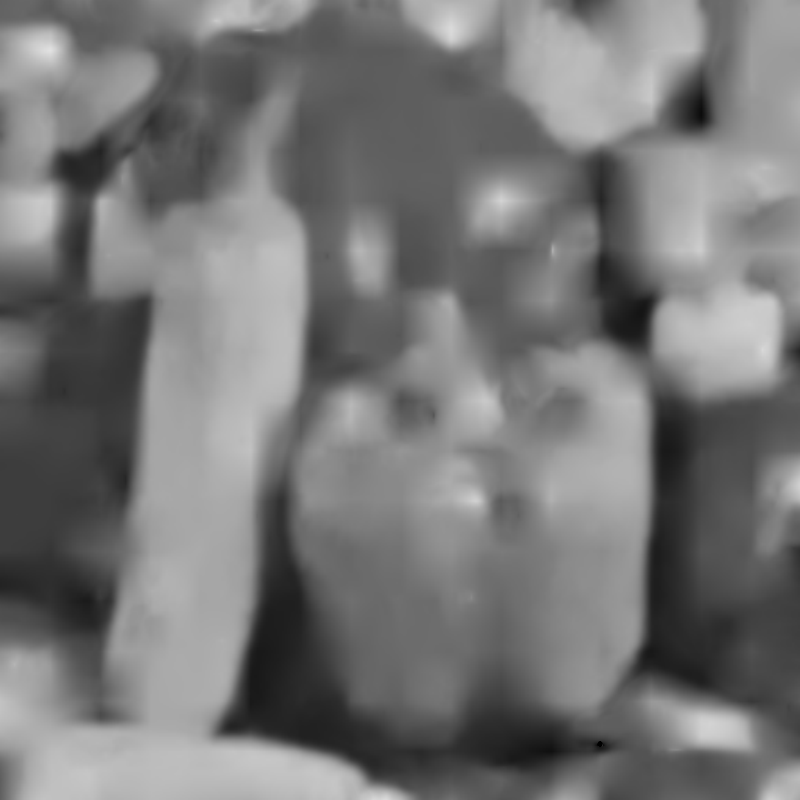}
\caption{GMs2, $32^2$}
\end{subfigure} \ 
\begin{subfigure}[b]{0.18\textwidth}
\includegraphics[page=2,width=1\linewidth]{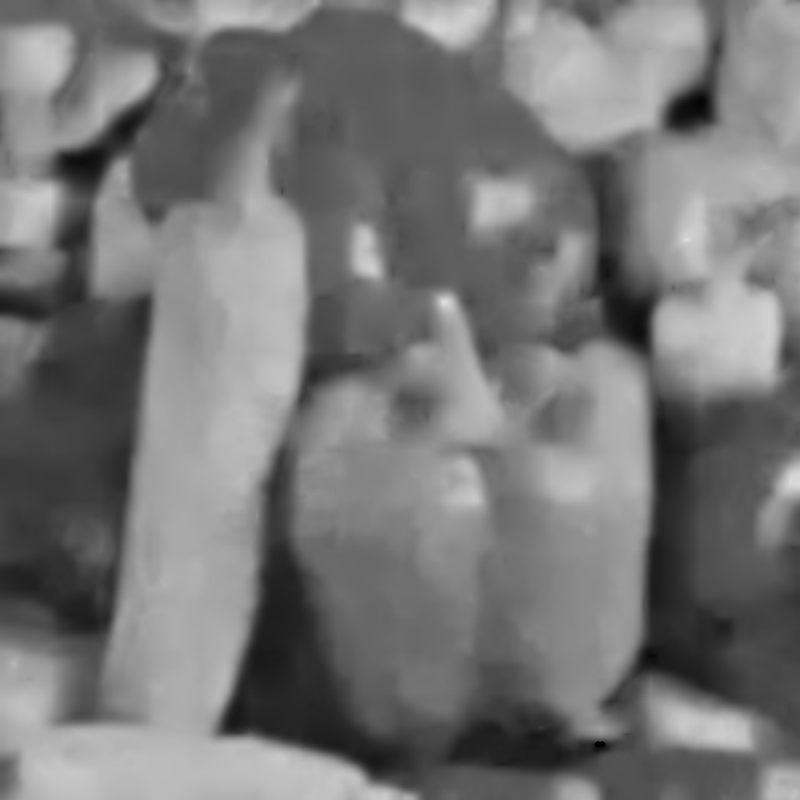}
\caption{GMs4, $32^2$}
\end{subfigure} \ 
\begin{subfigure}[b]{0.18\textwidth}
\includegraphics[page=2,width=1\linewidth]{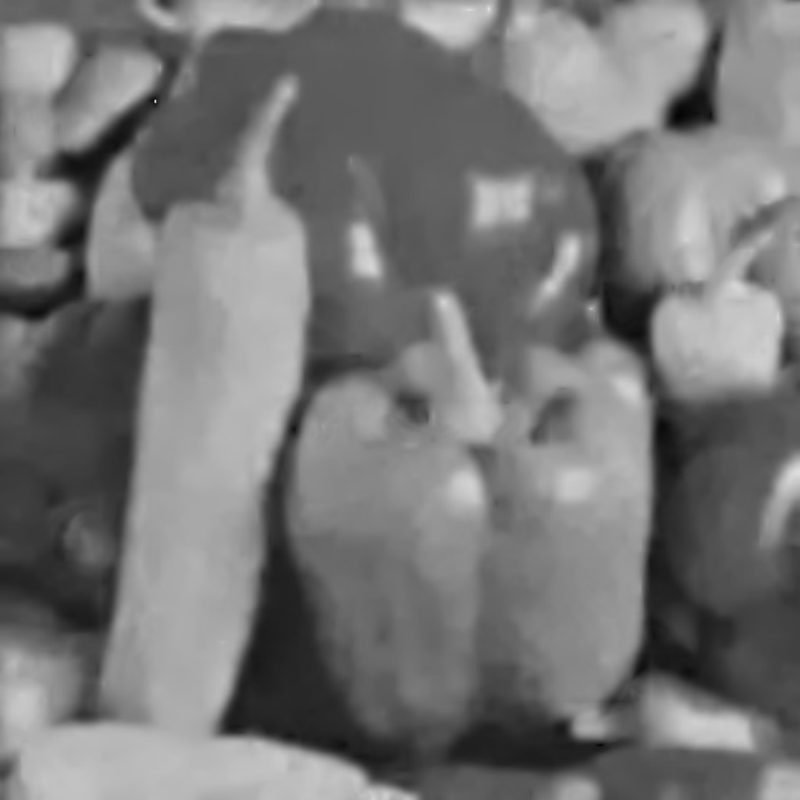}
\caption{GMs8, $32^2$}
\end{subfigure} \ 
\begin{subfigure}[b]{0.18\textwidth}
\includegraphics[page=2,width=1\linewidth]{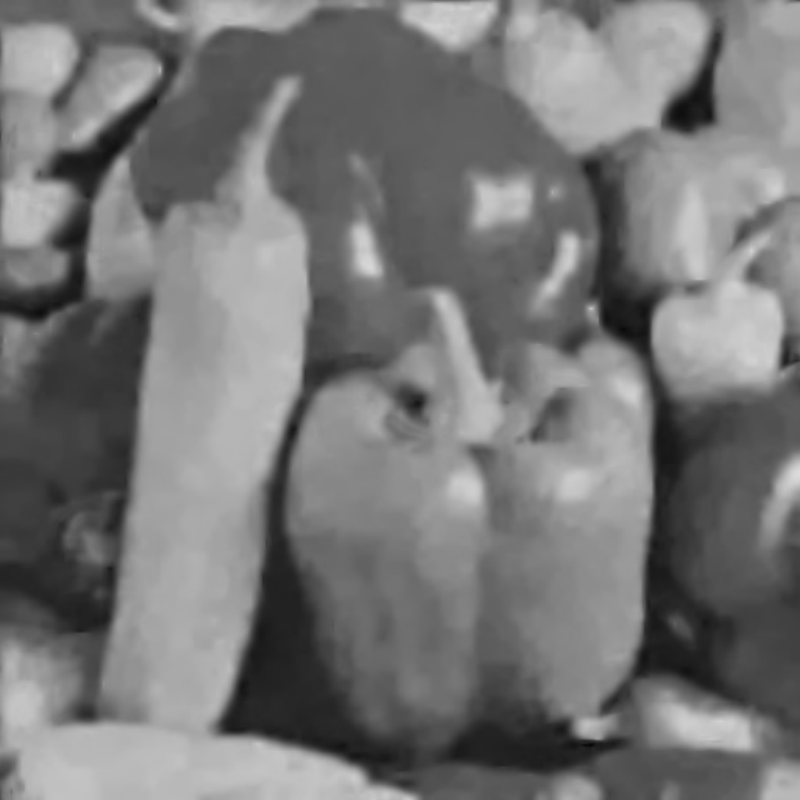}  
\caption{GMs16, $32^2$}
\end{subfigure}
\caption{\rev{Test 3. Image 2 (Peppers). Image denoising using multiscale method. 40 \% of noise}}
\label{fig:t3-noi40-2}
\end{figure}

\begin{figure}[h!]
\centering
\begin{subfigure}[b]{0.45\textwidth}  
\includegraphics[width=1\linewidth]{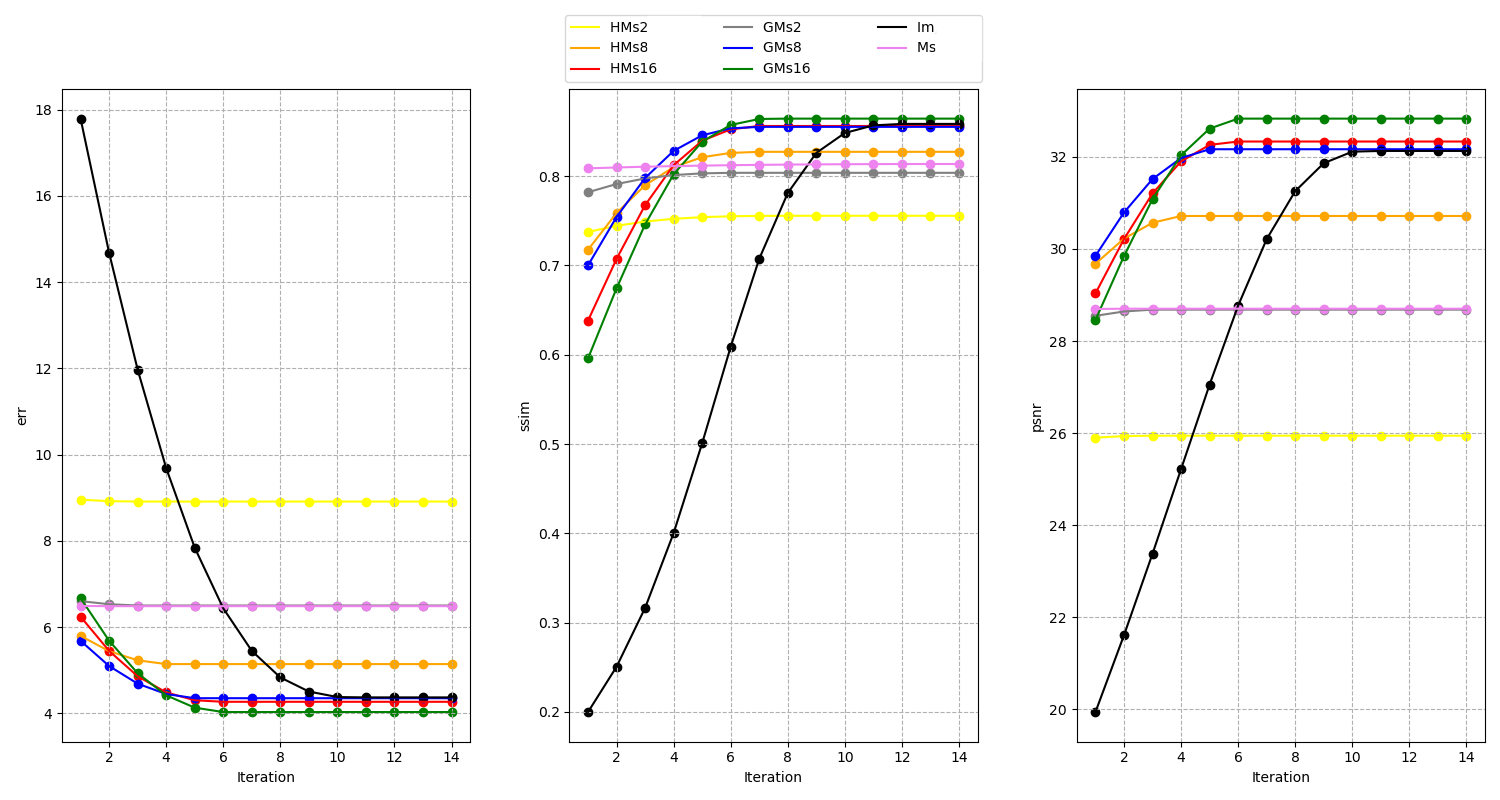}
\caption{Image 1 (House). 20 \% of noise}
\end{subfigure}
 \ \ \ \ \ \ \ \ 
\begin{subfigure}[b]{0.45\textwidth}  
\includegraphics[width=1\linewidth]{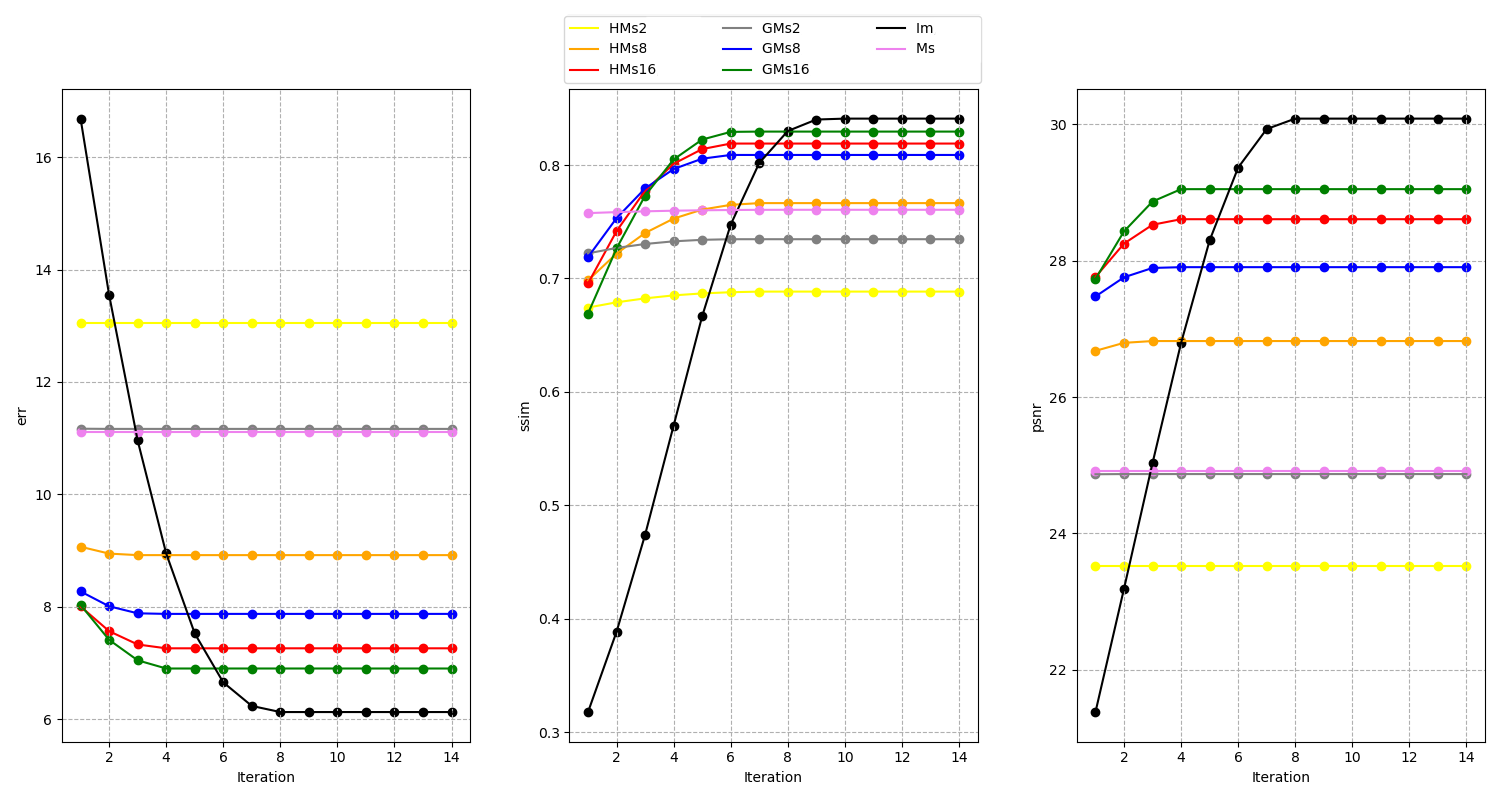}
\caption{Image 2 (Peppers). 20 \% of noise}
\end{subfigure}
\begin{subfigure}[b]{0.45\textwidth}  
\includegraphics[width=1\linewidth]{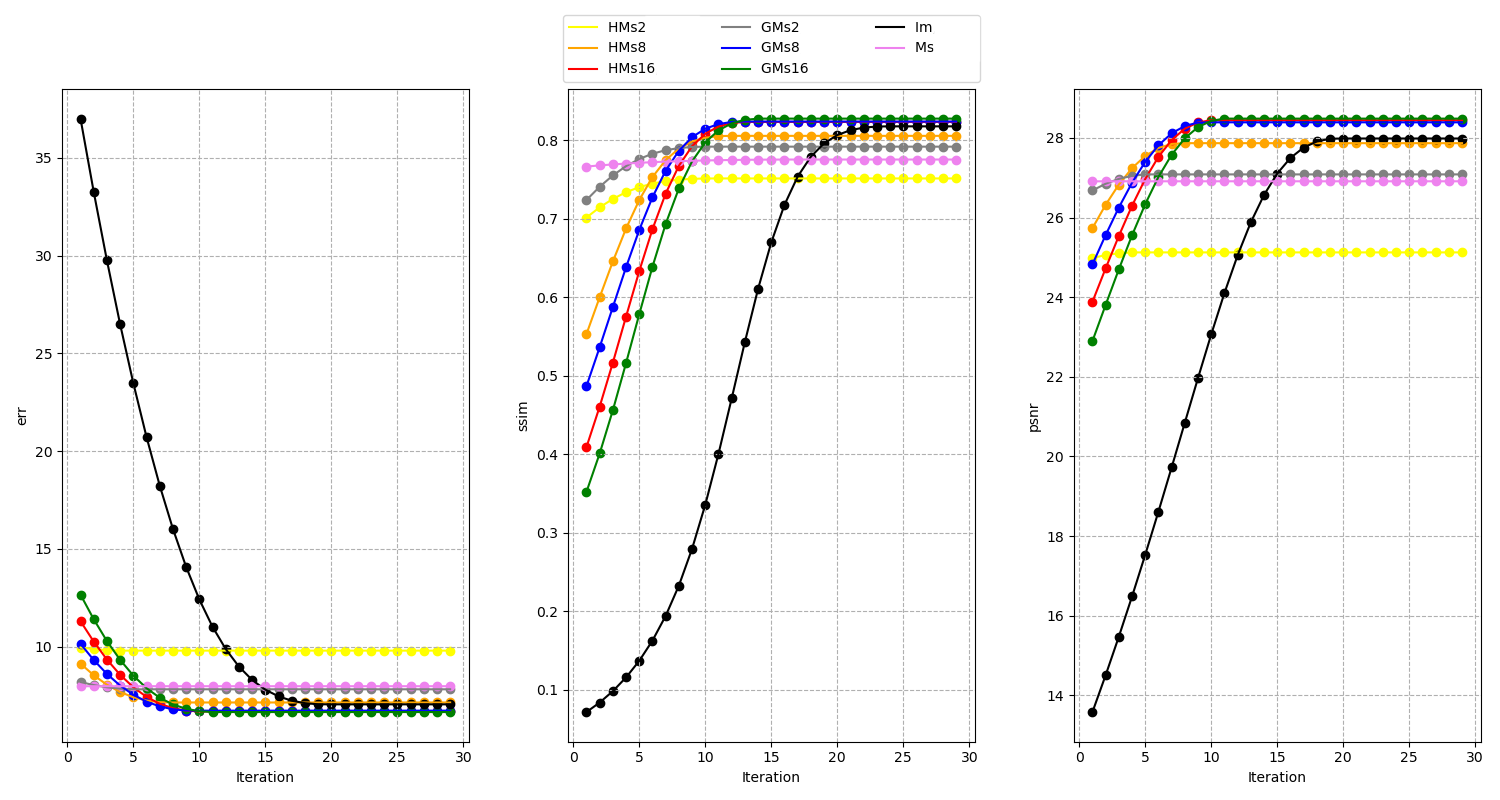}
\caption{Image 1 (House). 40 \% of noise.}
\end{subfigure}
 \ \ \ \ \ \ \ \ 
\begin{subfigure}[b]{0.45\textwidth}  
\includegraphics[width=1\linewidth]{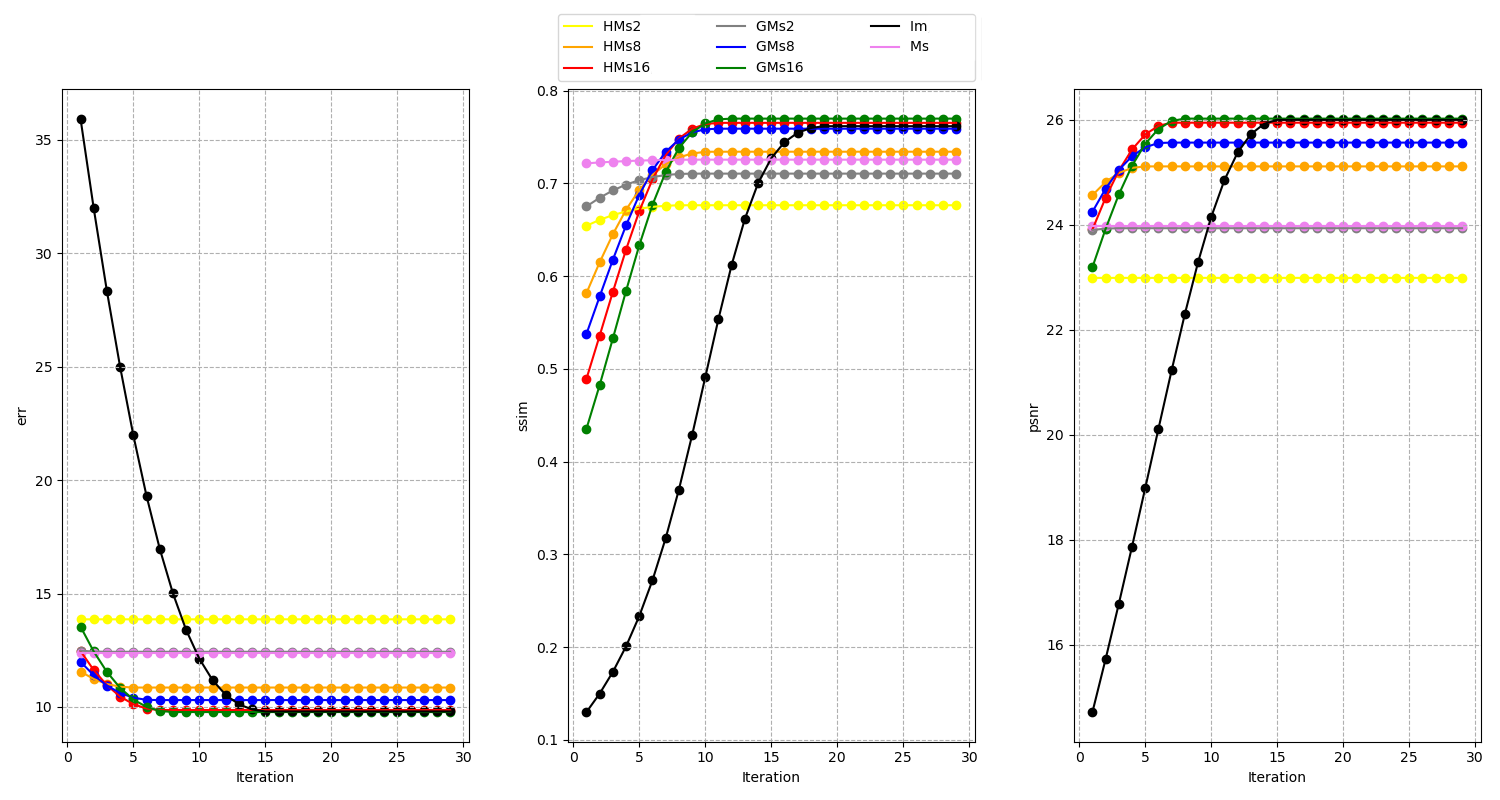}
\caption{Image 2 (Peppers). 40 \% of noise.}
\end{subfigure}
\caption{Test 3. Dynamics of the error and similarity for two images using different number of multiscale basis functions in AMs and GMs algorithms.  Relative $L_2$ error ($RRMSE$), similarity index ($SSIM$), peak signal-to-noise ratio ($PSNR$)}
\label{fig:t3-err}
\end{figure}

In the previous section, we proposed three multiscale algorithms: GMs, AMs, and Ms (MS$^0$). In GMs and Ms methods, we use linear local denoising to smooth (denoise) an initial local noised image, where the smoothed image is directly used for basis construction in the Ms algorithm, and local spectral problems with heterogeneous diffusion based on the denoised image is solved to enrich the coarse space with feature-preserving multiscale basis functions. In the AMs-algoritm, we use an analytical formula for eigenvectors relative to the linear denoising (Laplace) operator. 

In Tables \ref{table:t3-1} and \ref{table:t3-2}, we present the results of the denoising process for images with 10 and 40 \% of noise. To compare results, we show RRMSE, SSIM, and PSNR for an initial noised image, $I^0(x)$. Additionally, we add results for implicit (Im-ROF) and explicit (Ex-ROF) schemes using ROF nonlinear coefficient with 15 and 30 time iterations in the implicit scheme for $\tau = 3$ and explicit scheme with $\tau = 0.1$ leading to 450 and 900 time iterations. We investigate the influence of the coarse grid size on the performance of the multiscale methods. We consider $8 \times 8$ and $16 \times 16$ coarse grids with local cell size $N_h^{K_i} = 16 \times 16$ and $N_h^{K_i} = 32 \times 32$ leading to maximum of $N_h^{\omega_i} = 32 \times 32$ and  $N_h^{\omega_i} = 64 \times 64$ local domain size. Ms$^0$ and Ms methods use a local denoised image to define a multiscale basis function. In this method, we have one basis function per each local domain $\omega_i$, resulting in a coarse scale system with $DOF_H = N_c = 1089$ and $289$ for $32 \times 32$ and $16 \times 16$ coarse grids for $N_h^{K_i} = 16 \times 16$ and $N_h^{K_i} = 32 \times 32$, respectively. For the fine grid system on $512 \times 512$ fine grid, we have $DOF_h = 262,144$. Here, Ms$^0$ is related to projecting a noised image to a multiscale space without time iterations. In Ms-algorithms, we observe a good denoising behavior, and we see that the time iterations don't affect the denoising output. For the GMs and AMs algorithm, we present results with varying numbers of basis functions $M = 1, 2, 4, 8, 12$ and $16$ with corresponded coarse scale system size $DOF_H = N_c \times M$ ($N_c = 1089$ for $N_h^{\omega_i} = 16^2$ and $N_c = 289$ for $N_h^{\omega_i} = 32^2$). We observe better performance of the GMs approach due to the usage of the feature-preserved operator with local denoising. We note that the case with $M=1$ (AMs1 and GMs1) are identical because the first eigenvector for both cases is a constant vector. For the case with 20 \% of noise in Table \ref{table:t3-1}, we observe a good performance of the GMs algorithm with a sufficient number of basis functions. We also see that we have better denoising results with a finer coarse grid. For the case with a 40 \% error, we see an excellent performance of the GMs and AMs algorithm with a sufficient number of basis functions and can obtain the same results as a fine-grid solver. 

Figures \ref{fig:t3-noi40-1} and  \ref{fig:t3-noi40-2} depict the results of the denoising process to visually represent the method performance. We choose Image 1 (House) and Image 2 (Peppers) for illustration. We observe that using one basis function is insufficient to obtain good results in all methods. We improve denoising performance by enriching a multiscale space with spectral basis functions. Furthermore, we observe that the feature-preserved multiscale basis functions in the GMs method give visually better results, especially for the case with a smaller number of basis functions. Furthermore, we see how a choice of coarse grid resolution affects the denoising results. We observe that more basis functions are needed to encapsulate small-scale behavior in a multiscale space.  
We illustrate the dynamics of the denoising process for Image 1 (House) and Image 2 (Peppers) in Figure \ref{fig:t3-err}. From the presented results, we also observe the effect of the number of basis functions on the accurate image representation and convergence of the time iterations. 
We observe very slow dynamics using original fine-scale resolution and a fast error reduction for the multiscale method with a sufficient number of basis functions. We see that the first projection step significantly reduces the error of the initial noised image. We also know that we cannot obtain good results with a small number of basis functions because they cannot capture a fine-scale behavior, which clearly describes the poor performance of the Ms-method. Multiscale space enrichment with local spectral basis functions better encapsulate the multiscale feature. 
}


\subsection{Color images \rev{and high-resolution datasets (Test 4)}}

\rev{Finally, we consider color image denoising and study the effect of image resolution on algorithm performance}. We extend the proposed method to color images using YCrCb representation. The YCrCb color space separates image luminance (brightness) from chrominance (color) information, allowing for more efficient compression. Y is the luminance that represents the image's brightness, Cr is the chrominance red that represents the difference between the red component and the luminance, and Cb is chrominance blue that represents the difference between the blue component and the luminance.

\rev{
Numerical experiments are performed on classic and high-resolution datasets. 
We consider the following datasets:
\begin{itemize}
\item Set12(g): set contains 12 greyscale images from the Set12 benchmark dataset \cite{zhang2017beyond}.
Image resolution $256 \times 256$ pixels that give $N_h = 65,536$ of unknowns.

\item BSD68(g): set contains 68 greyscale images from the BSD68 dataset (Berkeley Segmentation Dataset) containing 68 natural images commonly used to evaluate image denoising algorithms \cite{martin2001database}. 
Image resolution $\approx 480 \times 320$ pixels gives $N_h \approx 153,000$ of unknowns.

\item Set14(c): set contains 14 color images form the Set14 benchmark dataset that commonly used to evaluate image algorithms \cite{zeyde2012single}. 
Image resolution $\approx 250 \times 250$ to $\approx 700 \times 600$ pixels that gives $N_h \approx 62,000$ to $420,000$ of unknowns.

\item BSDS100(c): set contains 100 color images from the BSDS100 (Berkeley Segmentation Dataset) dataset that is commonly used for evaluating image processing algorithms \cite{MartinFTM01}. 
Image resolution $\approx 500 \times 300$ pixels gives $N_h \approx 150,000$ of unknowns.

\item urban100(c): set contains 100 high-resolution (1k) color images from the Urban100 dataset that contains 100 images of urban scenes \cite{huang2015single}. 
Image resolution $\approx 600 \times 1,000$ pixels that give $N_h \approx 600,000$ of unknowns.

\item DIV2K(c): set contains 100 high-resolution color images from the DIVerse 2K resolution high-quality images \cite{Timofte_2017_CVPR_Workshops, Timofte_2018_CVPR_Workshops, Timofte_2018_CVPR_Workshops, Ignatov_2018_ECCV_Workshops}. 
Image resolution $\approx 1,500 \times 2,000$ pixels that give $N_h \approx 3,000,000$ of unknowns.

\item SIDD(c): set contains 10 high-resolution color images from the SIDD (Smartphone Image Denoising Dataset) benchmark that contains 10 scenes \cite{SIDD_2018_CVPR, Abdelhamed_2019_CVPR_Workshops}. 
Image resolution $\approx 3,000 \times 4,000$ pixels that give $N_h \approx 12,000,000$ of unknowns.
\end{itemize}

We start with five images from different datasets with different image resolutions to test the proposed methods' performance in the color image denoising process and investigate the method's performance depending on the size of the image. 
We consider the five test cases: Image 1 (Peppers) with $512 \times 512$ pixels resolution; Image 2 (Barbara) with $720 \times 576$ pixels resolution; Image 3 (Building) with $1,024 \times 768$ pixels resolution; Image 4 (Parrot) with $2,040 \times 1,536$ pixels resolution; and Image 5 (Markers) with $4,048 \times 3,044$ pixels resolution (see Table \ref{table:t4-time} for details). 
The images are presented in Figure \ref{fig:t4-ref}. We note that Image 1 and 2 are taken from the Set14, Image 3 from urban100 dataset, Image 4 from DIV2K and Image 5 from SIDD dataset. Next, we consider five image with different resolution to investigate methods performance in details and study an effect of the resolution to computational efficiency of the algorithms. 

\begin{figure}[h!]
\centering
\includegraphics[width=0.145\linewidth]{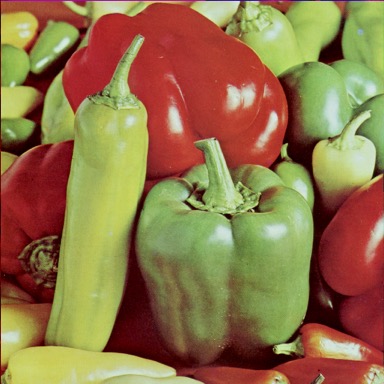} \ 
\includegraphics[width=0.182\linewidth]{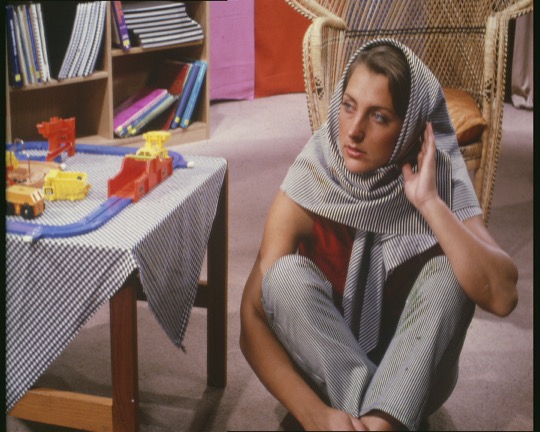}  \ 
\includegraphics[width=0.192\linewidth]{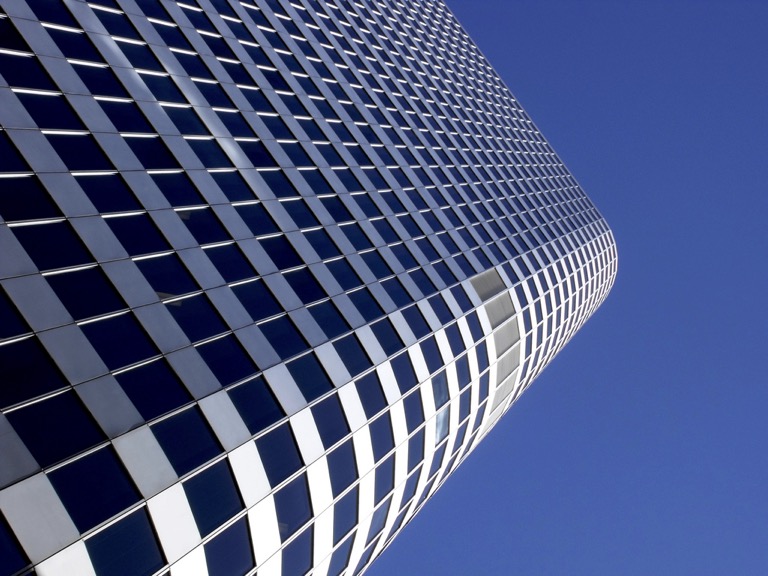} \
\includegraphics[width=0.218\linewidth]{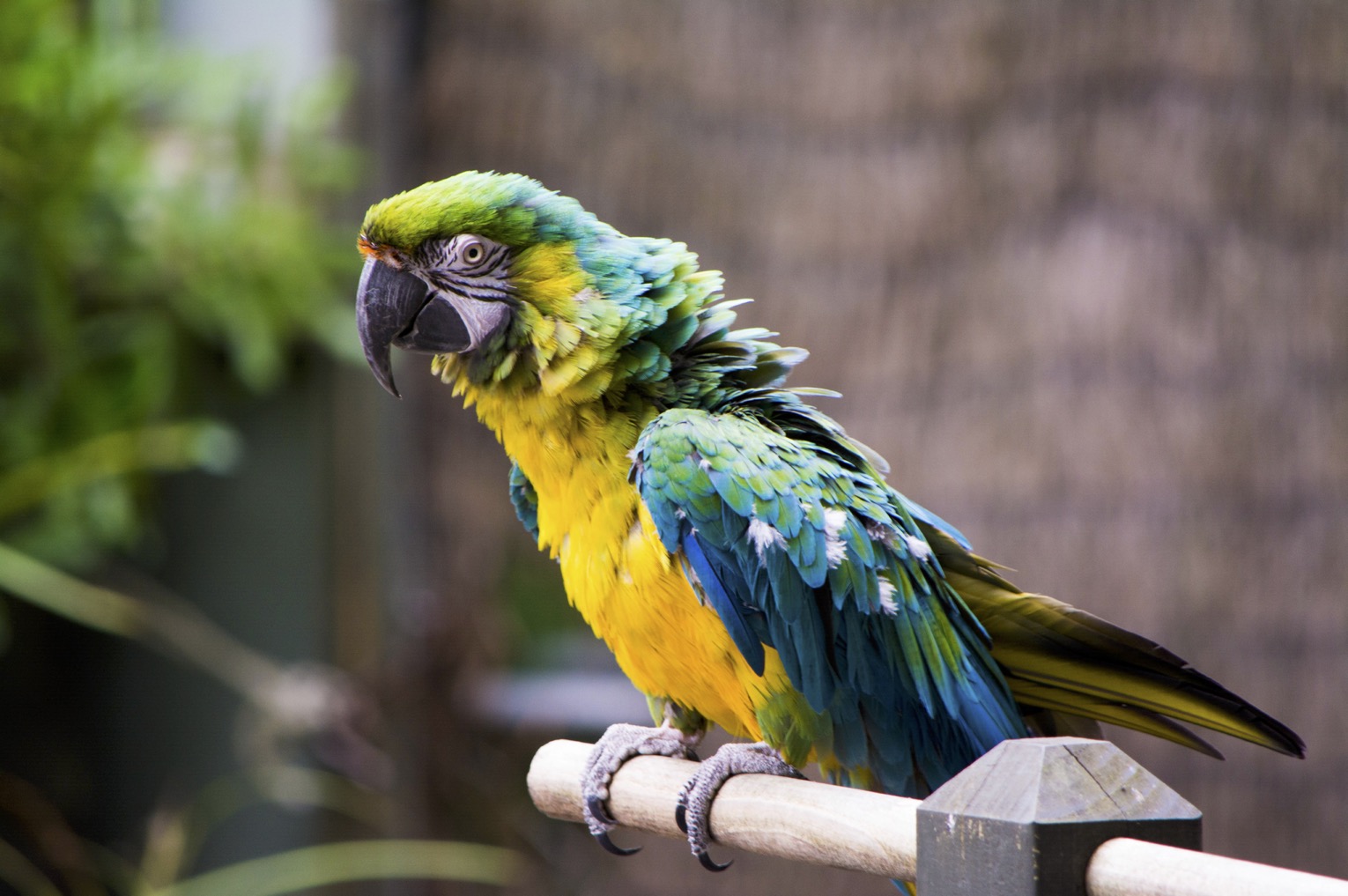} \ 
\includegraphics[width=0.193\linewidth]{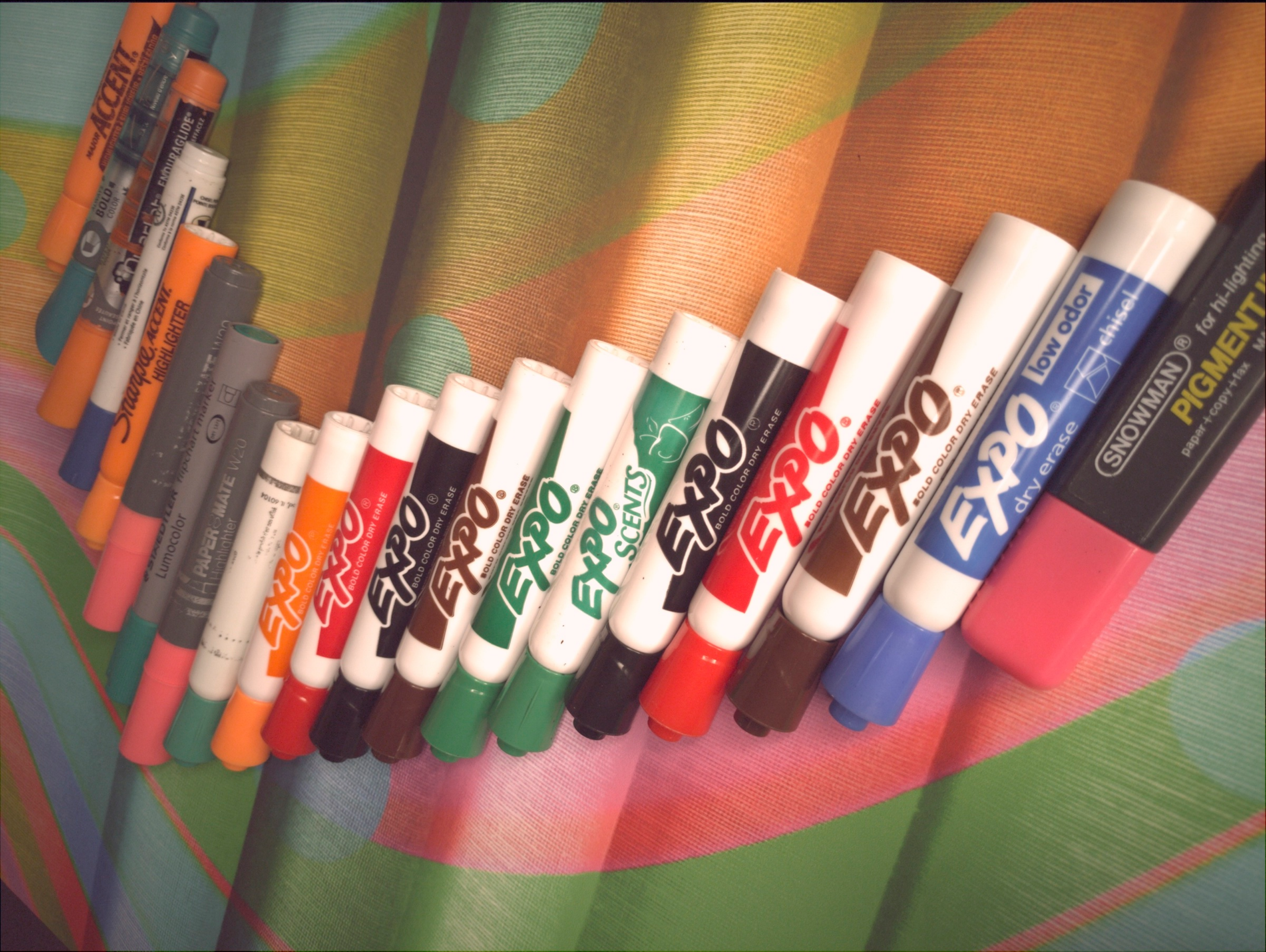}
\caption{\rev{Test 4. Reference image. Image 1 (Peppers), Image 2 (Barbara), Image 3 (Building), Image 4 (Parrot) and Image 5 (Markers) (from left to right) }}
\label{fig:t4-ref}
\end{figure}

\begin{table}[h!]
\begin{tabular}{| c | c | c | c | c | c | }
\hline
 & \tiny{RRMSE/SSIM/PSNR} & \tiny{RRMSE/SSIM/PSNR} & \tiny{RRMSE/SSIM/PSNR} & \tiny{RRMSE/SSIM/PSNR} & \tiny{RRMSE/SSIM/PSNR}\\
\hline
& 1.Peppers & 2.Barbara & 3.Building & 4.Parrot & 5.Markers \\
\hline
\multicolumn{6}{|c|}{Initial noised image with 20\% of noise}\\
\hline
$I_0(x)$  & 20.55/0.19/19.60 & 20.50/0.21/19.72 & 20.49/0.19/19.80 & 20.41/0.13/20.19 & 20.71/0.11/19.33 \\ 
\hline 
BM3D & 4.90/0.84/32.07 & 4.36/0.89/33.18 & 5.01/0.92/32.03 & 3.44/0.93/35.67 & 3.47/0.91/34.84 \\ 
Ex-ROF & 6.11/0.81/30.14 & 7.39/0.81/28.58 & 7.95/0.87/28.03 & 4.23/0.92/33.87 & 3.96/0.90/33.71 \\ 
Im-ROF & 6.17/0.81/30.05 & 6.98/0.83/29.09 & 8.15/0.88/27.82 & 4.28/0.93/33.76 & -/-/- \\
\hline
AMs1 & 11.60/0.70/24.57 & 11.72/0.71/24.58 & 20.01/0.76/20.01 & 7.76/0.88/28.59 & 5.74/0.87/30.48 \\ 
AMs2 & 9.79/0.74/26.05 & 10.57/0.73/25.48 & 17.32/0.77/21.27 & 6.97/0.89/29.52 & 5.12/0.88/31.47 \\ 
AMs4 & 8.36/0.76/27.42 & 9.65/0.75/26.27 & 14.92/0.79/22.56 & 5.61/0.90/31.41 & 4.49/0.89/32.62 \\ 
\hline
GMs1 & 11.60/0.70/24.57 & 11.72/0.71/24.58 & 20.01/0.76/20.01 & 7.76/0.88/28.59 & 5.74/0.87/30.48 \\ 
GMs2 & 8.59/0.76/27.18 & 9.66/0.76/26.26 & 15.78/0.79/22.07 & 6.53/0.90/30.08 & 4.61/0.89/32.37 \\ 
GMs4 & 7.57/0.78/28.27 & 8.92/0.78/26.95 & 13.82/0.81/23.23 & 5.26/0.91/31.96 & 4.28/0.90/33.03 \\ 
\hline
\multicolumn{6}{|c|}{Initial noised image with 40\% of noise}\\
\hline
$I_0(x)$  & 40.47/0.06/13.72 & 40.15/0.08/13.89 & 40.07/0.08/13.98 & 40.48/0.04/14.24 & 40.21/0.03/13.57 \\ 
\hline
BM3D & 11.24/0.77/24.84 & 9.41/0.80/26.49 & 12.11/0.84/24.37 & 8.35/0.88/27.95 & 8.42/0.87/27.15 \\ 
Ex-ROF & 12.80/0.74/23.72 & 11.46/0.75/24.78 & 15.31/0.81/22.34 & 9.30/0.88/27.02 & 9.09/0.87/26.49 \\ 
Im-ROF & 13.00/0.74/23.58 & 11.59/0.75/24.68 & 15.39/0.82/22.29 & 9.36/0.89/26.96 & -/-/- \\ 
\hline
AMs1 & 14.99/0.69/22.35 & 13.41/0.71/23.41 & 20.94/0.76/19.62 & 10.48/0.87/25.98 & 9.58/0.87/26.02 \\ 
AMs2 & 13.94/0.72/22.97 & 12.71/0.72/23.88 & 19.27/0.77/20.34 & 10.25/0.88/26.17 & 9.65/0.87/25.96 \\ 
AMs4 & 13.40/0.73/23.32 & 12.26/0.73/24.19 & 17.90/0.78/20.98 & 9.69/0.88/26.66 & 9.72/0.87/25.90 \\ 
\hline
GMs1 & 14.99/0.69/22.35 & 13.41/0.71/23.41 & 20.94/0.76/19.62 & 10.48/0.87/25.98 & 9.58/0.87/26.02 \\ 
GMs2 & 13.38/0.73/23.33 & 12.21/0.73/24.22 & 18.16/0.78/20.85 & 9.80/0.88/26.56 & 9.22/0.88/26.36 \\ 
GMs4 & 13.10/0.74/23.51 & 11.88/0.74/24.47 & 17.16/0.79/21.35 & 9.55/0.89/26.78 & 9.13/0.88/26.45 \\ 
\hline
\end{tabular}
\caption{\rev{Test 4. Color image denoising.  Image 1 (Peppers), Image 2 (Barbara), Image 3 (Building), Image 4 (Parrot) and Image 5 (Markers) (from left to right)}}
\label{table:t4-c}
\end{table}

\begin{figure}[h!]
\centering
\begin{subfigure}[b]{1\textwidth}
\includegraphics[width=0.145\linewidth]{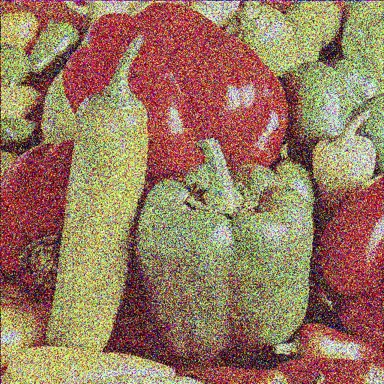} \ \ 
\includegraphics[width=0.182\linewidth]{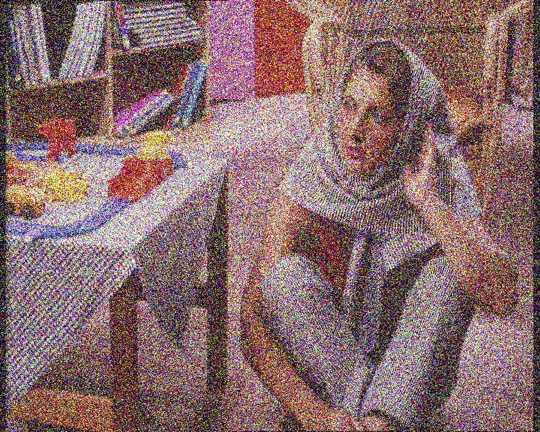}  \ \ 
\includegraphics[width=0.192\linewidth]{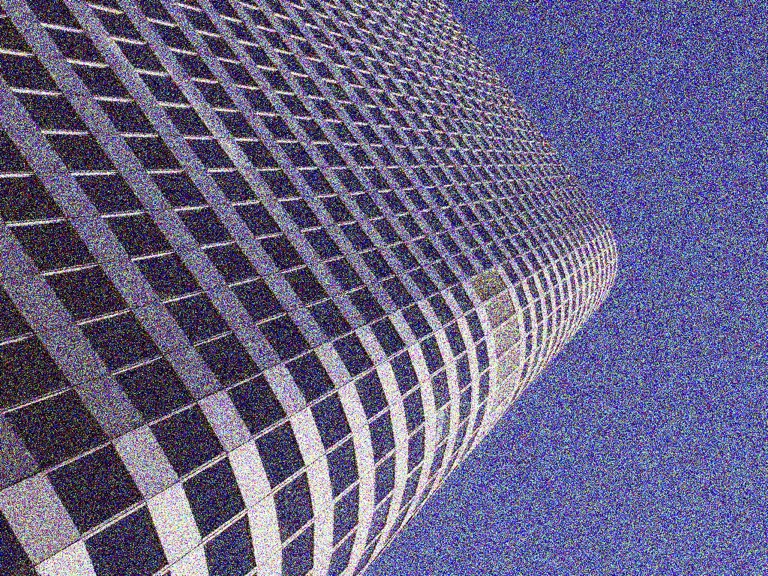} \ \ 
\includegraphics[width=0.218\linewidth]{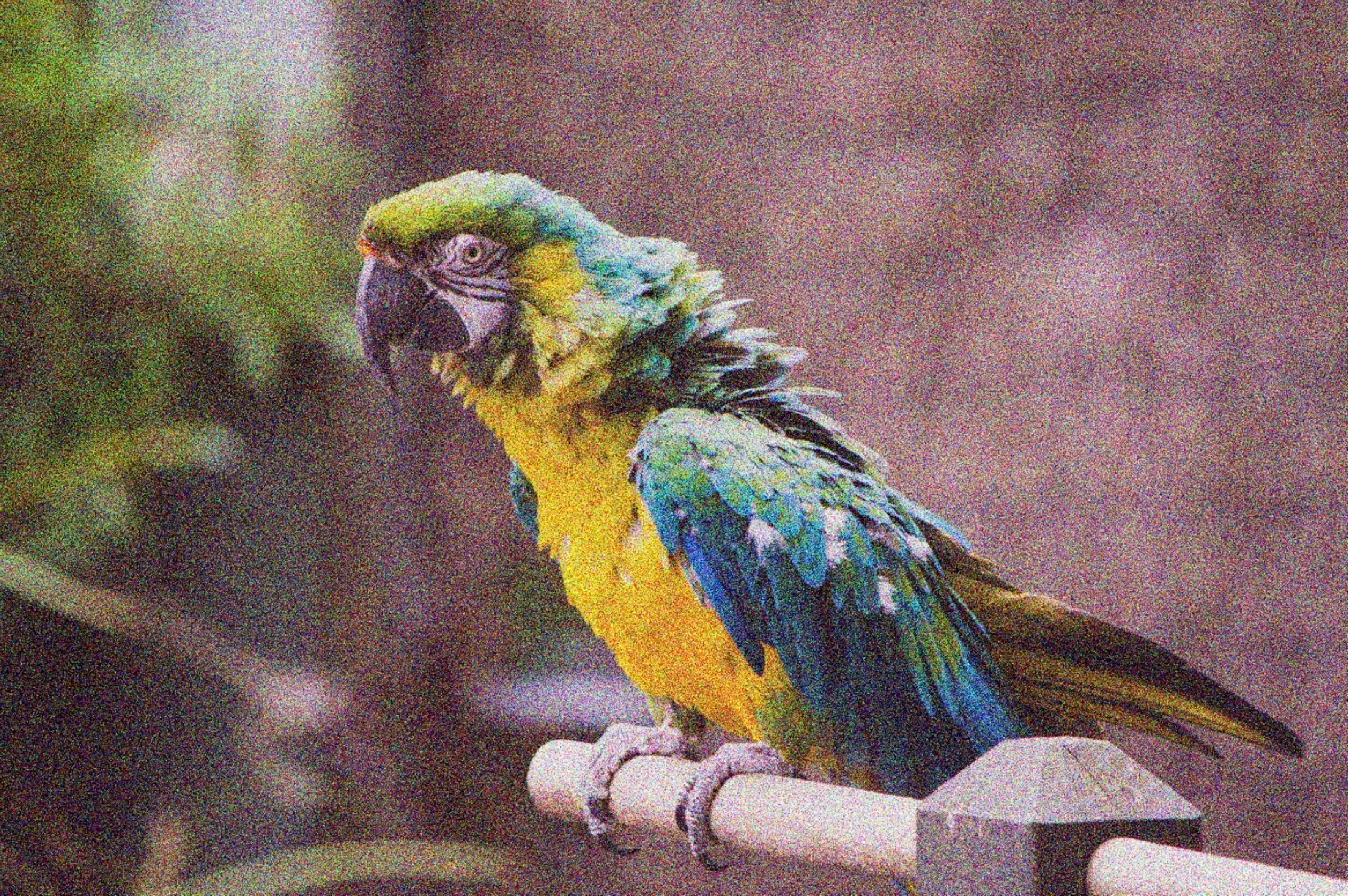} \ \ 
\includegraphics[width=0.193\linewidth]{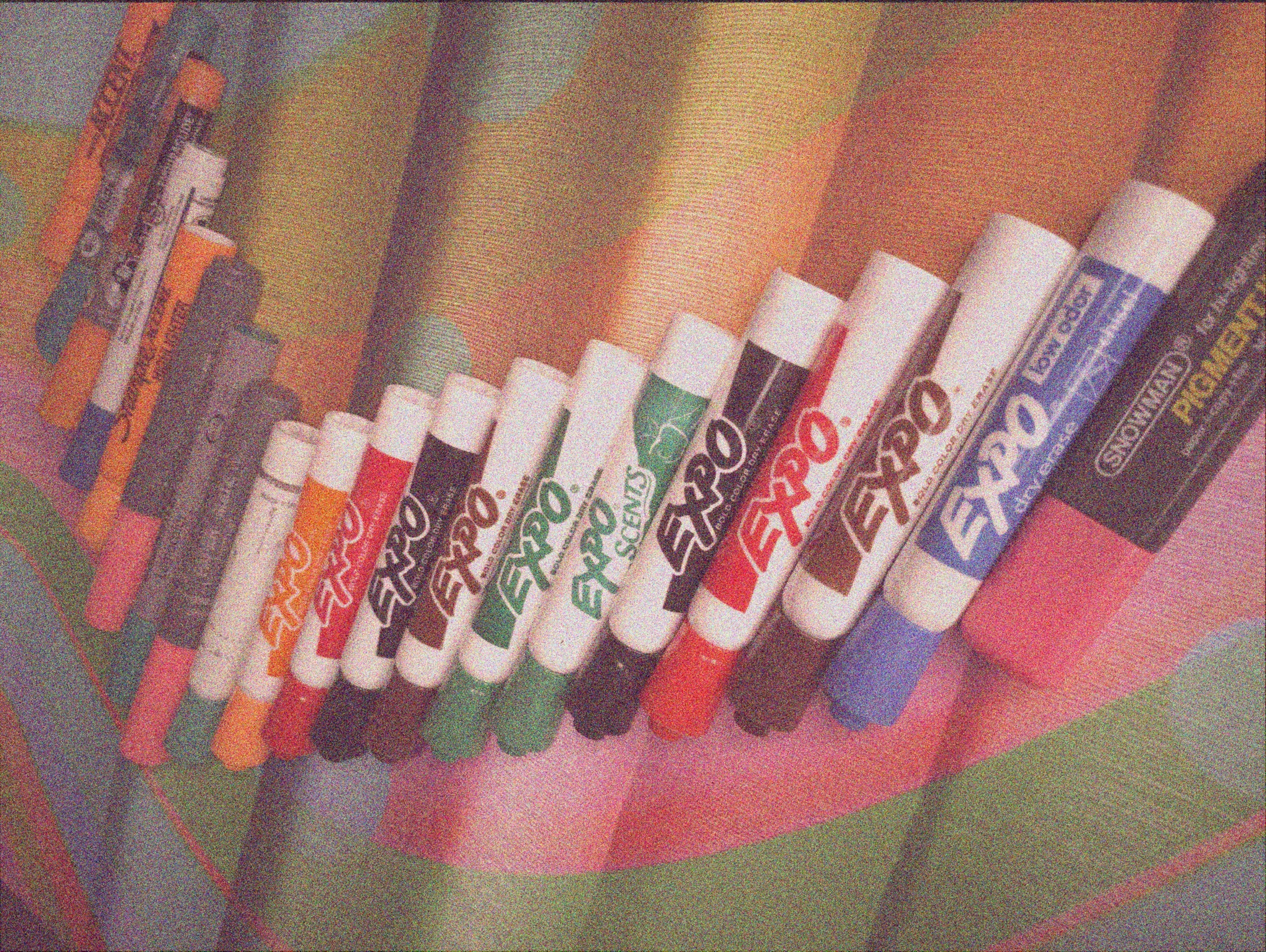}
\caption{Initial noised image, 40 \%}
\end{subfigure}
\begin{subfigure}[b]{1\textwidth}
\includegraphics[width=0.145\linewidth]{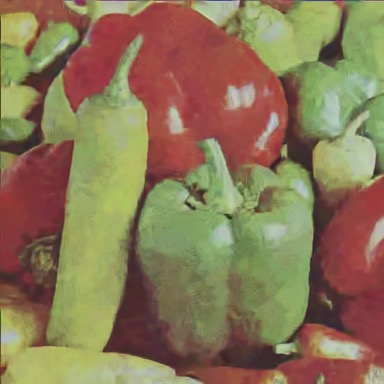} \ \ 
\includegraphics[width=0.182\linewidth]{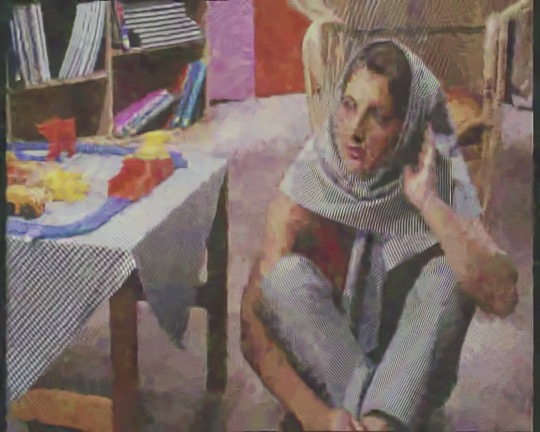}  \ \ 
\includegraphics[width=0.192\linewidth]{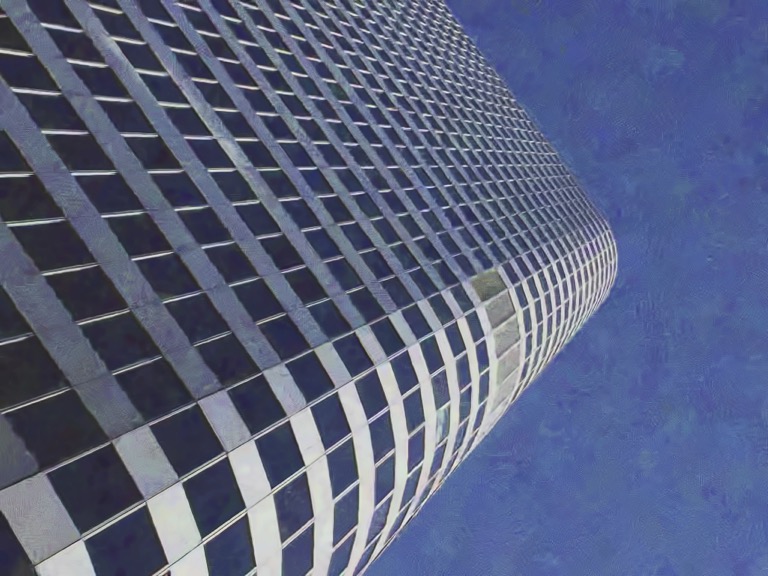} \ \ 
\includegraphics[width=0.218\linewidth]{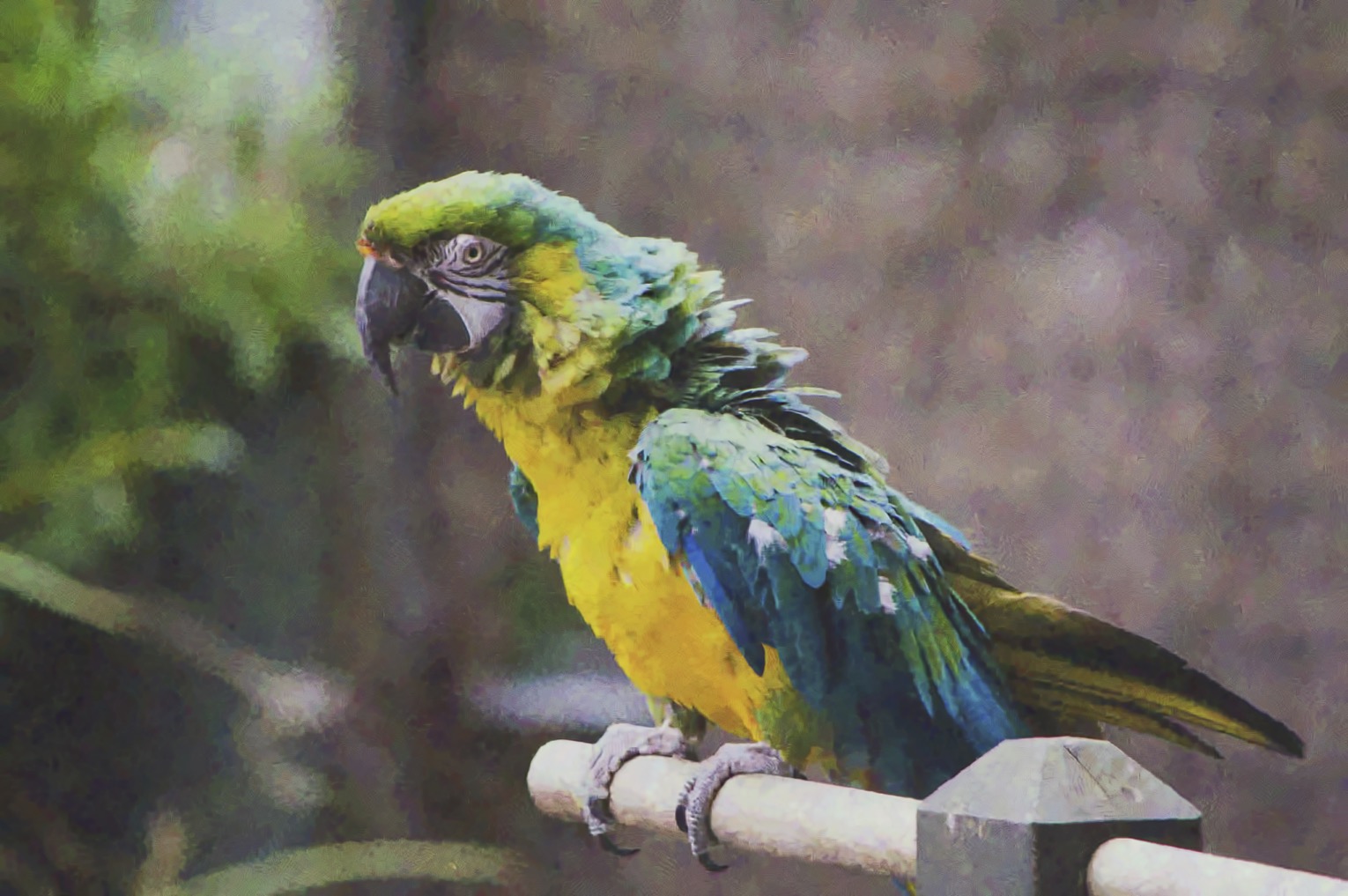} \ \ 
\includegraphics[width=0.193\linewidth]{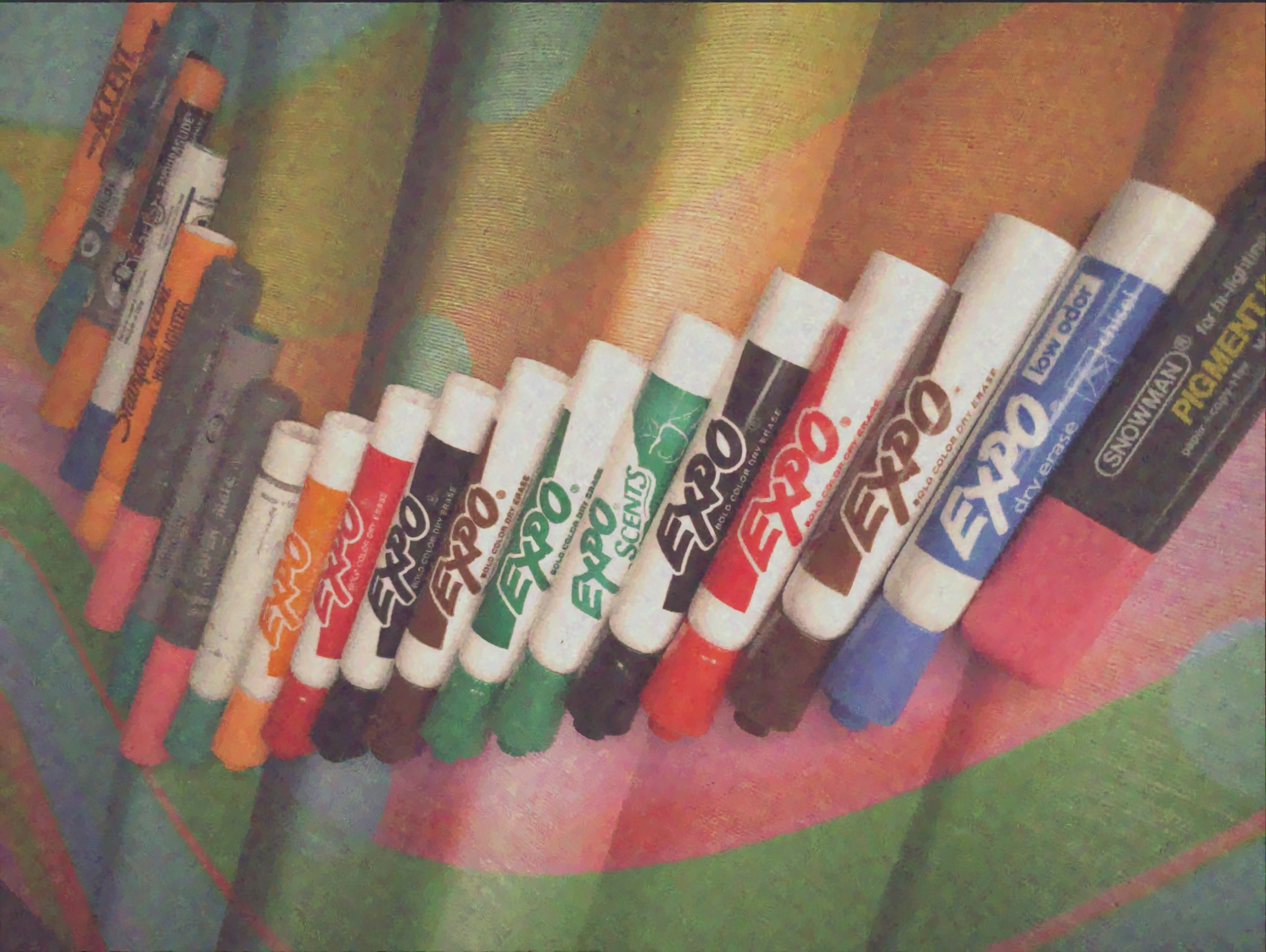}
\caption{BM3D}
\end{subfigure}
\begin{subfigure}[b]{1\textwidth}
\includegraphics[width=0.145\linewidth]{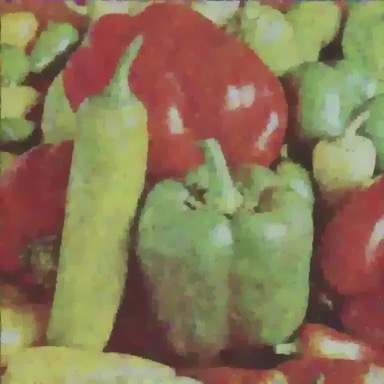} \ \ 
\includegraphics[width=0.182\linewidth]{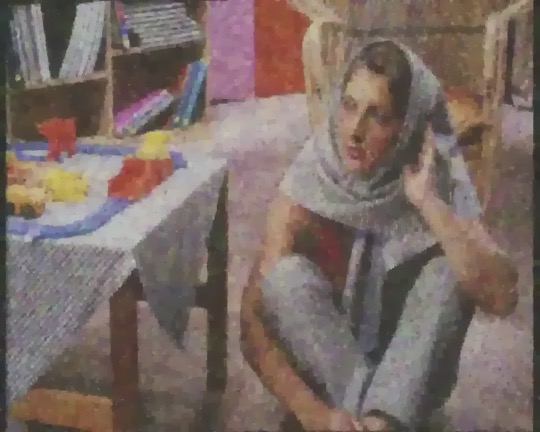}  \ \ 
\includegraphics[width=0.192\linewidth]{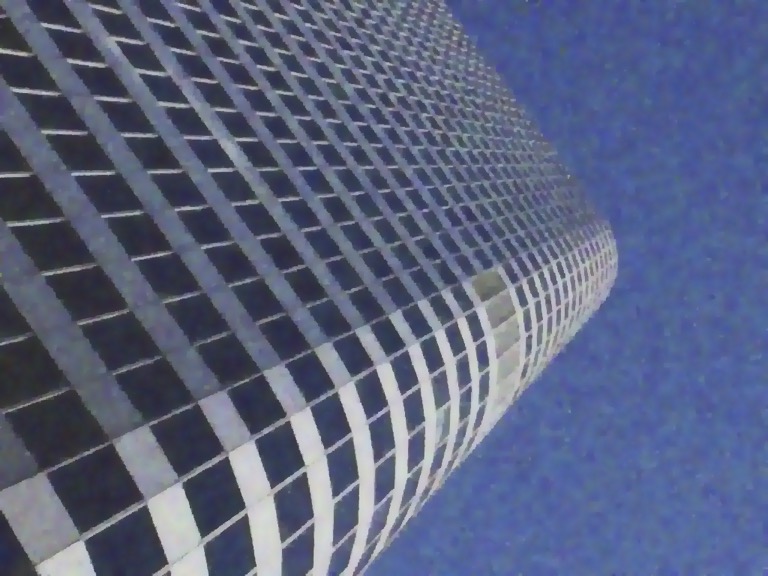} \ \ 
\includegraphics[width=0.218\linewidth]{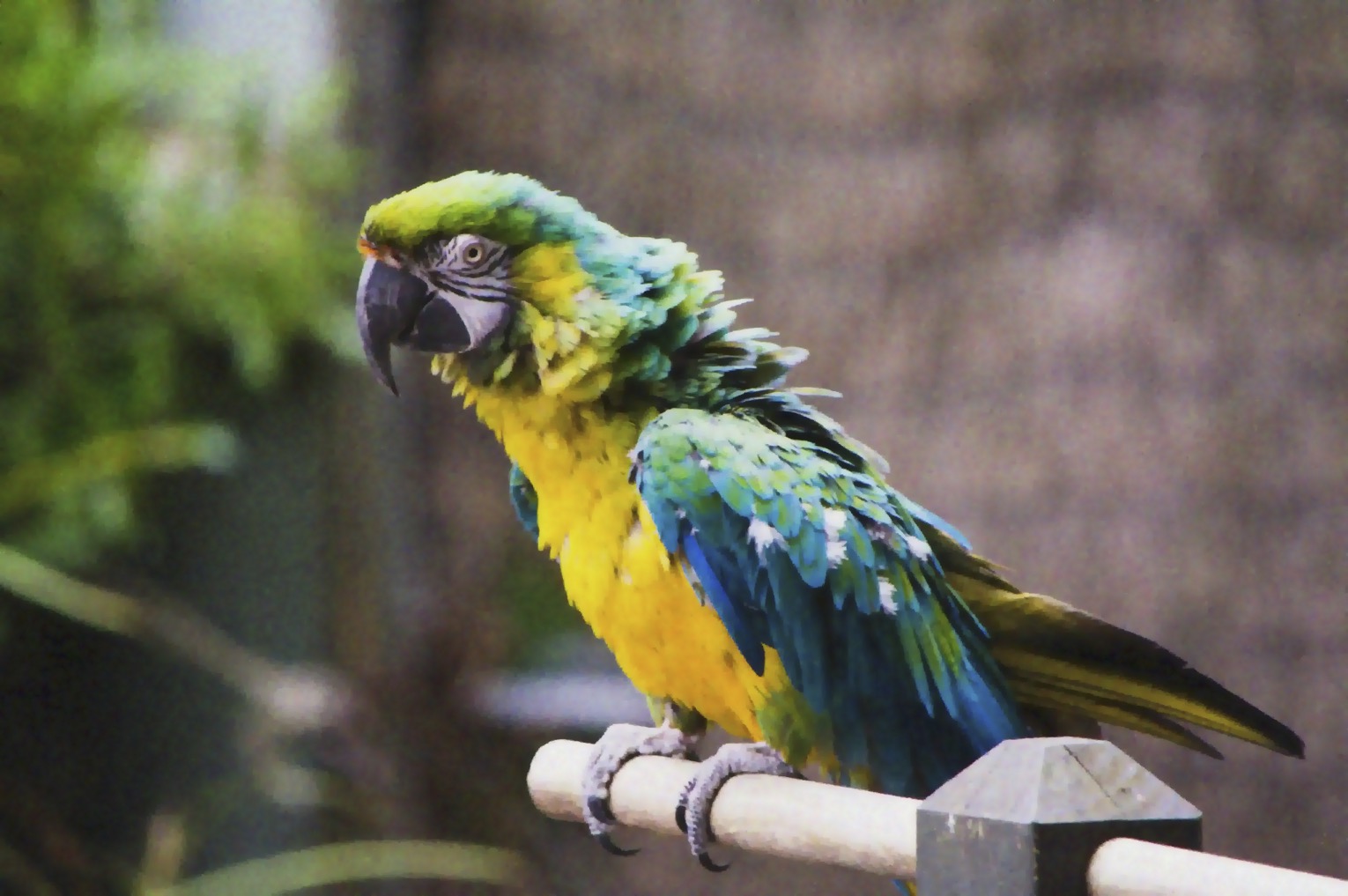} \ \ 
\includegraphics[width=0.193\linewidth]{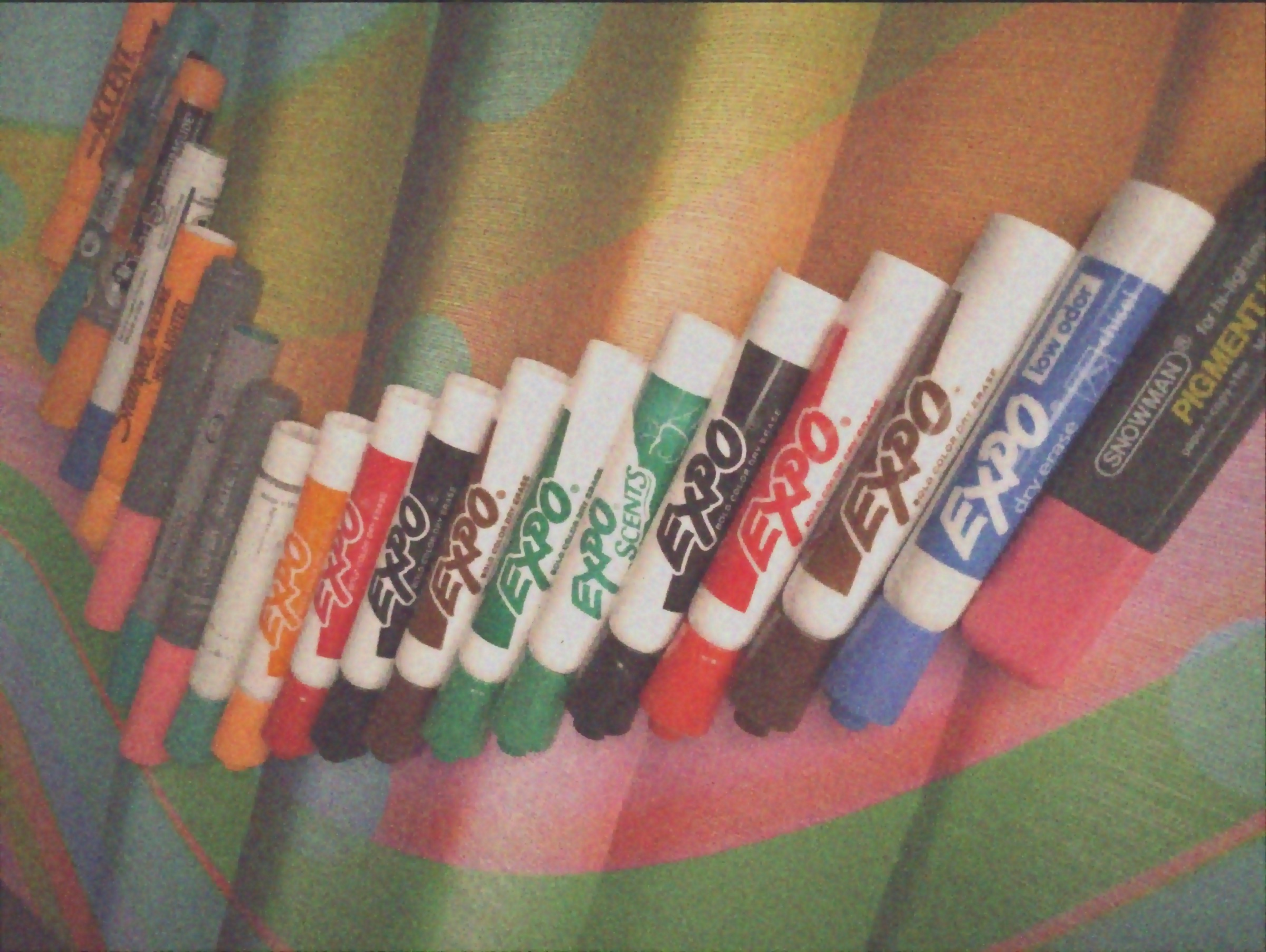}
\caption{Ex-ROF}
\end{subfigure}
\begin{subfigure}[b]{1\textwidth}
\includegraphics[page=5,width=0.145\linewidth]{img-c-2-t2-4-nl3-40} \ \ 
\includegraphics[page=1,width=0.182\linewidth]{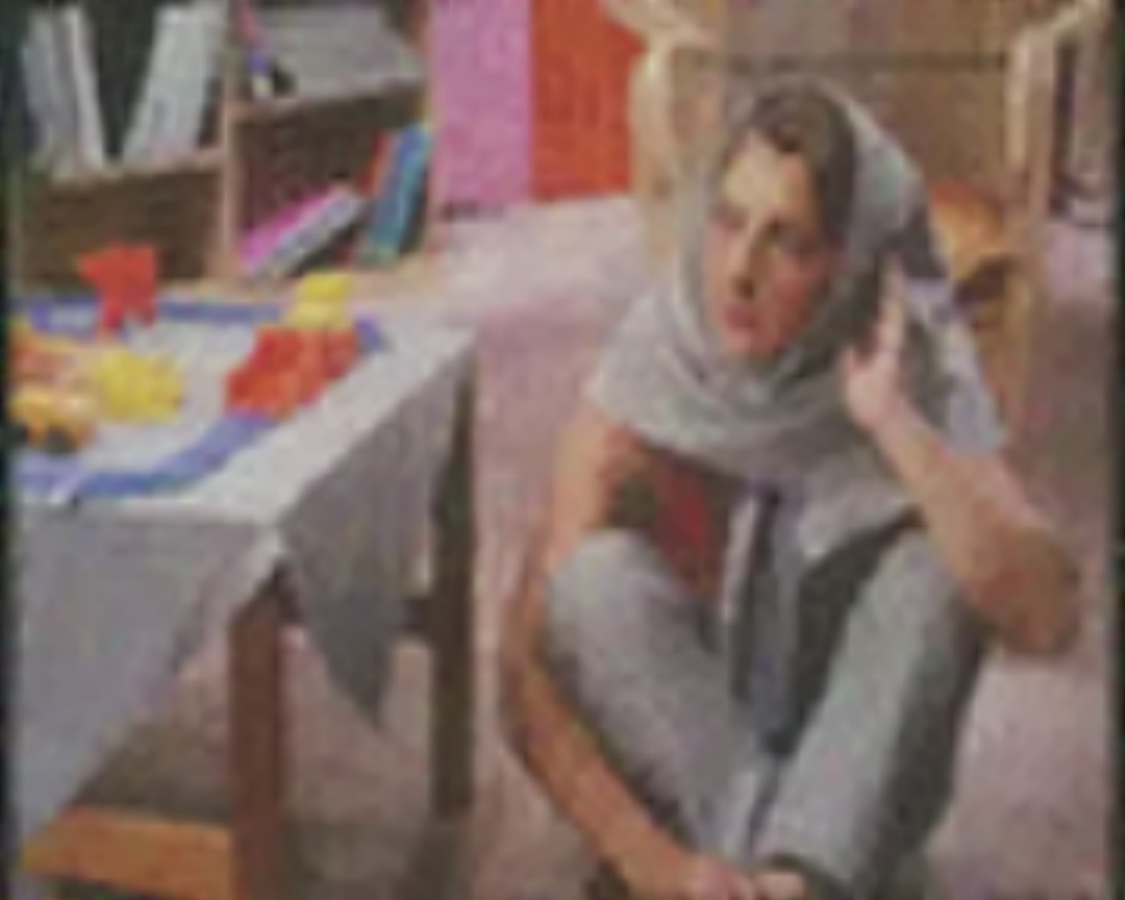}  \ \ 
\includegraphics[page=2,width=0.192\linewidth]{img-c-2-t2-4-nl3-40} \ \ 
\includegraphics[page=4,width=0.218\linewidth]{img-c-2-t2-4-nl3-40} \ \ 
\includegraphics[page=3,width=0.193\linewidth]{img-c-2-t2-4-nl3-40}
\caption{AMs4-Im}
\end{subfigure}
\begin{subfigure}[b]{1\textwidth}
\includegraphics[page=5,width=0.145\linewidth]{img-c-2-t3-4-nl3-40} \ \ 
\includegraphics[page=1,width=0.182\linewidth]{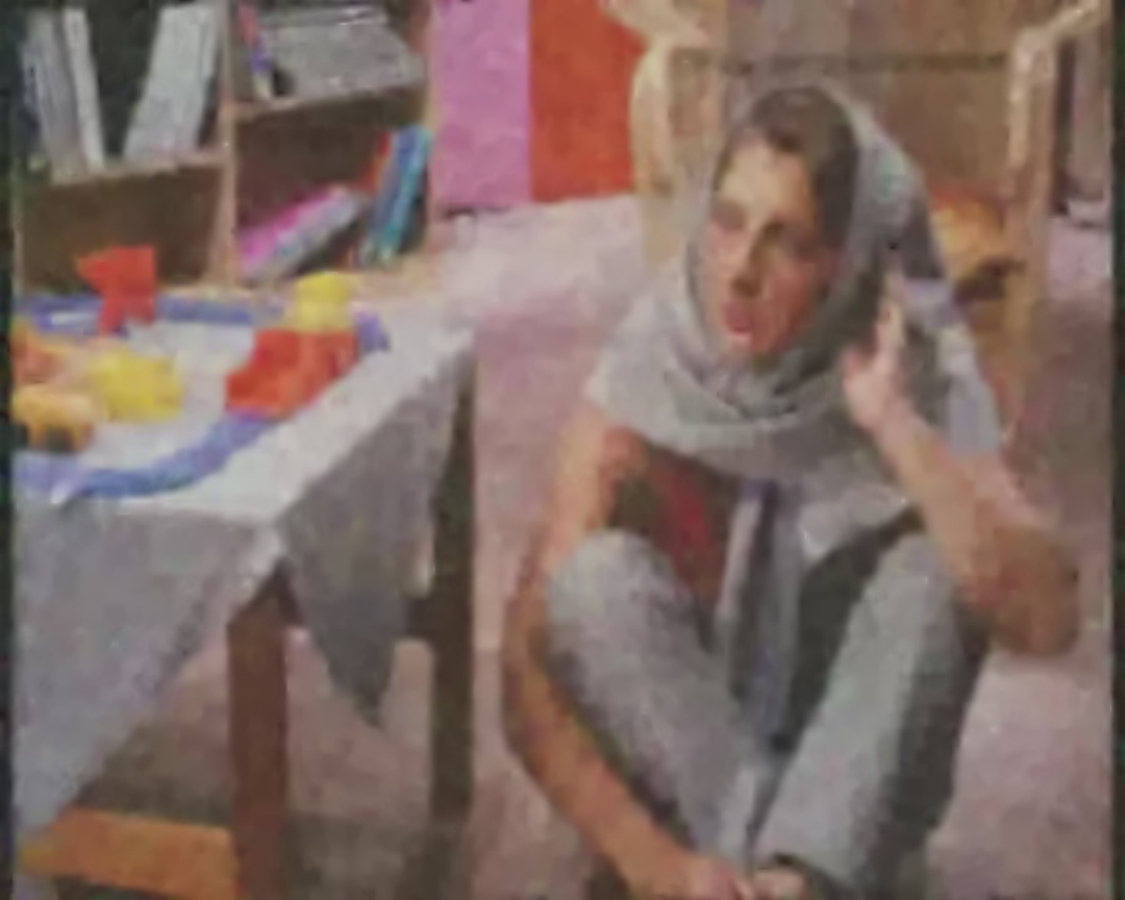}  \ \ 
\includegraphics[page=2,width=0.192\linewidth]{img-c-2-t3-4-nl3-40} \ \ 
\includegraphics[page=4,width=0.218\linewidth]{img-c-2-t3-4-nl3-40} \ \ 
\includegraphics[page=3,width=0.193\linewidth]{img-c-2-t3-4-nl3-40}
\caption{GMs4-Im}
\end{subfigure}
\caption{\rev{Test 4. Color image denoising. Image 1 (Peppers), Image 2 (Barbara), Image 3 (Building), Image 4 (Parrot) and Image 5 (Markers) (from left to right). 40 \% of noise}}
\label{fig:t4-noi40}
\end{figure}

\begin{table}[h!]
\centering
\begin{tabular}{| c | c | c | c | c | c | }
\hline
& 1.Peppers & 2.Barbara & 3.Building & 4.Parrot & 5.Markers \\
\hline
$N_1 \times N_2$  & 512$\times$512  & 720$\times$576  & 1,024$\times$768  & 2,040$\times$1536  & 4,048$\times$3,044 \\ 
$DOF_h$ & 262,144 & 414,720 & 786,432 & 3,133,440 & 12,322,112 \\
$N_c$ & 1,089 & 1,702 & 3,185 & 11,094 & 48,769 \\
\hline
\hline
 & {DOF/$t_{run}$(sec)}  & {DOF/$t_{run}$(sec)}  & {DOF/$t_{run}$(sec)}  & {DOF/$t_{run}$(sec)}  & {DOF/$t_{run}$(sec)}\\
\hline
BM3D & -/6.1 & -/8.7 & -/14.9 & -/49.0 & -/215.2 \\ 
Ex-ROF & 262,144/3.8 & 414,720/9.6 & 786,432/15.0 & 3,133,440/69.8 & 12,322,112/290.7 \\
Im-ROF & 262,144/8.1 & 414,720/13.8 & 786,432/36.2 & 3,133,440/459.5 & 12,322,112/- \\
\hline
AMs1 & 1,089/0.7 & 1,702/1.3 & 3,185/2.3 & 11,094/8.4 & 48,769/65.4 \\ 
AMs2 & 2,178/2.1 & 3,404/2.9 & 6,370/5.4 & 22,188/26.1 & 97,536/195.2 \\
AMs4 & 4,356/5.2 & 6,808/7.6 & 12,740/18.5 & 44,376/108.8 & 195,072/416.6 \\ 
\hline
\end{tabular}
\caption{\rev{Test 4. Running time for color image denoising with various resolution. 40 \% of noise}}
\label{table:t4-time}
\end{table}

}

\rev{
In order to extend the multiscale method to a three-channel (color) image, we use the following approach: (Step 1) Solve time-dependent problem \eqref{eq:fine} on a fine grid (full resolution) using the implicit or explicit scheme for each channel Y, Cr, and Cb separately; and (Step 2) Combine solutions for each channel and save the solution.
For the multiscale approach, we have similar steps: (Step 1) Construct a multiscale basis function for a one-channel image representation (grayscale) for a given coarse grid size $N_c$l (Step 2) Construct a global projection matrix $R$ and solve time-dependent problem \eqref{eq:coarsep} on a coarse grid using the multiscale method for each channel Y, Cr, and Cb separately; and (Step 3) Combine multiscale solutions for each channel and save the solution.
}
For multiscale simulations using the proposed method, we generate a coarse grid of the following sizes to have coarse cells with the resolution \rev{$N_h^{K_i} = 16^2$ pixels.} Note that, the local calculations can be done in parallel. The time of the local solution depends on the size of the local domain. The coarse grid should be chosen to balance local calculations and the compression ability of the multiscale coarse-grid representation.

\rev{In Table \ref{table:t4-c}, we present convergence results for the multiscale method and fine-scale denoising processes. In Table \ref{table:t4-time}, we present the corresponding size of the fine-scale $DOF_h$, reduced coarse-scale model $DOF_H$ and running of calculations in seconds. Figure \ref{fig:t4-noi40} gives a visual representation of the results. We use explicit and implicit schemes (Ex-ROF and Im-ROF) and implicit schemes with multiscale solvers based on AMs and GMs methods with four basis functions (AMs4-Im and GMs4-Im). Similarly to previous tests, we consider explicit scheme (Ex-ROF) with 180 and 360 time iterations with $\tau = 0.25$ for 20 and 40 \% of noise, respectively. In implicit solving with and without multiscale methods, we use 7 and 15-time iterations with $\tau = 3$ for 20 and 40 \% of noise. Additionally, we present results for the BM3D method. In multiscale methods, we use about three times fewer iterations than in the fine-scale scheme IM-ROF (3 and 5 iterations). 
We observed good performance for fine-scale solvers (Ex-ROF and Im-ROF) with almost the same results. However, we could not solve a fine grid system in the Im-ROF method due to the very large size of the system. As we mentioned before, the size of the fine grid system is equal to the number of image pixels $DOF_h = N_1 \times N_2$, resulting in 12.3 million unknowns in the last Image 5 with 4k resolution. This results in the most expensive calculations in all methods. By comparing the BM3D method and our PDE-based approach, we observe similar denoising results for images with large errors. The time of calculations of BM3D and Ex-ROF methods are similar. We have slower calculations for the implicit scheme even with a much larger time step size. We note that we used a direct solver from scipy.sparse library to solve an underlying linear system. Further investigation of the preconditioned iterative methods and parallel solver should be performed in future works. 
We obtain a faster algorithm by utilizing the multiscale method for approximating and solving the system from implicit time approximation. We observe a good denoising performance. However, the results are slightly worse than for the fine-scale solver (Im-ROF).  
} 
Note that because we consider the decoupled process for each channel, we use the same values of the parameters for each channel. The multiscale basis functions are calculated for a grayscale version of the image and used for each Y, Cr, and Cb channel. Denoising of each channel is done separately. In general, we can consider a color image as a 3D array and construct multiscale basis functions that incorporate eigenvectors of a color image. However, in this work, we concentrate on the simplified decoupled version of the image into separate channels that do not interact. \rev{Further improvement can be made by considering a coupled PDE-based approach without color channel separation on fine and multiscale solvers.} 
We also note that the time of calculations presented in tables for the multiscale method does not include the time of multiscale basis construction that can be done in parallel. In future work, we will consider hybrid approaches to reduce the time of basis construction, for example, using machine learning techniques \cite{vasilyeva2019convolutional, vasilyeva2020learning, vasilyeva2021machine}. \rev{Moreover, we used a projection approach and constructed a global projection matrix using a multiscale method to project a system onto the course grid. Further optimization of the computational efficiency should be done by efficiently constructing the coarse scale system locally in each coarse cell, especially for the AMs method. We may obtain significantly faster implementation. We will consider it in future works in a couple with optimizing the multiscale basis construction process.}

\begin{table}[h!]
\centering
\begin{tabular}{| c | c |  c | c | c | }
\hline
 & \tiny{RRMSE/SSIM/PSNR} & \tiny{RRMSE/SSIM/PSNR} & \tiny{RRMSE/SSIM/PSNR} & \tiny{RRMSE/SSIM/PSNR}\\
 \hline
& $I^0(x)$ & Ex-ROF  & AMs4-Im  & BM3D \\
\hline
\multicolumn{5}{|c|}{Initial noised image with 20\% of noise}\\
\hline
Set12(g) 	& 20.03/0.37/19.61 & 7.51/0.83/28.14 & 7.79/0.82/27.81 & 6.18/0.83/29.82 \\
BSD68(g) 	& 21.20/0.30/19.65 & 8.64/0.80/27.45 & 8.89/0.79/27.20 & 7.15/0.82/29.10 \\
\hline
Set14(c) 	& 20.78/0.22/19.67 & 7.73/0.85/28.26 & 6.59/0.90/29.65 & 6.21/0.90/30.16 \\
BSDS100(c) 	& 20.53/0.17/19.82 & 6.16/0.90/30.27 & 6.53/0.90/29.77 & 5.02/0.92/32.06 \\
Urban100(c) & 20.04/0.29/20.35 & 7.59/0.83/28.37 & 7.62/0.82/28.29 & 5.76/0.88/30.82\\  
DIV2K(c) 	& 20.30/0.12/19.74 & 6.06/0.86/30.66 & 6.19/0.85/30.58 & 5.25/0.89/31.95\\ 
SIDD(c) 	& 20.71/0.11/19.33 & 4.00/0.92/34.90 & 4.08/0.92/34.90 & 3.03/0.94/37.34\\
\hline
\multicolumn{5}{|c|}{Initial noised image with 40\% of noise}\\
\hline
Set12(g) 	& 40.62/0.18/13.47 & 12.34/0.72/23.81 & 12.61/0.71/23.63 & 10.30/0.72/25.38 \\
BSD68(g) 	& 40.27/0.16/14.08 & 12.65/0.71/24.13 & 12.80/0.70/24.04 & 10.77/0.73/25.54 \\
\hline
Set14(c) 	& 40.03/0.09/13.98 & 10.83/0.82/25.34 & 10.79/0.83/25.36 & 10.38/0.81/25.70 \\
BSDS100(c) 	& 40.26/0.06/13.97 & 11.50/0.85/24.85 & 11.74/0.85/24.67 & 10.25/0.85/25.85 \\
Urban100(c) & 40.44/0.11/14.25 & 14.35/0.75/22.81 & 14.43/0.74/22.59 & 12.45/0.78/24.03\\
DIV2K(c) 	& 40.32/0.04/13.78 & 12.41/0.80/24.36 & 12.58/0.79/24.21 & 10.89/0.82/25.71\\
SIDD(c) 	& 40.21/0.03/13.57 &  9.05/0.90/27.88 & 9.17/0.90/26.93 & 8.38/0.89/28.55 \\ 
\hline
\end{tabular}
\caption{\rev{Test 4. Numerical experiments on classic and high-resolution datasets}}
\label{table:t4-ds}
\end{table}

\rev{
Finally, we present results for classic and high-dimensional datasets. 
In Table \ref{table:t4-ds}, we present numerical results for two greyscale image datasets (Set12 and BSD68) and five color image datasets (Set14, BSDS100, Urban100, DIV2K, SIDD). We present the average RRMSE, SSIM, and PSNR results of different denoising methods for 20 \% and 40 \% of noise. Similarly, we fix local domain size $N_h^{K_i} = 16^2$ with varying coarse grid size depending on the given image resolution in each dataset. The result shows that the proposed method performed well for image-denoising tasks. It has some open problems related to the further optimization of the computational cost, the development of various choices of the multiscale space representation, the interplay between the cost of local calculations and global solution, and adaptive techniques in choosing time step size and nonlinear diffusion parameters. It opens many interesting questions about the intersection of optimization method development and acceleration of the iterations with techniques from a solution of PDEs with implicit time approximation, high-order time approximation schemes, decoupling techniques, and iterative methods with preconditioning. Furthermore, it would be interesting to consider hybrid approaches that combine various denoising techniques to obtain a better denoising performance and computationally efficient algorithms, especially with multiple types of machine learning techniques, including supervised, unsupervised, and semi-supervised learning. 
}

\section{Conclusion}

We considered a time-dependent nonlinear diffusion equation with application to the image-denoising process. In this approach, the noised image is given as the initial condition for the time-dependent problem, and time iterations with a nonlinear diffusion operation are performed to produce a denoised image at the final time. 
The traditional approximation scheme was constructed for a given image resolution using a finite volume approximation \rev{with explicit} and implicit time stepping and a linearization from the previous time layer. 
In order to construct a computationally efficient and accurate algorithm, we proposed a novel multiresolution approach based on the Multiscale Finite Element Method. In the proposed multiscale method, we construct basis functions for a coarse image resolution based on solving \rev{local} problems in each local domain of basis support. In the basis construction, we showed that the original noised image can lead to poor approximation properties that lead to the introduction of local denoising to preprocess the given local image and produce good coefficients for \rev{basis construction}. 
The resulting multiscale basis functions were used to construct an accurate projection operator and perform a denoising process with fewer iterations.
The presented numerical results show the applicability of the multiscale method for the image-denoising process. We considered \rev{classic and high-resolution image datasets}, applied varying noise, and performed image denoising using the proposed techniques. We showed that the method accurately approximates high-resolution images, efficiently reduces a denoised image in fewer number of iteration, and performs faster calculations  within reduced resolution. 
In future works, we will consider the construction of a hybrid framework for PDE-based image denoising by incorporating machine learning techniques for the fast construction of a coarse-scale model. Optimization techniques should be also incorporated for parameter tuning for better noise removal. Moreover, we will extend multiscale techniques to fourth-order PDEs that give better edge preservation properties. 


\bibliographystyle{unsrt}
\bibliography{lit}

\end{document}